\documentclass{article}
\usepackage{latexsym,amssymb,amsmath,amsthm,verbatim}
\usepackage{float,epsfig}
\usepackage{fig}
\def\begfig {
\begin{figure}
\small
}
\def\endfig {
\normalsize
\end{figure}
}

\newcommand{\path}{./figures/}

\def\FigHOneRhombus {
 \begfig
 \hspace{0.6in} 
 \epsfxsize = 3.3in 
 \epsffile{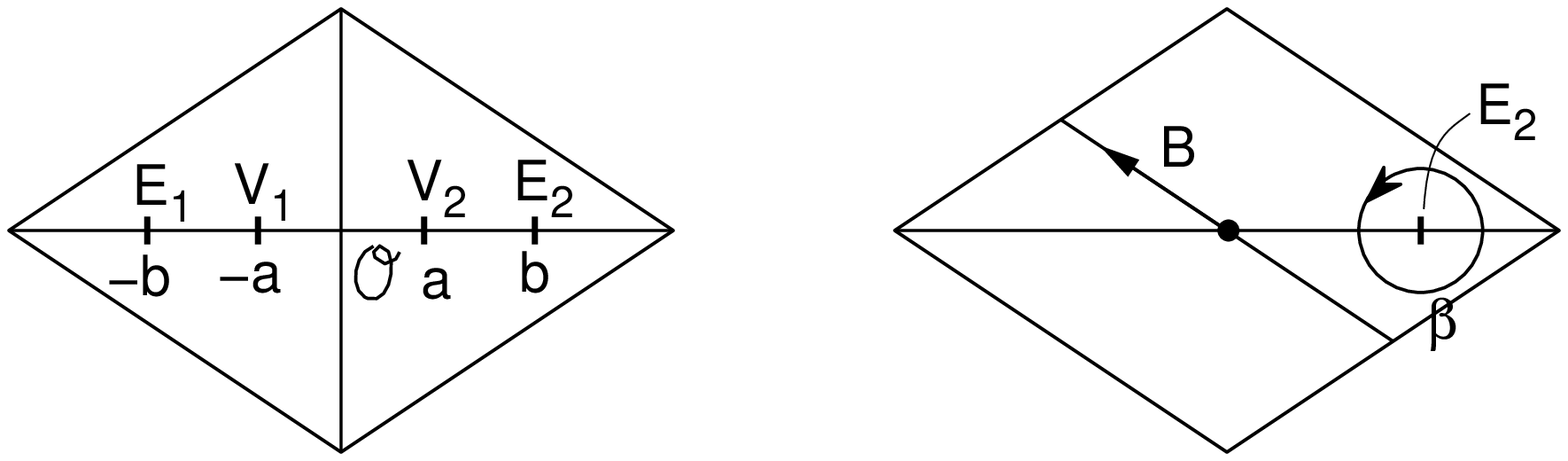} 
 \hfill 
 \begin{center} 
 \vspace{0.3cm}
 \parbox{4.0in}{
  \caption{
   \label{figure:honerhombus}
   The rhombus model for the underlying Riemann surface of 
   the quotient surface of \Hone/ modulo translations is drawn on the 
   left. The vertical points, $V_i$ and the horizontal points,
   $E_i$, on the
   horizontal diagonal are included. The cycles $\beta$ and $B$,
   illustrated on the right, and their reflections in
   the vertical diagonal, generate a homology basis for the
   punctured torus. 
  }
 }
 \end{center} 
 \endfig
}

\def\FigslitTone {
 \begfig
 \hspace{0.8in} 
 \epsfxsize = 3.0in 
 \epsffile{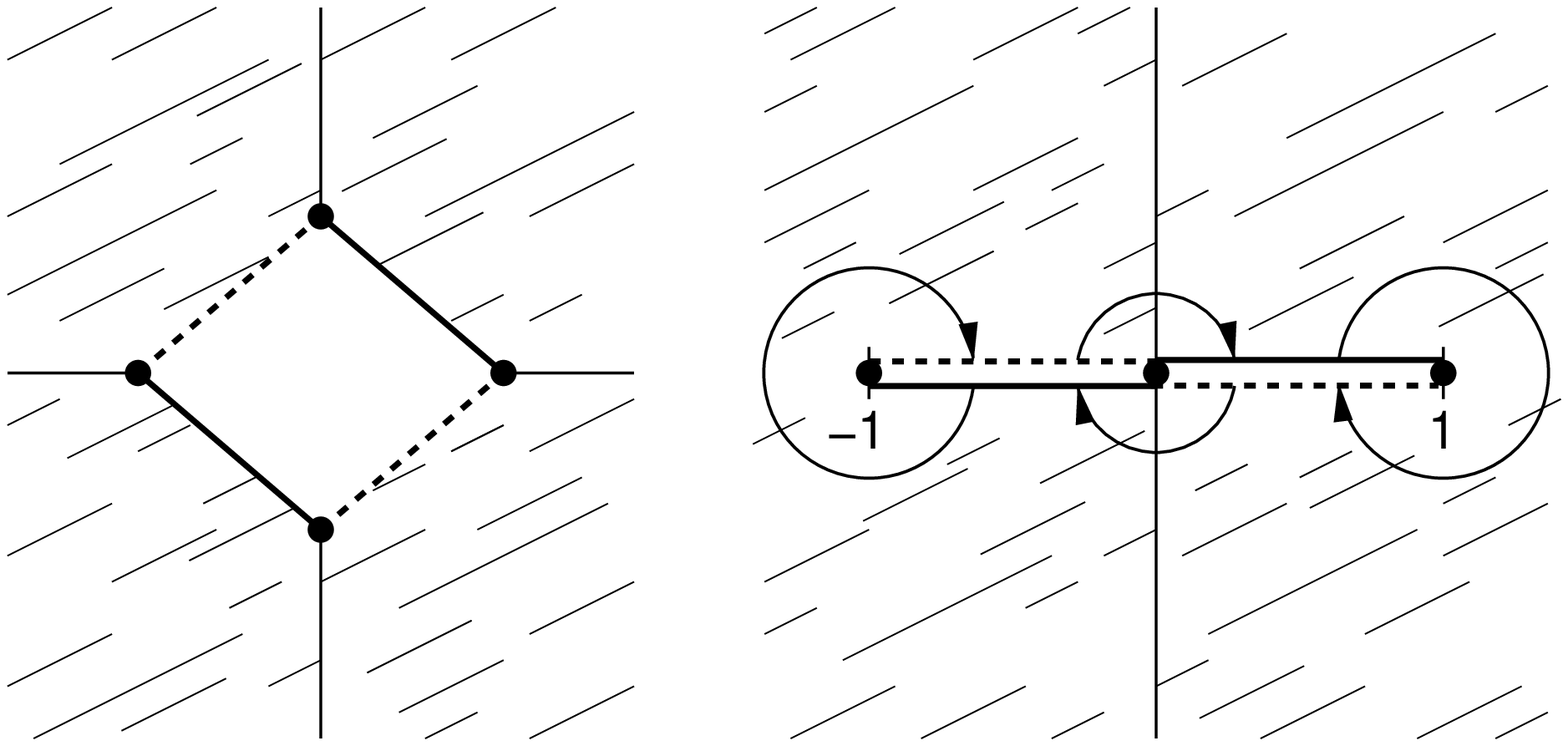} 
 \hfill 
 \begin{center} 
 \vspace{0.3cm}
 \parbox{4.0in}{
  \caption{
   \label{figure:slitTone}
   On the left: a rhombus removed from the extended complex plane.
   Identification of the opposite sides of the remaining set produces a 
   Riemann
   surface of genus one that is also a rhombic torus. On the right: a
   degenerate rhombus is removed from the extended plane. Identification 
   of opposite sides produces a non-degenerate rhombic torus.
   The one-form $dz$ on the extended complex plane has a double pole at
   infinity. On the torus we have constructed,  the one-form $d\zeta$ 
   descends to define a one-form 
   with a
   double pole at the point corresponding to $\infty$ and a double zero 
   at the
   vertex point of the removed rhombus. The periods of this induced form 
   on the edge cycles of the excised rhombus are
   given by the  complex numbers defined by the edges themselves.
   In the case of the torus on the 
   right,
   these periods are real.
  }
 }
 \end{center} 
 
 \endfig
}

\def\Figslitandrhombus {
 \begfig
 \hspace{0.6in} 
 \epsfxsize = 3.0in 
 \epsffile{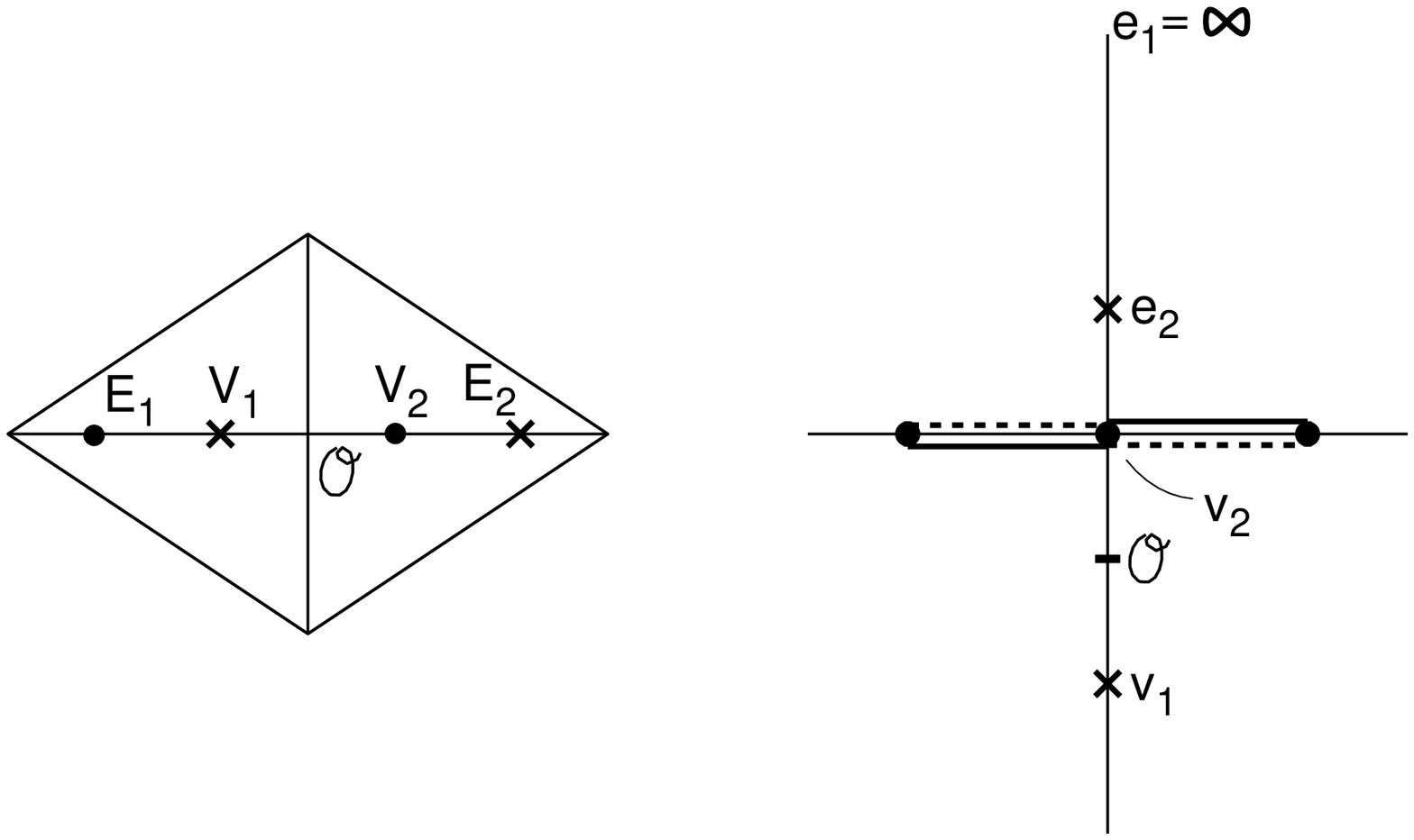} 
 \hfill 
 \begin{center} 
 \vspace{0.3cm}
 \parbox{4.0in}{
  \caption{
   \label{figure:slitandrhombus}
   The slit model for the  rhombic torus has two distinguished points, $e_1$ 
   and
   $v_2$. We may identify the slit model with a rhombus so that the 
   imaginary
   axis is identified with the horizontal diagonal. Labelling the points 
   on the
   rhombus corresponding to $e_1$ and $v_2$ in the same manner, but with 
upper-case
   letters, we may
   translate the rhombus horizontally to place the center, \O/, of the 
   rhombus
   to the right of $E_1$ and to the left of $E_2$. The points $V_1$ and 
   $E_2$
   are defined to be the points on the horizontal diagonal symmetric
   (with respect to \O/) to $V_2$ and $E_1$, respectively.  
  }
 }
 \end{center} 
 \endfig
}

\def\Fig_a_and_b {
 \begfig
 \hspace{1.2in} 
 \epsfxsize = 1.9in 
 \epsffile{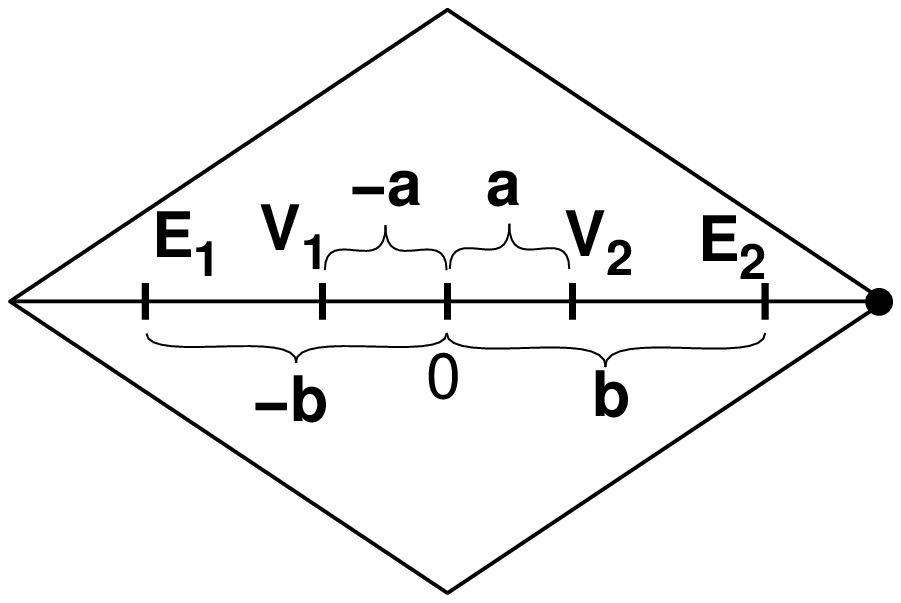} 
 \hfill 
 \begin{center} 
 \vspace{0.3cm}
 \parbox{4.0in}{
  \caption{
   \label{figure:_a_and_b}
   We label the points $V_i$ and $E_i$ according to their signed and
   scaled distance from \O/. From Abel's Theorem we know that $a+b=1$.
  }
 }
 \end{center} 
 \endfig
}

\def\FigRectangle_R {
 \begfig
 \hspace{0.8in} 
 \epsfxsize = 3.2in 
 \epsffile{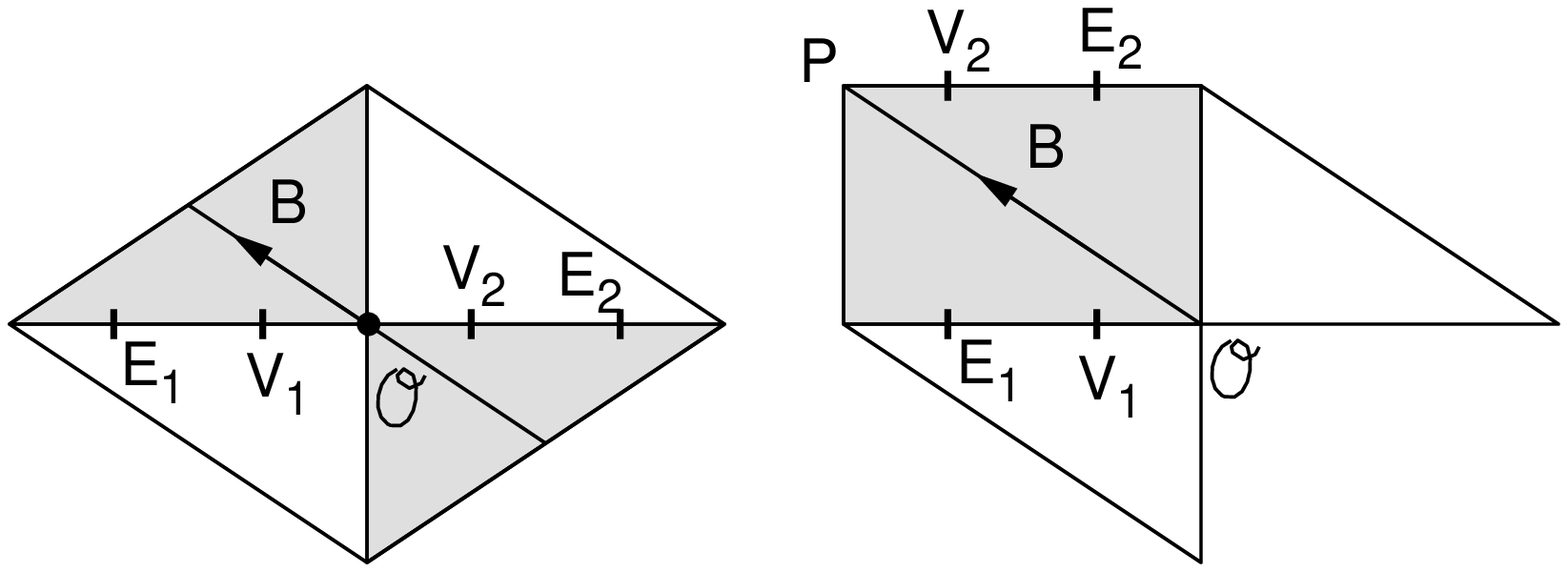} 
 \hfill 
 \begin{center} 
 \vspace{0.3cm}
 \parbox{4.0in}{
  \caption{
   \label{figure:Rectangle_R}
   On the left, the path, $B$, is illustrated on \T/. 
   On the right, an equivalent path, also labelled
   $B$, is drawn on the shaded rectangle \aR/.
  }
 }
 \end{center} 
 \endfig
}

\def\FigDhonRone {
 \begfig
 \hspace{0.5in} 
 \epsfxsize = 3.0in 
 \epsffile{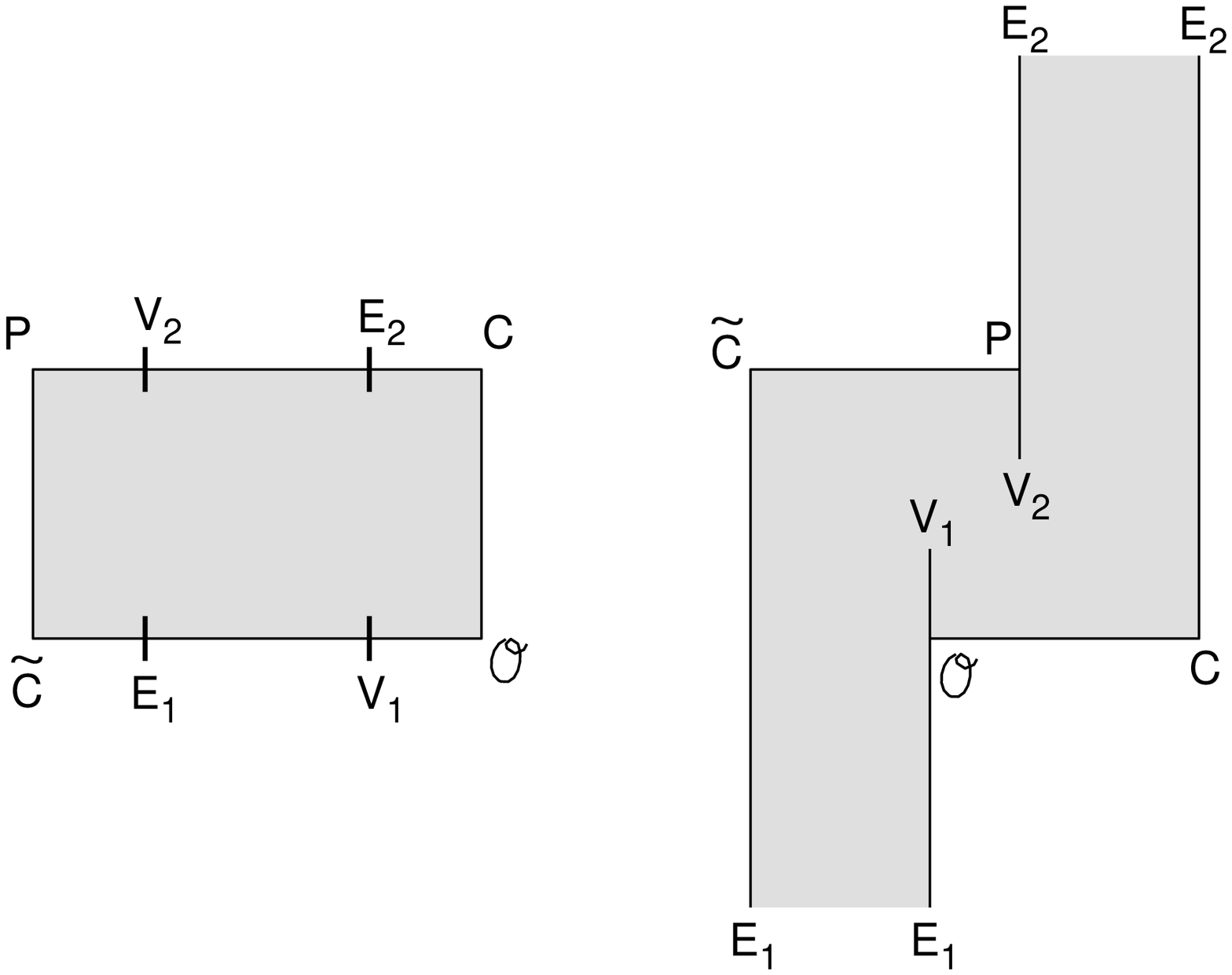} 
 \hfill 
 \begin{center} 
 \vspace{0.3cm}
 \parbox{4.0in}{

  \caption{
   \label{figure:DhonRone}
   On the left is the rectangle \aR/. The points $V_i$ and $E_i$ are  
   placed 
   for a generic value of $a$ between $0$ and $1/2$. On the right is the 
   image
   of \aR/ under the mapping $F_a$.
  }
 }
 \end{center} 
 \endfig
}

\def\FigDhonRtwo{
 \begfig
 \hspace{0.4in} 
 \epsfxsize = 3.3in 
 \epsffile{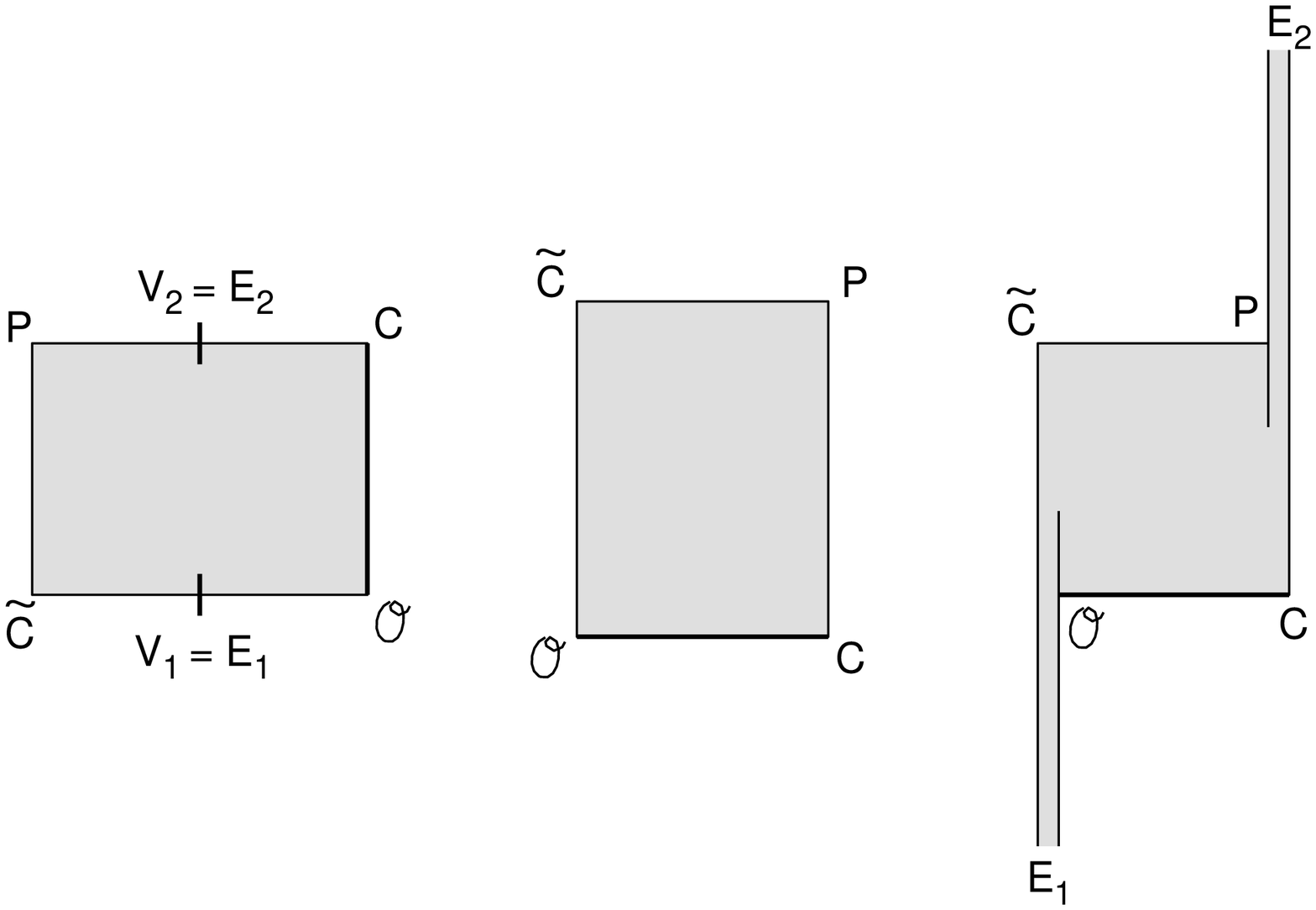} 
 \hfill 
 \begin{center} 
 \vspace{0.3cm}
 \parbox{4.0in}{
  \caption{
   \label{figure:DhonRtwo}
   The case $a=\f12$.  When $a=\f12$, 
   we have $E_i=V_i$ as is illustrated in the 
   picture of \aR/ on the left. In this case, 
   the one-form $dh$ is regular and the image of 
   \aR/ 
   under $F_{\f12}$ is illustrated in the center image. Up to scaling, 
   $F_{\f12}$ is
   clockwise rotation by $90$ degrees. The image of \aR/ under $F_a$ for 
   $a$
   near $\f12$ is illustrated on the right.
  }
 }
 \end{center} 
 \endfig
}

\def\FigDhonRthree {
 \begfig
 \hspace{0.2in} 
 \epsfxsize = 4.4in 
 \epsffile{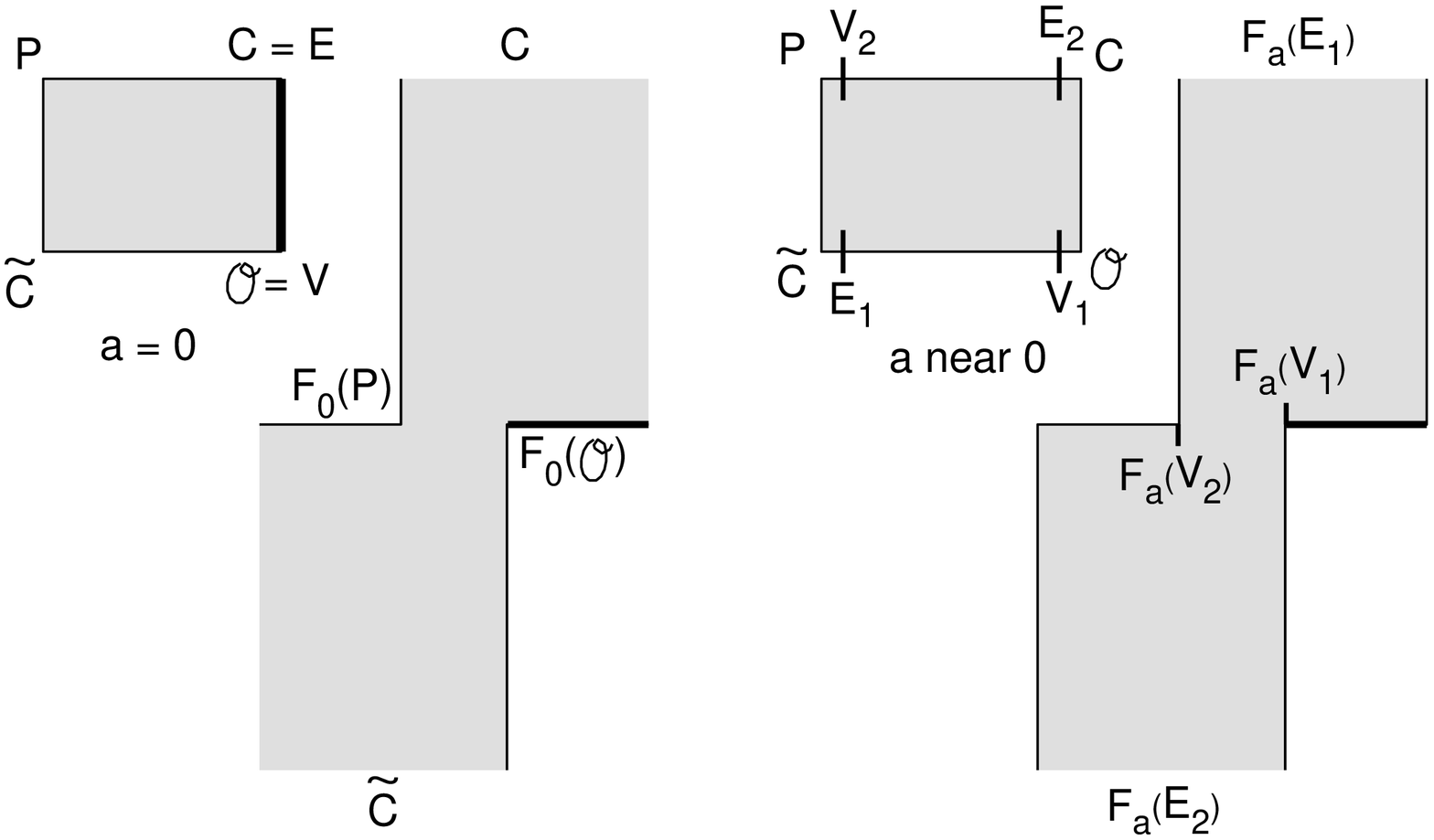} 
 \hfill 
 \begin{center} 
 \vspace{0.3cm}
 \parbox{4.0in}{
  \caption{
   \label{figure:DhonRthree}
   On the extreme left, the rectangle ${\cal R}$  is drawn with $V_1=V_2=$\O/ 
   and $E_1=E_2$ at $C$ and $\tilde{C}$. The image of ${\cal R}$ under 
   $F_0$ is 
   drawn just to the right, illustrating that $F_0($\O/$)$ lies to the 
   right of
   $F_0(P)$ when $a=0$. On the right-hand-side, ${\cal R}$  and $F_a(\aR/)$ 
   are drawn for $a\neq 0$ small.
  }
 }
 \end{center} 
 \endfig
}

\def\FigTjigsaw {
 \begfig
 \hspace{0.8in} 
 \epsfxsize = 3.0in 
 \epsffile{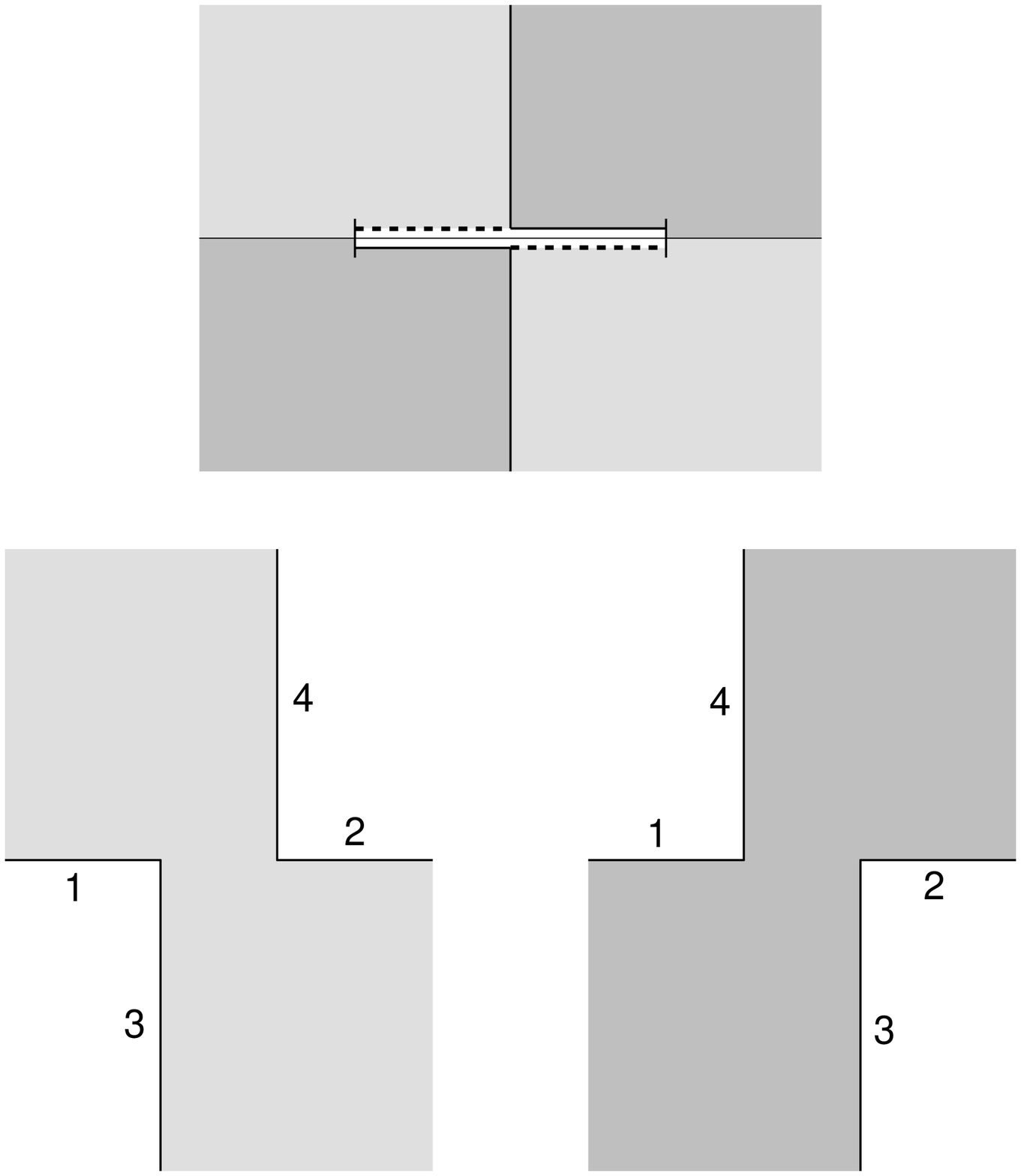} 
 \hfill 
 \begin{center} 
 \vspace{0.3cm}
 \parbox{4.0in}{
  \caption{
   \label{figure:Tjigsaw}
   \aR/  in  the slit model. The slit model of \T/ (in 
   Figure~\ref{figure:slitandrhombus}) can be cut along symmetry lines---the 
   real and imaginary axes---to produce two conformal rectangles. In the 
   illustration, 
   they are the two differently shaded regions. As described in the text, 
   one of them, the one on the right, must be (up to scaling) the image of 
   \aR/ 
   under $F_0$.
  } 
 }
 \end{center} 

 \endfig
}

\def\FigHelicoidZeroOne {
 \begfig
 \hspace{1.0in} 
 \epsfxsize = 0.9in 
 \epsffile{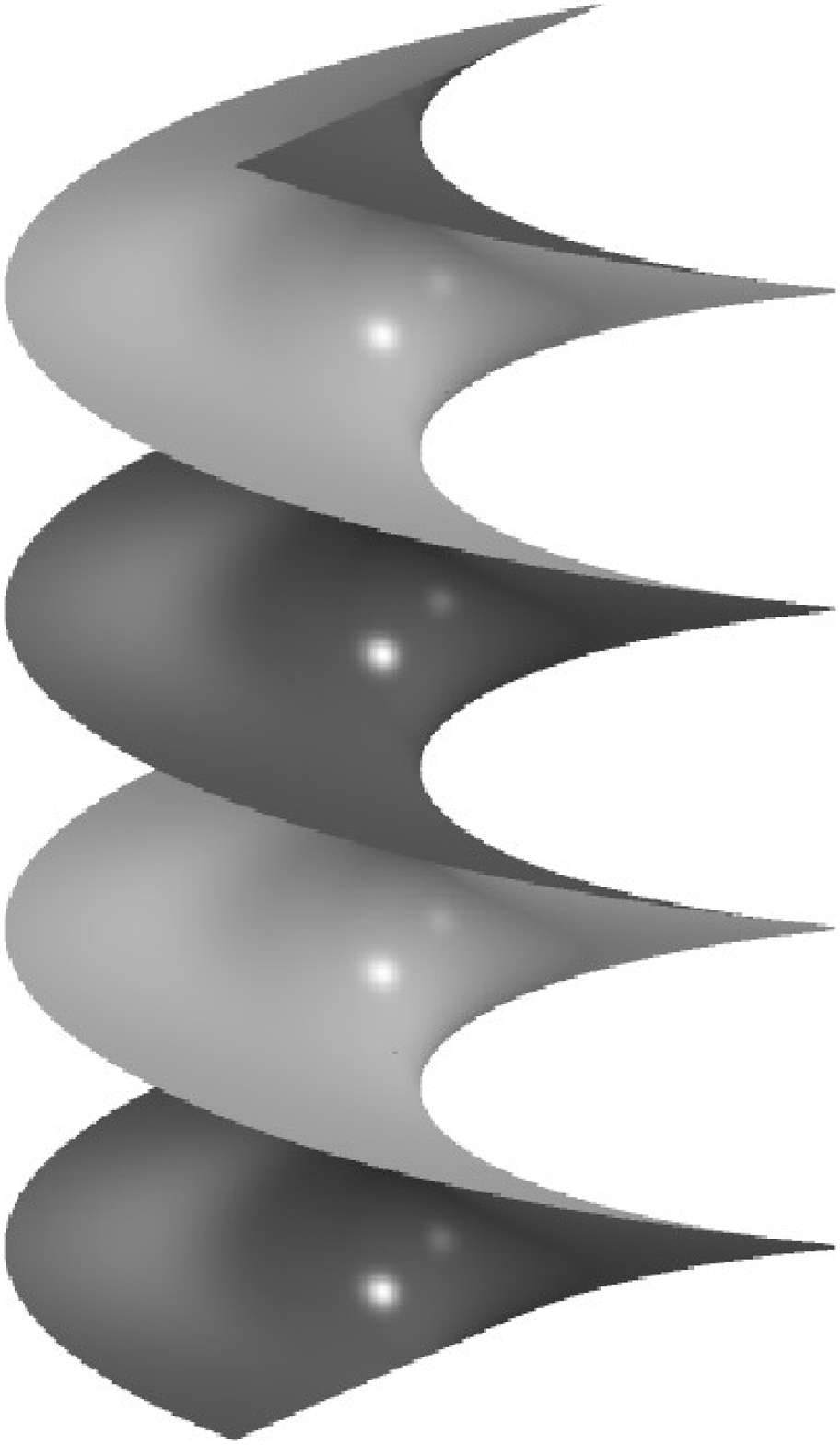}
 \hspace{0.2in} 
 \epsfxsize = 0.9in 
 \epsffile{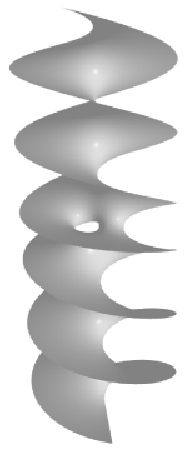}
 \hfill
 \begin{center} 
 \vspace{0.3cm}
 \parbox{4.0in}{
  \caption{
   \label{figure:helicoidzeroone}
   The helicoid and the \HEone/ of \cite{howe3}}
  }
 \end{center} 
 \endfig
}

\def\FigPerihel {
 \begfig
 \hspace{0.2in} 
 \epsfxsize = 1.6in 
 \epsffile{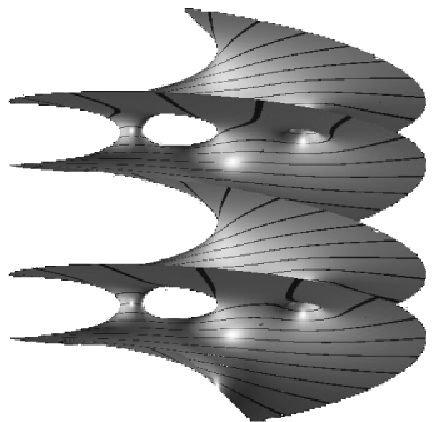}
 \hspace{0.2in} 
 \epsfxsize = 1.6in 
 \epsffile{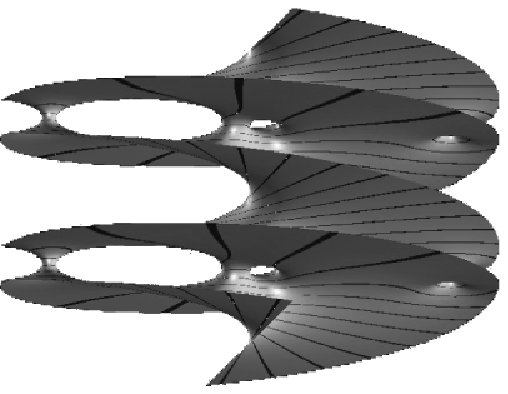}
 \hfill
 \begin{center} 
 \vspace{0.3cm}
 \parbox{4.0in}{
  \caption {
   \label{figure:perihel}
   Images of higher genus examples of translation-invariant surfaces 
   with genus two (left) and three (right) in the quotient.   They satisfy 
the other  geometric
   conditions of (4). These images were computed by Martin Traizet using 
MESH . There is at the time of writing no proof that they exist. See 
   however,  the recent work of Traizet and Weber \cite{TraizetWeber03}. The 
software package MESH is described in \cite{chh1}. A Java version and the 
manual are 
available online at {\small 
www.msri.org/publications/sgp/jim/software/jmesh/dist/indexc.html.} An 
archive
of computed minimal surfaces including  the \Hk/ and \HEone/ can be found at 
{\small \dots /jmesh/surfacearchive/indexc.html} .
  }
 }
 \end{center} 
 \endfig
}

\def\FigG2Heli {  
 \begfig
 \hspace{1.4in} 
 \epsfxsize = 1.1in 
 \epsffile{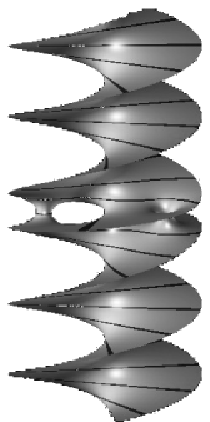}
 \hfill
 \begin{center} 
 \vspace{0.3cm}
 \parbox{4.0in}{
  \caption{
   \label{figure:g2heli}
   A computed image of a  conjectured 
   genus-two analog of \HEone/. 
   It satisfies all the other conditions of Theorem~1, also computed by 
Traizet 
   using MESH.
  }
 }
 \end{center} 
 \endfig
}

\def\FigHelPOne {
 \begfig
 \hspace{1.4in} 
 \epsfxsize = 1.3in 
 \epsffile{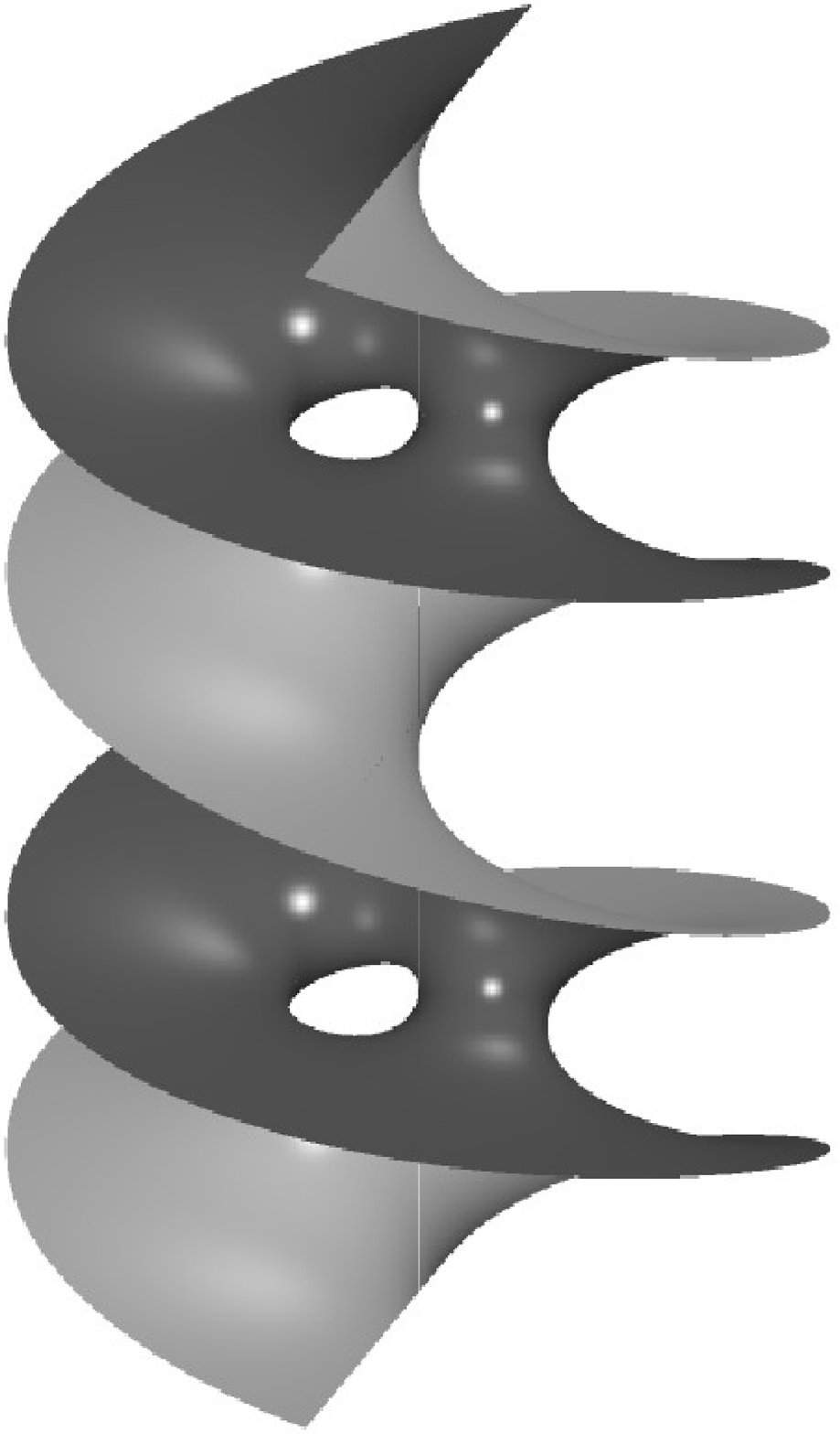}
 \hfill
 \begin{center} 
 \vspace{0.3cm}
 \parbox{4.0in}{
  \caption{
   \label{figure:help1}
   Singly Periodic Genus-one helicoid, \Hone/.}
 }
 \end{center} 
 \endfig
}

\def\FigFundDomainhK {
 \begfig
 \begin{center} 
 \vspace{0.3cm}
 \hspace{0.8in}
 \epsfxsize = 3.2in
 \epsffile{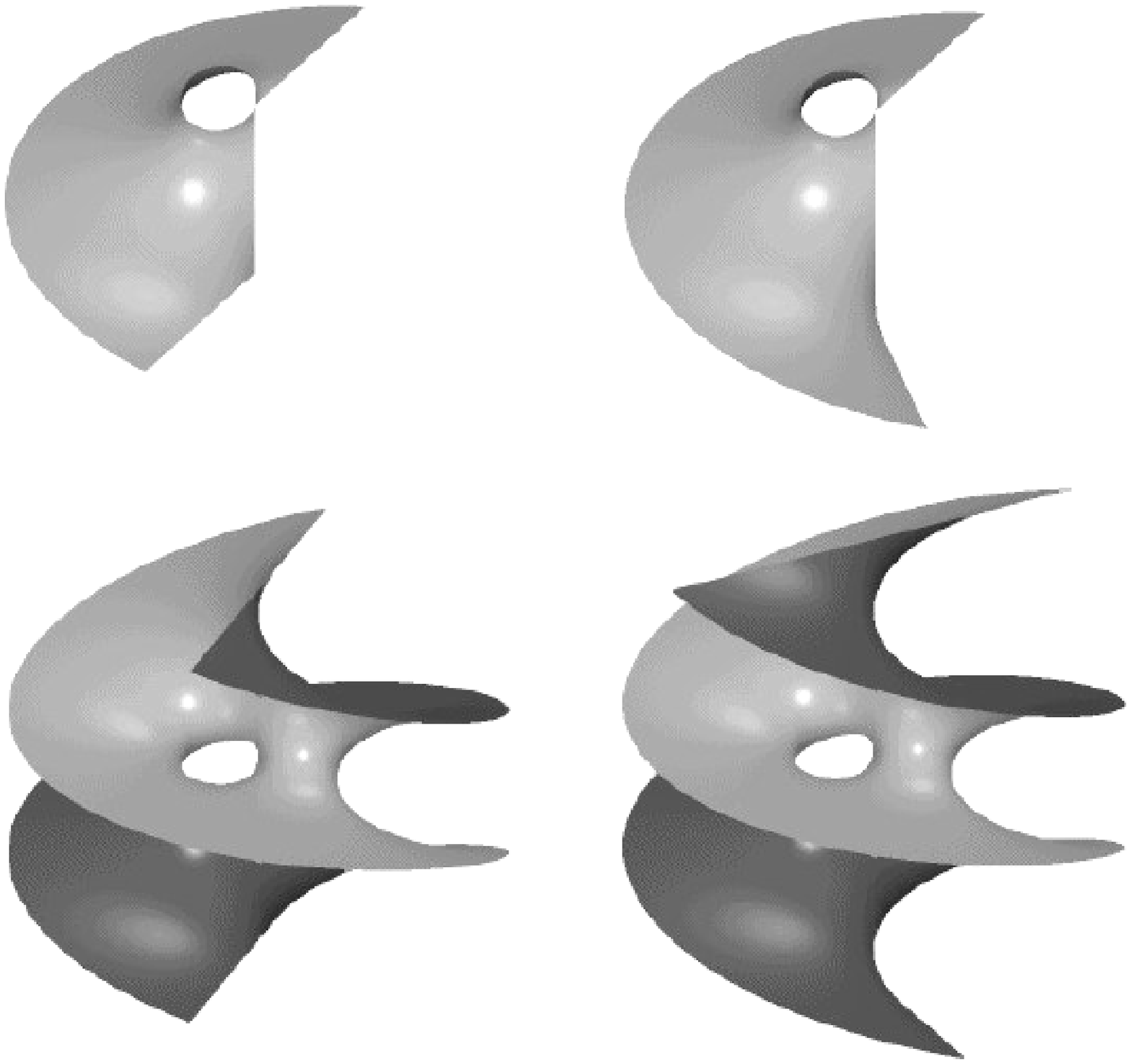}
 \hfill
 \parbox{4.0in}{
  \caption{
   \label{figure:funddomainhk}
   A fundamental domain of \Hk/ can be imagined
   qualitatively as a region of the helicoid bounded by two
   horizontal lines between which the helicoid turns by an angle
   of $\theta=2\pi k$ and in the middle of which there is a
   handle. The boundary lines are identified in the quotient as
   a single line. The only other line, besides the vertical
   axis, that survives the surgery necessary to insert the
   handle is a horizontal line in the middle, at the level of
   the handle. The surface on the bottom, left, is a fundamental domain of 
   \Hone/. The surface on the bottom, right, is a fundamental domain of \Hk/ 
   for $k\sim 1.25$. 
   The two images on top are are each one quarter of the surfaces below them. 
   They are bounded by a segment of the vertical axis, two horizontal half 
   lines and a closed loop that that is not contractible in \Hk/$/\sigma_k$. 
   See also \cite{howe5}.
  }
 }
 \end{center} 
 \endfig
}

\def\FigSPGONEhalf {
 \begfig
 \hspace{0.8in}
 \epsfxsize = 1.8in
 \epsffile{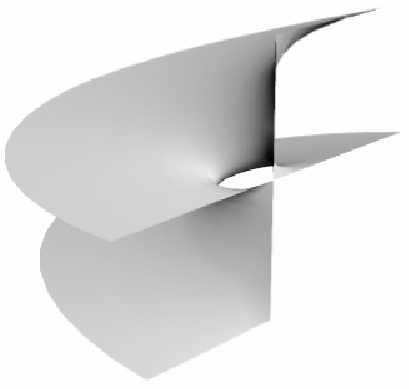}
 \hfill
 \begin{center}
 \vspace{0.3cm}
 \parbox{5.0in}{
  \caption{
   \label{figure:spgonehalf}
   A fundamental domain of \Hone/ can be imagined
   qualitatively as a region of the helicoid bounded by two
   horizontal lines between which the helicoid turns by an angle
   of $2\pi$ and in the middle of which there is a
   handle.  Illustrated here is one-half of the fundamental domain;
   this half is on the side of the vertical
   plane that contains the vertical axis and the two
   horizontal lines. The other half of the fundamental domain is
   produced by reflection through the vertical axis.
   The boundary lines, top and bottom, are identified in the quotient as
   a single line. The only other line, besides the vertical
   axis, that survives the surgery necessary to insert the
   handle is a horizontal line in the middle, at the level of
   the handle. 
   The normal symmetry line $L$ is not illustrated here but
   can be easily visualized
   as the horizontal line perpendicular to the vertical and 
   horizontal line through the center point. 
(See (\cite{howe5}).)
  }
 }
 \end{center}
 \endfig
}

\def\FigTorus1 {
 \begfig
 \hspace{1.0in} 
 \epsfxsize = 1.9in 
 \epsffile{\path torus1.ps}
 \hfill
 \begin{center} 
 \vspace{0.3cm}
 \parbox{4.0in}{
  \caption{
   \label{figure:torus1}
   Diagram of the lattice
  }
 }
 \end{center} 
 \endfig
} 

\def\FigDivisorsZero {
 \begfig
 \begin{center}
 $
 \begin{array}{c|lll}
            & z=0        & z=\infty &       \\
 \hline
 dh         & \infty    &   \infty    &   \\
 g          & 0^k      & \infty^k   &       \\
 g dh       & 0^{k-1} & \infty^{k+1}      &     \\
 \frac1g dh &    \infty^{k+1}        & 0^{k-1}   & \\
 \end{array}
 $
 \vspace{0.3cm}
 \parbox{4.0in}{
  \caption{
   \label{figure:divisors0}
   The divisors of $g$, $dh$, $gdh$ and $\frac1g dh$ for the heliciod   
modulo $\sigma_k$, a vertical screw motion with twist angle $2\pi k$.
  }
 }
 \end{center} 
 \endfig
 }
\def\FigDivisorsOne {
 \begfig
 \begin{center}
 $
 \begin{array}{c|llll}
            & E_1        & V_1    & V_2 & E_2    \\
 \hline
 dh         & \infty     & 0      & 0   & \infty \\
 g          & \infty     & \infty & 0   & 0      \\
 g dh       & \infty^{2} & *      & 0^2 & *      \\
 \frac1g dh & *          & 0^2    & *   & \infty^{2} \\
 \end{array}
 $
 \vspace{0.3cm}
 \parbox{4.0in}{
  \caption{
   \label{figure:divisors1}
   The divisors of $g$, $dh$, $gdh$ and $\frac1g dh$.
  }
 }
 \end{center} 
 \endfig
}

\def\FigDivisorsThree {
 \begfig
 \begin{center} 
 $ 
 \begin{array}{cllll}
          &E_1          & V_1    & V_2 & E_2 \\
 \\
 g        & \infty^k    &\infty  & 0   & 0^k \\
 dh       & \infty      &0       & 0   & \infty.
 \end{array}
 $ 
 \parbox{4.0in}{
  \caption{
   \label{figure:divisors3}
   The divisors of $g$ and $dh$.
  }
 }
 \end{center} 
 \endfig
}

\def\FigDivisorsFive {
 \begfig
 \begin{center} 
 $ 
 \begin{array}{c|llll}
            & E_1          & V_1    & V_2    & E_2          \\
 \hline
 g          & \infty^k     & \infty & 0      & 0^k          \\
 dh         & \infty       & 0      & 0      & \infty       \\
 g dh       & \infty^{k+1} & -      & 0^2    & 0^{k-1}      \\
 \frac1g dh & 0^{k-1}      & 0^2    & -      & \infty^{k+1} \\
 \frac{dg}g  & \infty      & \infty & \infty & \infty       \\
 \text{residue}\frac{dg}g  &   k    & 1      &-1    &-k     \\
 \end{array}
 $ 
 \parbox{4.0in}{ 
  \caption{
   \label{figure:divisors5}
   Divisor and residue requirements for the Weierstrass data 
   of \Hk/.
  }
 }
 \end{center} 
 \endfig
}

\def\FigDivisorsSix {
 \begfig
 \begin{center} 
 $ 
 \begin{array}{c|llll}
          & E_1    & V_1 & V_2 & E_2       \\
 \hline
 dh       & \infty & 0   & 0   & \infty    \\
 \end{array}
 $ 
 \parbox{4.0in}{
  \caption{
   \label{figure:divisors6}
   The divisor of $dh$.
  }
 }
 \end{center} 
 \endfig
}

\def\FigTkd {
 \begfig
 \hspace{0.8in}
 \epsfxsize = 2.5in
 \epsffile{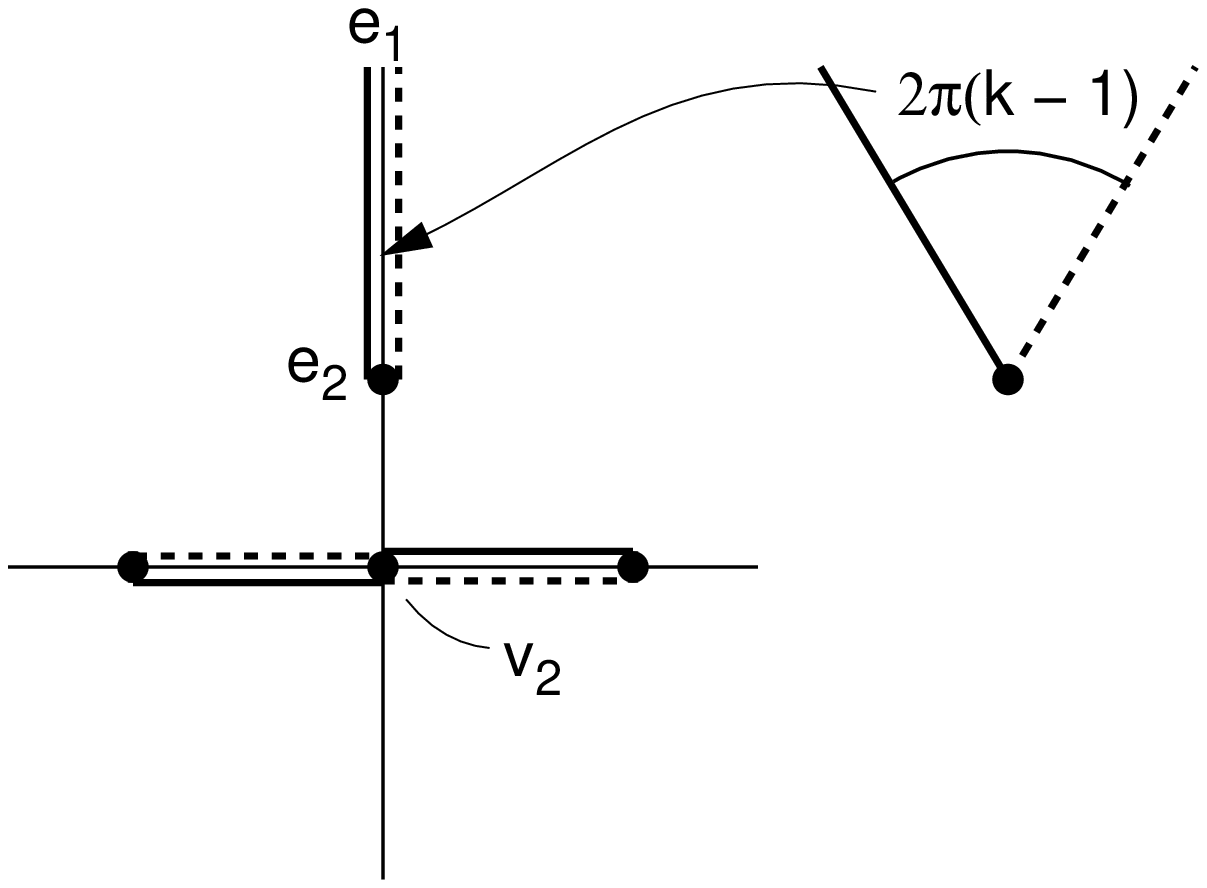}
 \hfill 
 \begin{center} 
 \vspace{0.3cm}
 \parbox{4.0in}{
  \caption{
   \label{figure:tkd}
   The torus \tkd/ is constructed by sewing into the torus \T/
   the cone metric $S_{k-1}$ described in Section~2. The 
   sewing is done along the ray in the imaginary axis beginning at $di$ 
   and terminating at $\infty$. The point $di$ is labeled $e_2$, the point
   at $\infty$ is labeled $e_1$, and the origin is labeled $v_2$. 
  Note that when $k=1$ we sew nothing in: 
   ${\cal T}_1(d)=$\T/ with a distinguished point at
   $\zeta=di$.  We may also remove 
   an $S_k$ from \T/, for $0<k<\f12$, allowing us to define \tkd/ for 
   $k>\f12$.
  }
 }
 \end{center} 
 \endfig
}

\def\FigKdrectangle {
 \begfig
 \hspace{1.1in}
 \epsfxsize = 2.2in
 \epsffile{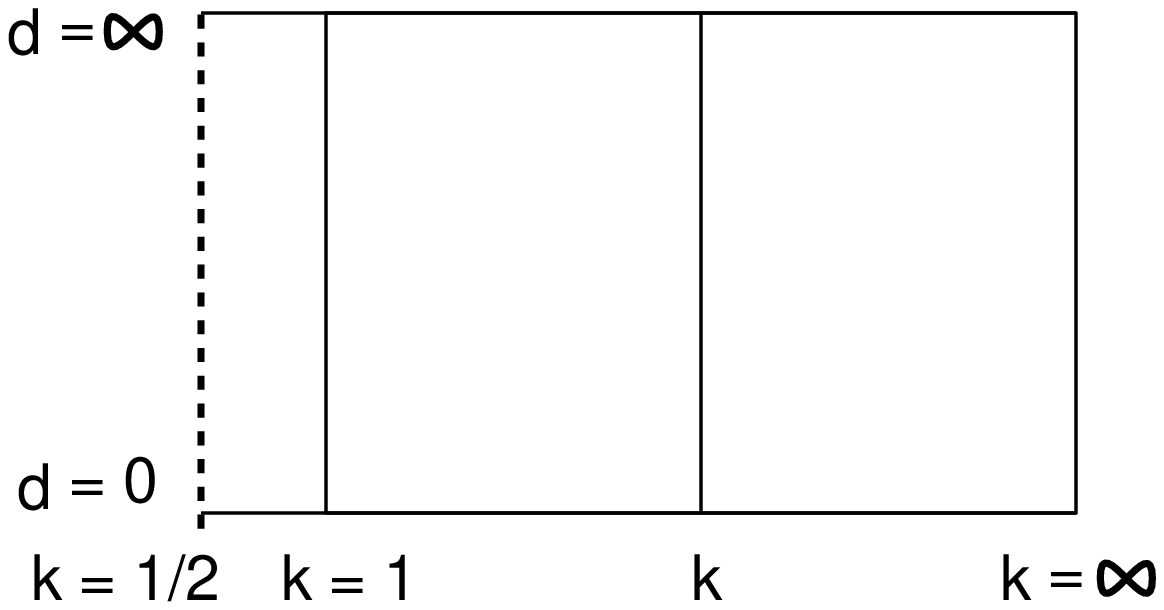}
 \hfill 
 \begin{center} 
 \vspace{0.3cm}
 \parbox{4.0in}{
  \caption{
   \label{figure:kdrectangle}
   The $(k,d)$ rectangle. Each point corresponds to a \tkd/, on which we 
   have defined $gdh$, $\f1g dh$, $g$ and $dh$.
  }
 }
 \end{center} 
 \endfig
}

\def\FigKdrectangletwo {
 \begfig
 \hspace{1.1in}
 \epsfxsize = 2.2in
 \epsffile{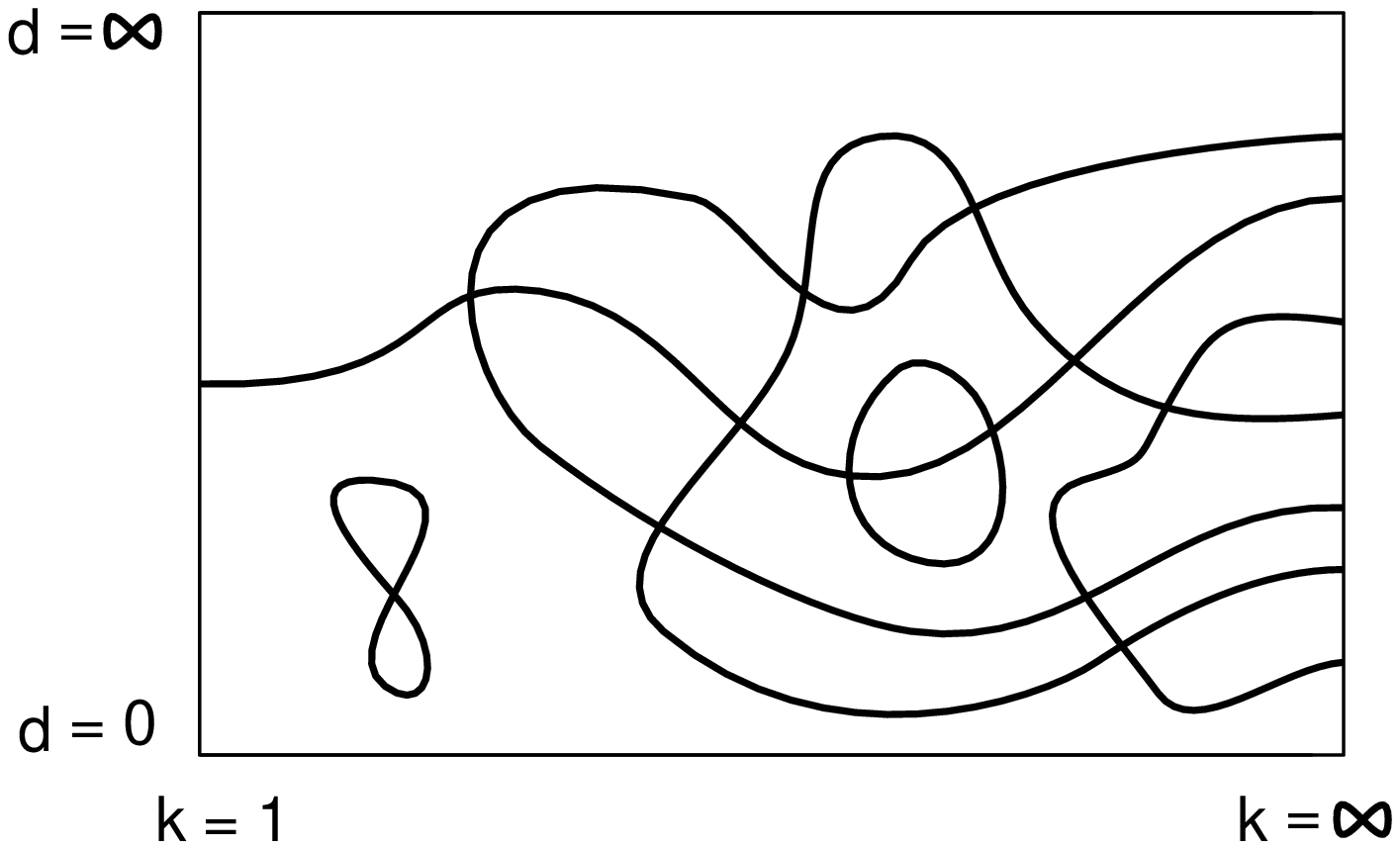}
 \hfill 
 \begin{center} 
 \vspace{0.3cm}
 \parbox{4.0in}{
  \caption{
   \label{figure:kdrectangle2}
   On the $(k,d)$ rectangle, we have defined a function $h(k,d)$ whose 
   value is
   the real part of the integral of $dh$ on the path $B$ of the underlying 
   \tkd/. 
   Illustrated here is a possible locus of the zeros of $h(k,d)$. We know 
   that the locus
   is an analytic set and that it has a unique point on the vertical 
   line
   $k=1$ (See Proposition~\ref{rmk2} in Section~3.3.) 
   corresponding to \Hone/. 
   We show that there is at 
   least one path connecting  this point to a point on the right-hand side of 
   the rectangle, where $k=\infty$. Note: the illustration in this 
   figure is of the closed rectangle $[1,\infty] \times [0,\infty]$
   which is properly contained in $\SP$. 
  }
 }
 \end{center} 
 \endfig
}

\def\FigGammaPaths {
 \begfig
 \hspace{0.6in}
 \epsfxsize = 3.4in
 \epsffile{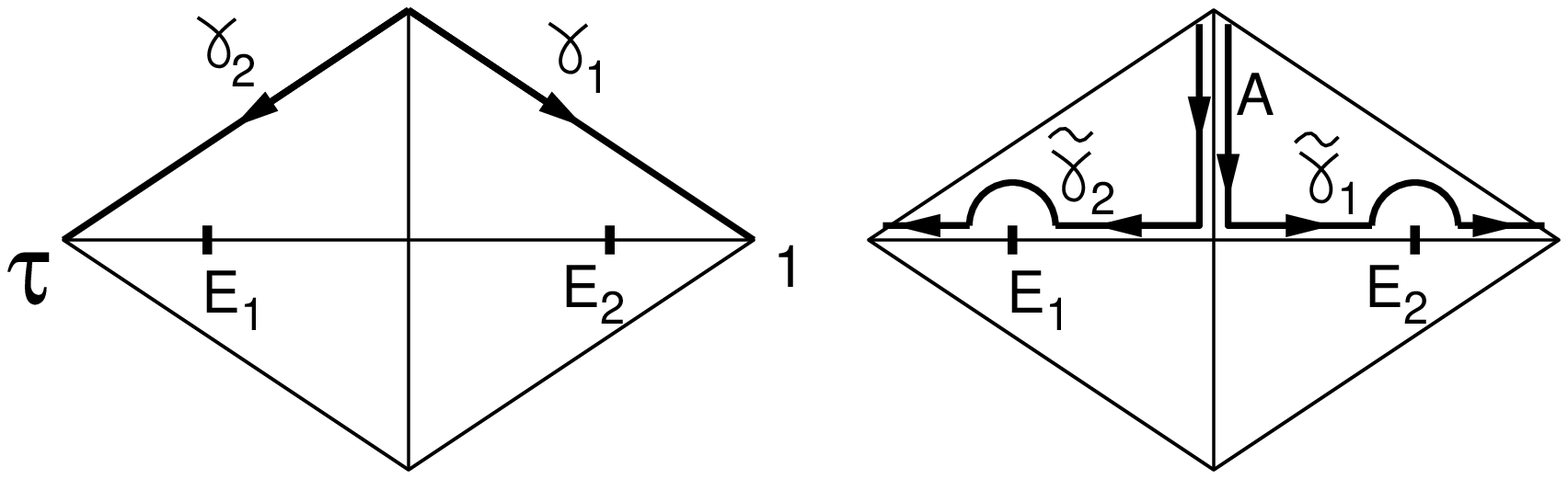}
 \hfill 
 \begin{center} 
 \vspace{0.3cm}
 \parbox{4.0in}{
  \caption{
   \label{figure:gammapaths}
   On the left, a rhombic torus is 
   represented as a planar rhombus with sides
   $1$ and $\tau$. The paths $\gamma_1$ and $\gamma_2$ correspond to the 
   vectors
   $\overline{0,1}$ and $\overline{0,\tau}$. The ends $E_1$ and $E_2$ are 
   labelled. On the right, the paths $\tilde{\gamma_1}$ and 
   $\tilde{\gamma_2}$ 
   are drawn.
   Each begins by descending $A$, the top half of the vertical diagonal, 
   then follows their respective halves of the horizontal diagonal, avoiding 
   the ends, $E_i$, by  traversing semicircular paths. Each 
   $\tilde{\gamma_i}$ is 
   homotopic to $\gamma_i$.
  } 
 }
 \end{center} 
 \endfig
}

\def\FigHk {
 \begfig
 \hspace{0.2in} 
 \epsfxsize = 4.2in 
 \epsffile{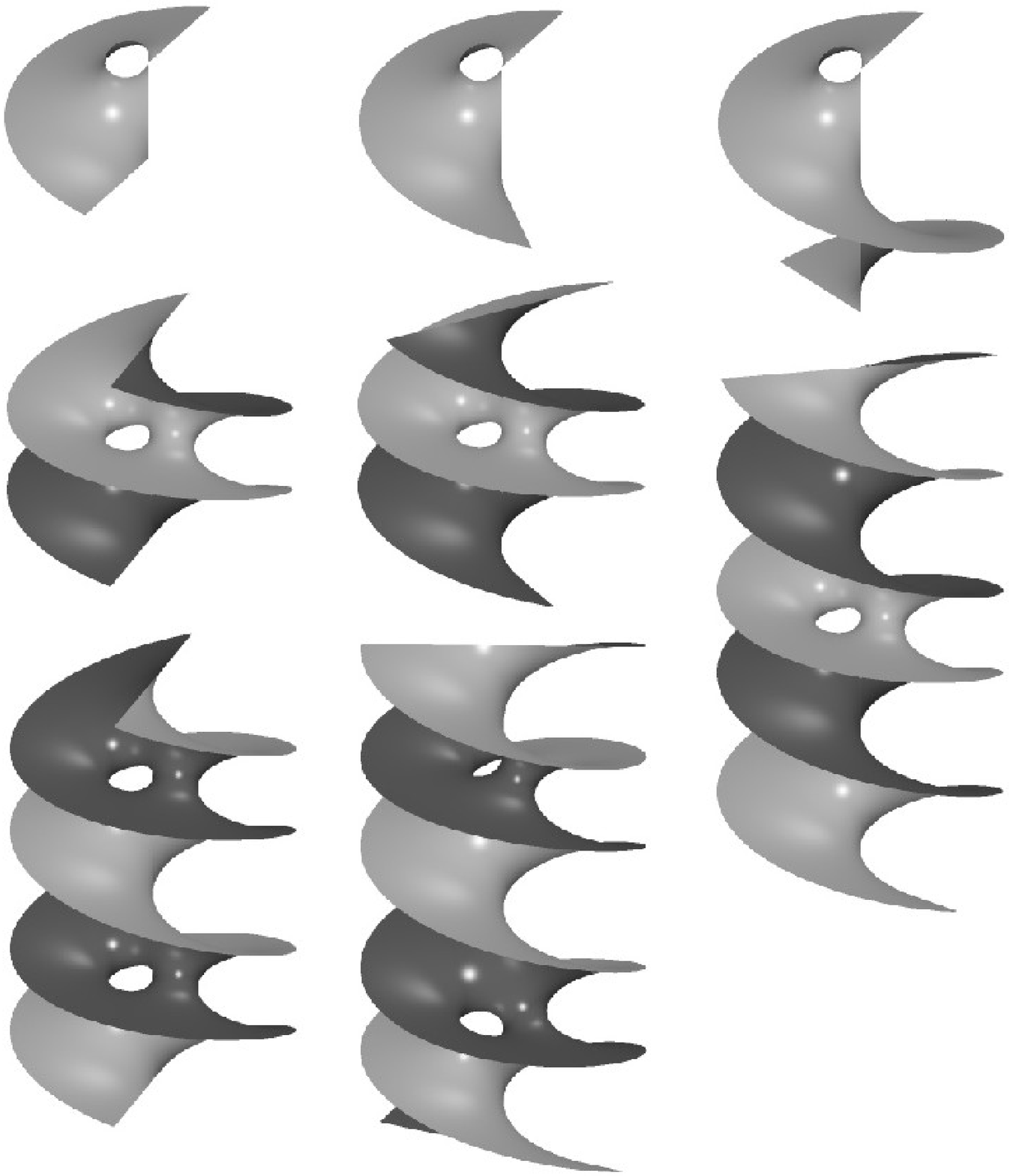} 
 \hfill 
 \begin{center} 
 \vspace{1.0cm}
 \parbox{4.0in}{
  \caption{
   \label{figure:hk}
   Left Column. The surface \Hone/; on top, one quarter of the
   fundamental domain of \Hone/ modulo $\sigma_1$;
   in the middle, one full fundamental domain containing four copies of 
   the region above; on the bottom, two fundamental domains.
   Middle column: The surface \Hk/  for $k\sim 1.25$, with images that
   correspond to those in the first column.
   Right column: The surface \Hk/  for $k\sim 2.5$, with images 
   corresponding
   to the top and middle images of the other two columns.
  } 
 }
 \end{center} 
 \endfig
}

\def\FigPieceHeOne {
 \begfig
 \hspace{0.2in}
 \epsfxsize = 2.0in
 \epsffile{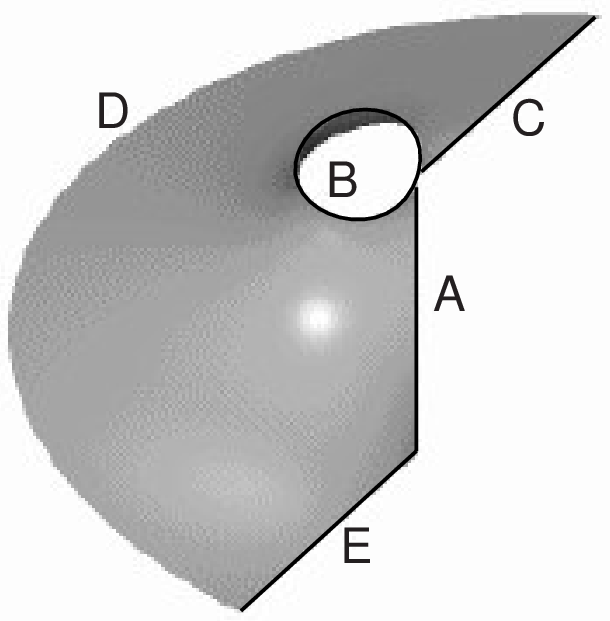}
 \hspace{0.2in}
 \epsfxsize = 2.2in
 \epsffile{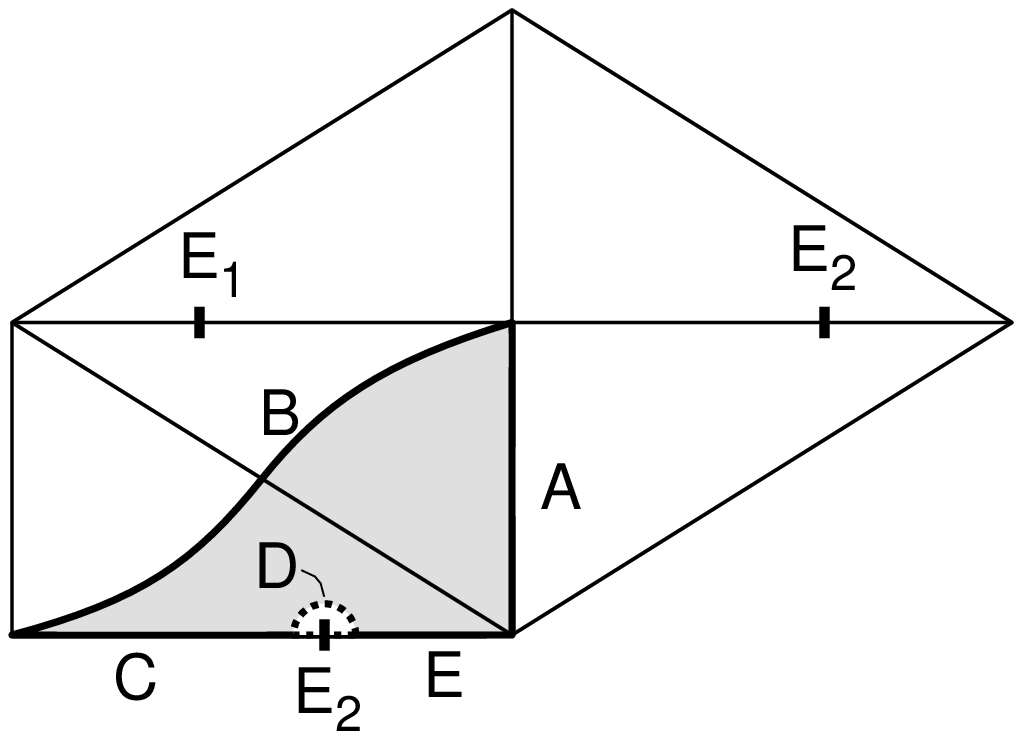}
 \hfill
 \begin{center}
 \vspace{0.3cm}
 \parbox{4.0in}{
  \caption{
   \label{figure:pieceheone}
   One quarter of a fundamental domain of \Hone/ modulo translation is 
   illustrated on the right. It is bounded by a vertical line segment, two 
   horizontal line 
   segments, a portion of a curve that approximates 
   part of a turn of a helix, and a closed loop
   that is a closed cycle on the quotient torus. On the right the associated 
   region of a rhombic torus is drawn, the desired image of which is the 
   minimal surface on the left. Corresponding curves in the two images are 
   similarly labelled. Note that the curve $B$ is drawn to pass through a 
fixed 
   point of $180-\circ$-rotation about the center of the rhombus.  On the 
   minimal surface, this point is a fixed point of the {\em normal symmetry} 
   described in the text.
  }
 }
 \end{center}
 \endfig
}

\def\FigOneOverG {
 \begfig
 \hspace{1.4in}
 \epsfxsize = 1.6in
 \epsffile{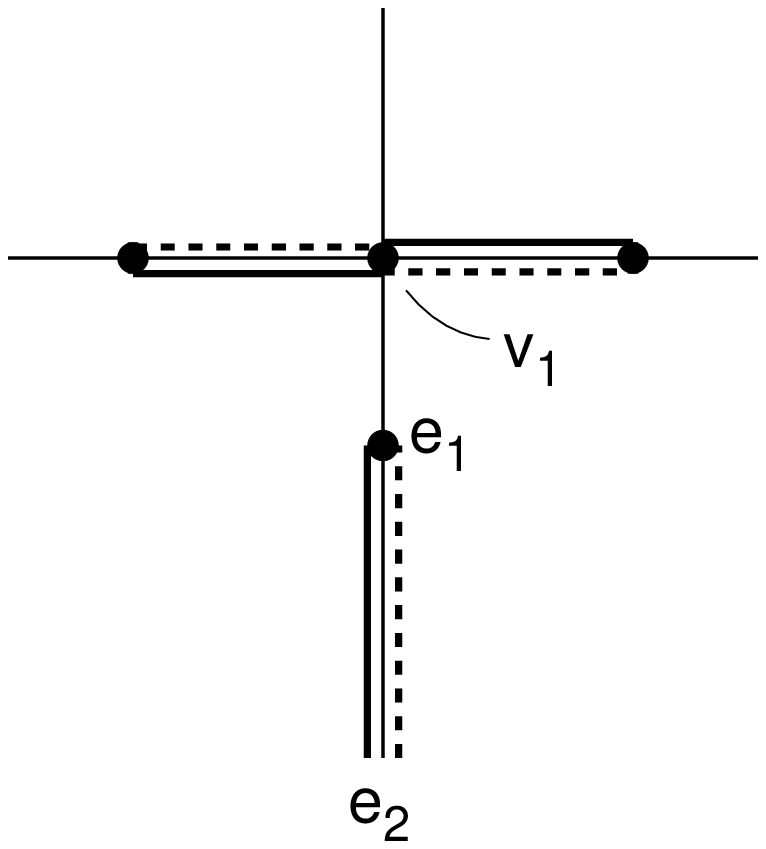}
 \hfill
 \begin{center}
 \vspace{0.3cm}
 \parbox{4.0in}{
  \caption{
   \label{figure:oneoverg}
   The one-form $\f1g dh$ can be considered as 
   the form descended from $d\z$ to 
the slit model of the torus that is constructed in the diagram above. The 
point $e_1$ is placed
   at $-di$.
  }
 }
 \end{center}
 \endfig
}

\def\FigRectangleRTwo {
 \begfig
 \hspace{0.2in} 
 \epsfxsize = 4.2in 
 \epsffile{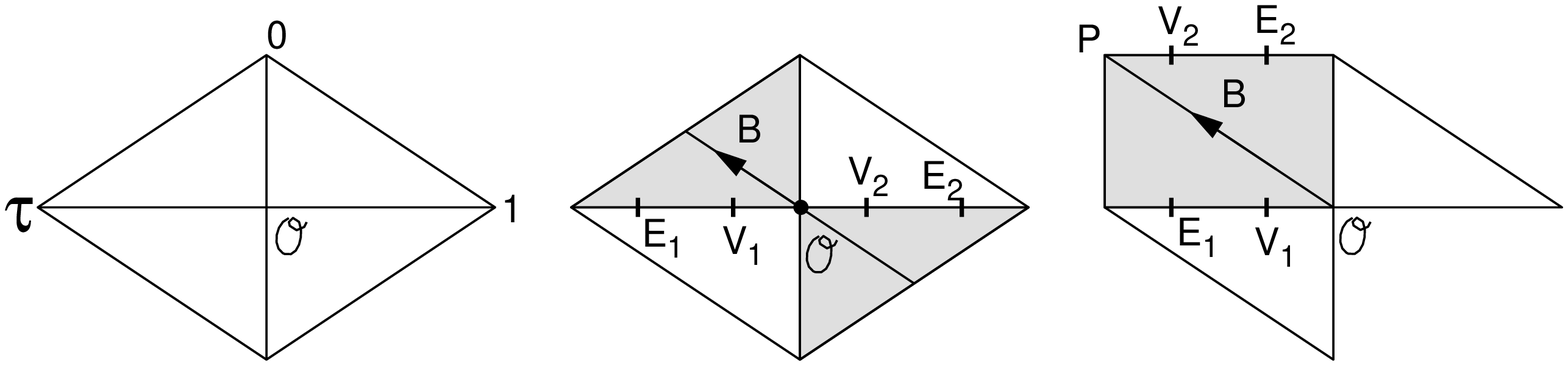} 
 \hfill 
 
 \begin{center}
 \vspace{0.3cm}
 \parbox{4.0in}{
  \caption{
   \label{figure:rectanglertwo}
   The rectangle.
   As in Section~3.4.1, we have the same situation for the
   formation of a rectangular domain that is half of \tkd/ and on which $B$ 
is 
   the path of integration. The only difference is that the underlying rhombic 
   torus is not, in general, equal to \T/.
  }
 }
 \end{center}
 \endfig
}

\def\FigGammaPathsTwo{
 \begfig
 \hspace{0.2in} 
 \epsfxsize = 4.2in 
 \epsffile{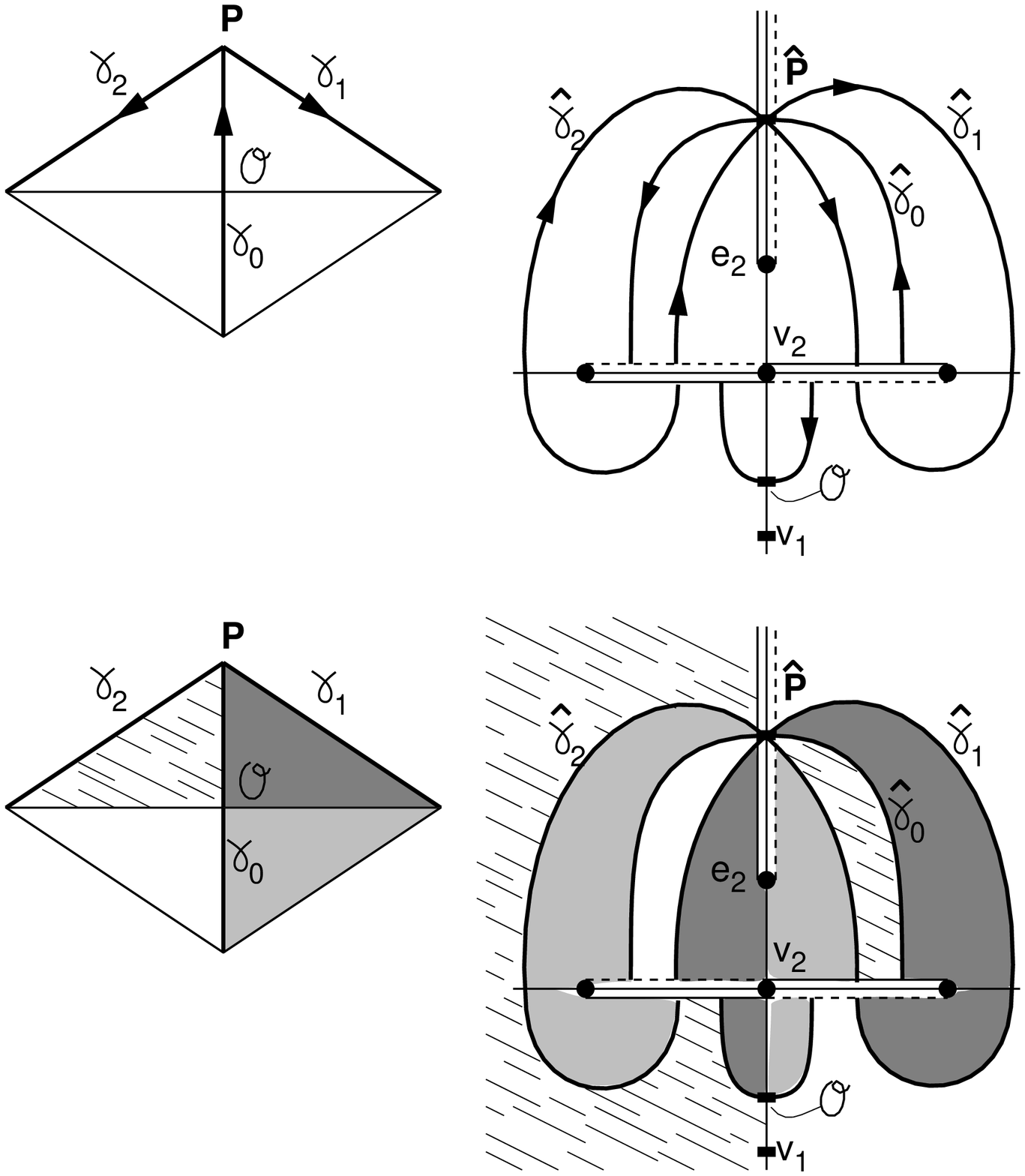} 
 \hfill 
 
 \begin{center}
 \vspace{0.3cm}
 \parbox{4.0in}{
  \caption{
   \label{figure:gammapaths2}
  Upper left: The rhombus with the paths $\gamma_1$ and $\gamma_2$
indicated. The vertical diagonal is labeled  $\gamma_0$. Upper right:
The paths  $\hat{\gamma_i}$---the images of the paths $\gamma_i$ under the 
conformal 
diffeomorphism---are drawn in the slit model. Note that 
$\hat{\gamma_1}$ and $\hat{\gamma_2}$ pass through the image of the vertex 
point 
in 
the 
slit model and no other point on the imaginary axis. The path 
$\hat{\gamma_3}$, 
the 
image of vertical diagonal of the rhombus, passes through both  ${\hat P}$ and
${\cal O}$. Lower left: The rhombus divided into four domains by the 
$\gamma_i$ 
and 
the horizontal diagonal. Lower right: The slit domain divided into four 
domains 
corresponding to those in the rhombus model under the conformal 
diffeomorphism.} 
 } 
 \end{center}
 \endfig
}

\def\FigPathsforaPluskb{
\begfig
\hspace{0.3in}
\epsfxsize = 3.7in
\epsffile{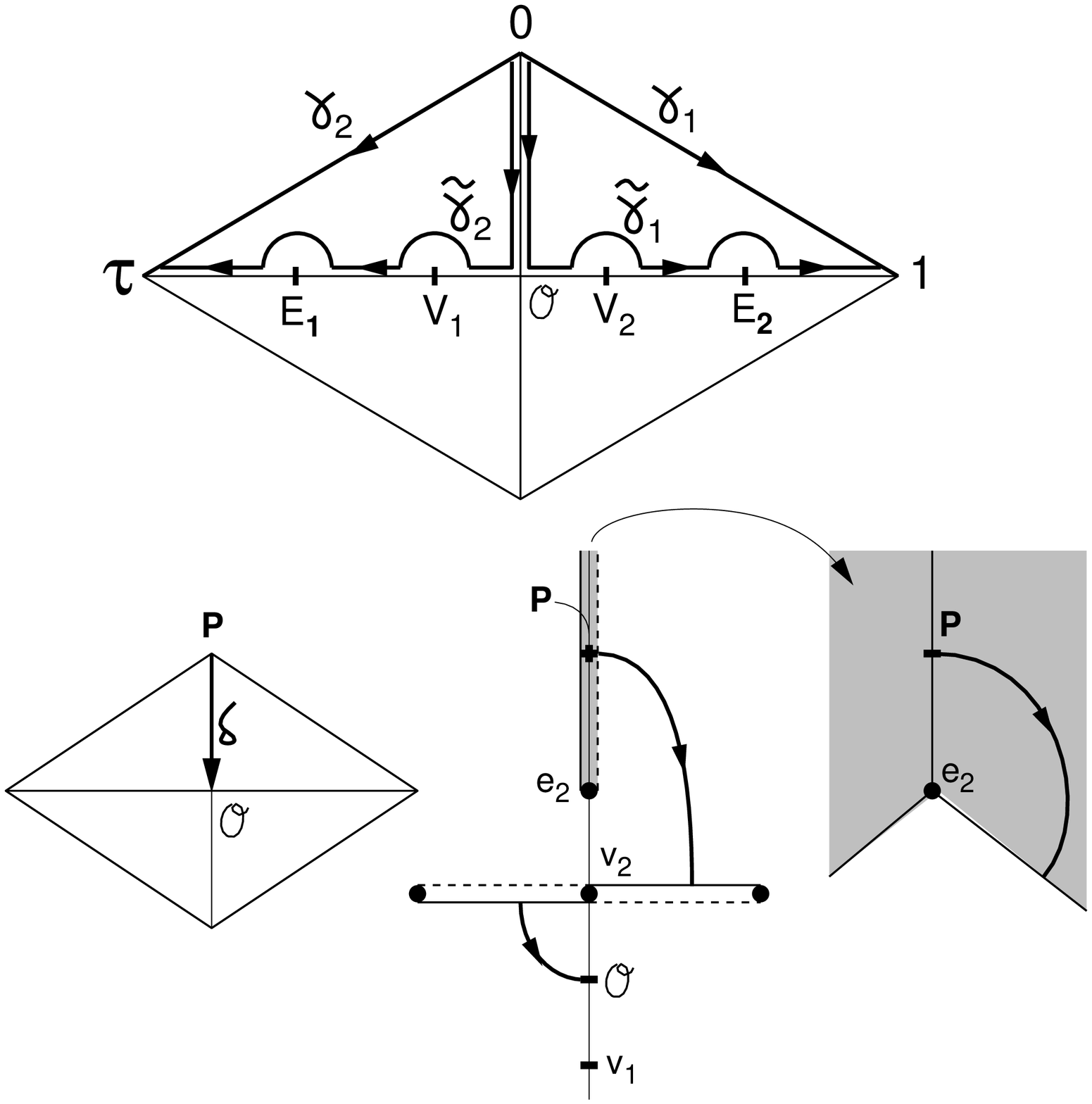}
\hfill
\begin{center}
\vspace{0.3cm}
\parbox{4.0in}{
\caption{The paths of integration in the proof of $a+kb=k$. Top: The rhombic 
model
with paths of integration realized on the torus \tkd/. The paths 
$\tilde{\gamma _i}$ are
homotopic to the paths $\gamma _i$. Each of the paths $\tilde{\gamma _i}$ 
consists of 
the top half of the vertical diagonal,
two semicicular arcs centered at $E_i$ and three line segments on the 
horizontal
diagonal. Bottom left: The part of the arcs  $\tilde{\gamma _i}$ consisting of 
the
top half of the vertical diagonal is labeled $\delta$. Bottom right: The image 
of 
$\delta$ in the slit model is illustrated here. It begins at a point, labeled 
$P$, 
in the 
sewn-in cone (this point corresponds to the vertex point in the rhombic 
model), then
descends to the slit, passing through it to end at the point ${\cal O}$  
corresponding to
the identically labeled center point in the rhombus.
\label{figure:pathsforapluskb}
}

}	
\end{center}
\endfig
}

\def\FigCompactnessLemma{
\begfig
\hspace{1.4in}
\epsfxsize = 1.5in
\epsffile{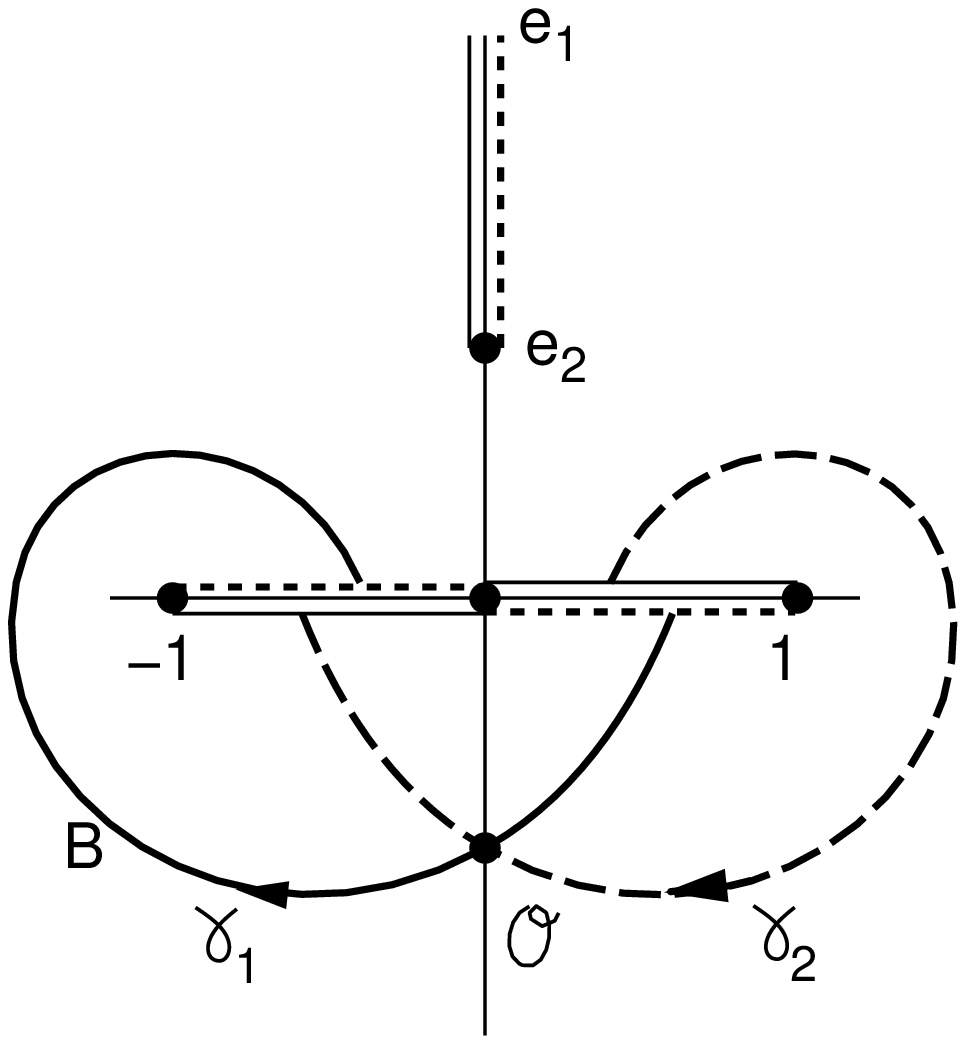}
\hfill
\begin{center}
\vspace{0.3cm}
\parbox{4.0in}{
\caption{
\label{figure:compactnesslemma}
The paths $\gamma_1$ and $\gamma_2$ used in the proof of 
Lemma~\ref{compactnesslemma}.
}
}
\end{center}
\endfig
}

\def\FigDeltaModel{
\begfig
\hspace{1.1in}
\epsfxsize = 2.5in
\epsffile{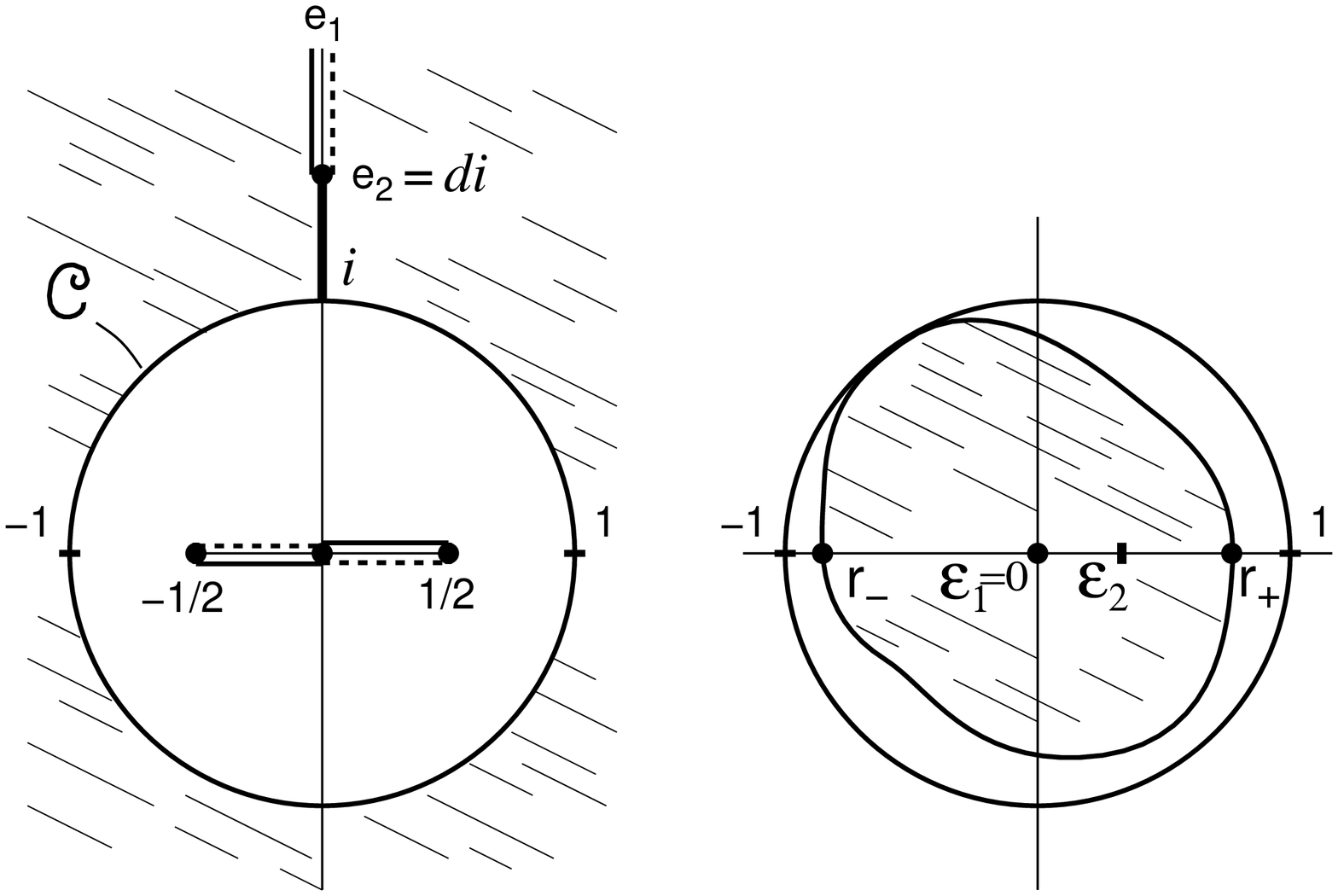}
\hfill
\begin{center}
\vspace{0.3cm}
\parbox{4.0in}{
\caption{
\label{figure:deltamodel}
The $\Delta$ model. On the left: The (scaled) rhombic model with the circle 
${\cal C}$, exterior to which is the  simply connected domain $D$. On the 
right:
the disk $\Delta$, isomorphic to $D$ (with the metric induced by $|d\zeta|$ in 
the
slit model). The disk $\Delta$ has been normalized so that the point 
corresponding
to $e_1$ is at the origin, the point corresponding to $e_2$ is on the positive 
real axis, and $\Delta$ lies inside the unit disk but not inside any smaller 
disk.
}
}
\end{center}
\endfig
}

\def\FigA1D{
\begfig
\hspace{1.4in}
\epsfxsize = 1.5in
\epsffile{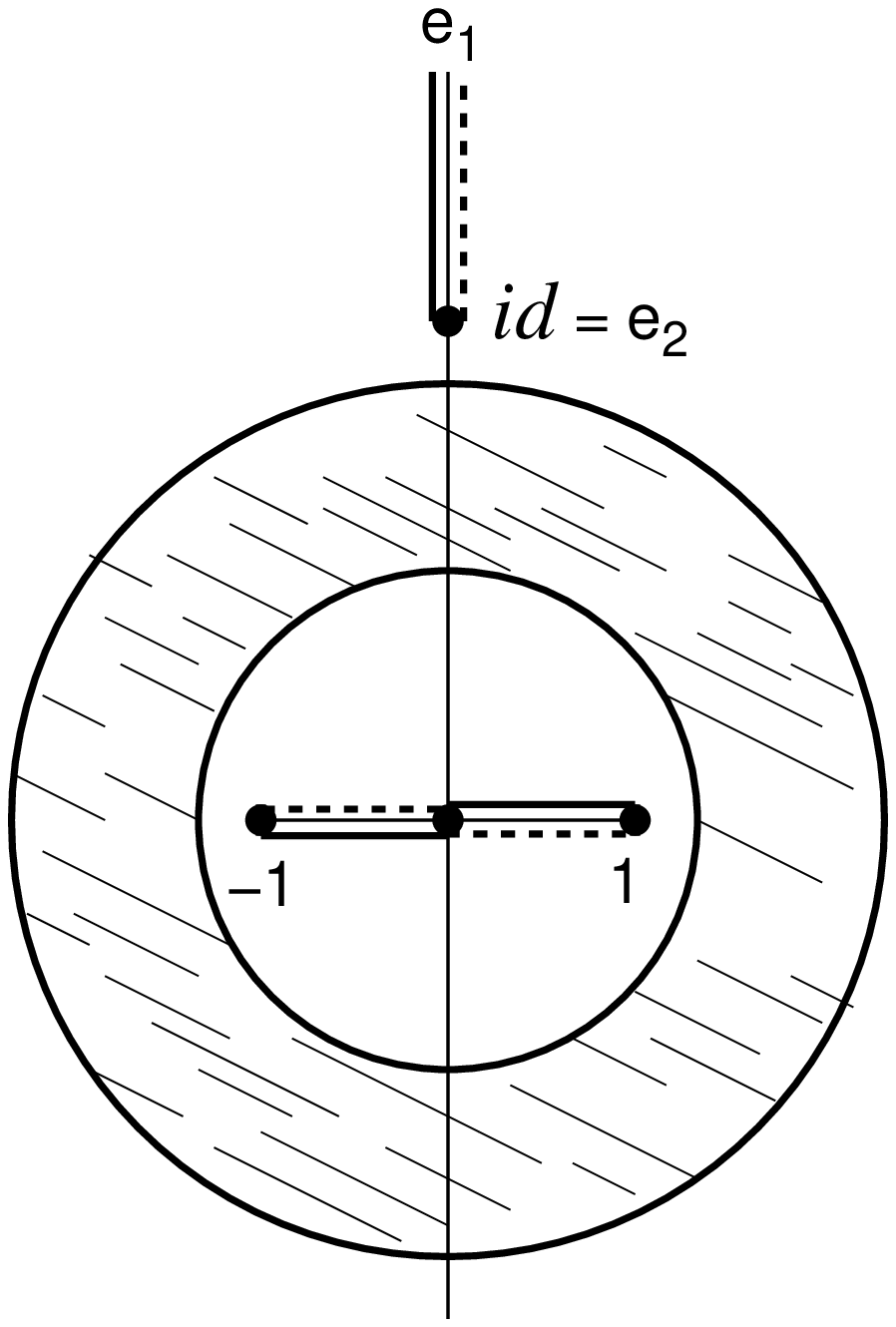}
\hfill
\begin{center}
\vspace{0.3cm}
\parbox{4.0in}{
\caption{
\label{figure:aoned}
The annulus $A(1,d)$.
}
}
\end{center}
\endfig
}

\def\FigTheCurveB{
\begfig
\hspace{0.8in}
\epsfxsize = 2.3in
\epsffile{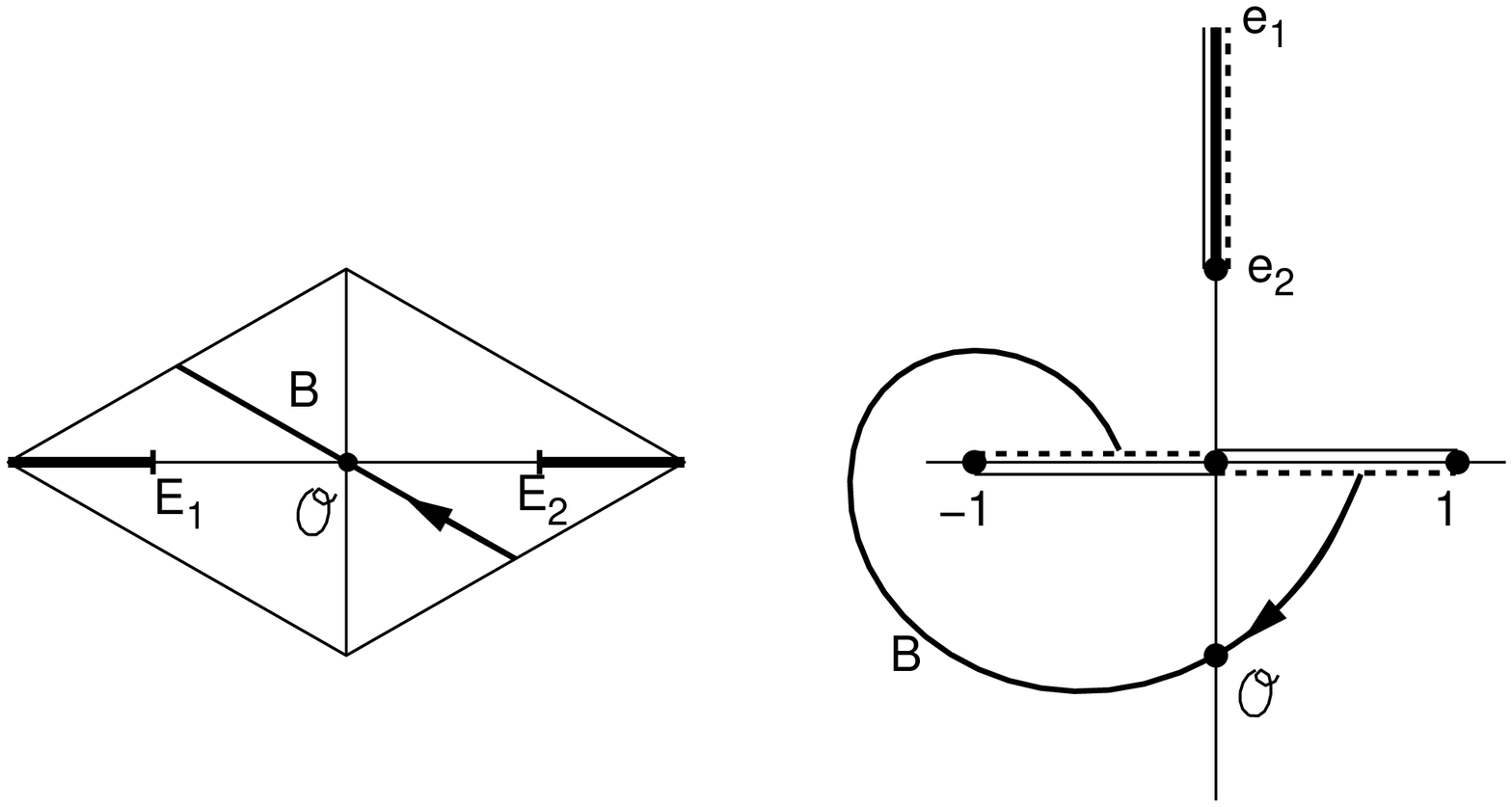}
\hfill
\begin{center}
\vspace{0.3cm}
\parbox{4.0in}{
\caption{
\label{figure:thecurveB}
 The path $B$ of integration for the vertical period problem illustrated
 on the left in the rhombic model and on the right in the slit model. The cut
 referred to in the convention runs in the rhombic model from $E_1$ to $E_2$,
 passing through the vertex of the rhombus. It is indicated by a
 bold line in the
 illustration on the left. Its counterpart runs from $e_1$ to $e_2$ in the
 slit model on the right, and is also indicated by a bold line.
}
}
\end{center}
\endfig
}

\def\FigTheCurveC{
\begfig
\hspace{0.45in}
\epsfxsize = 3.5in
\epsffile{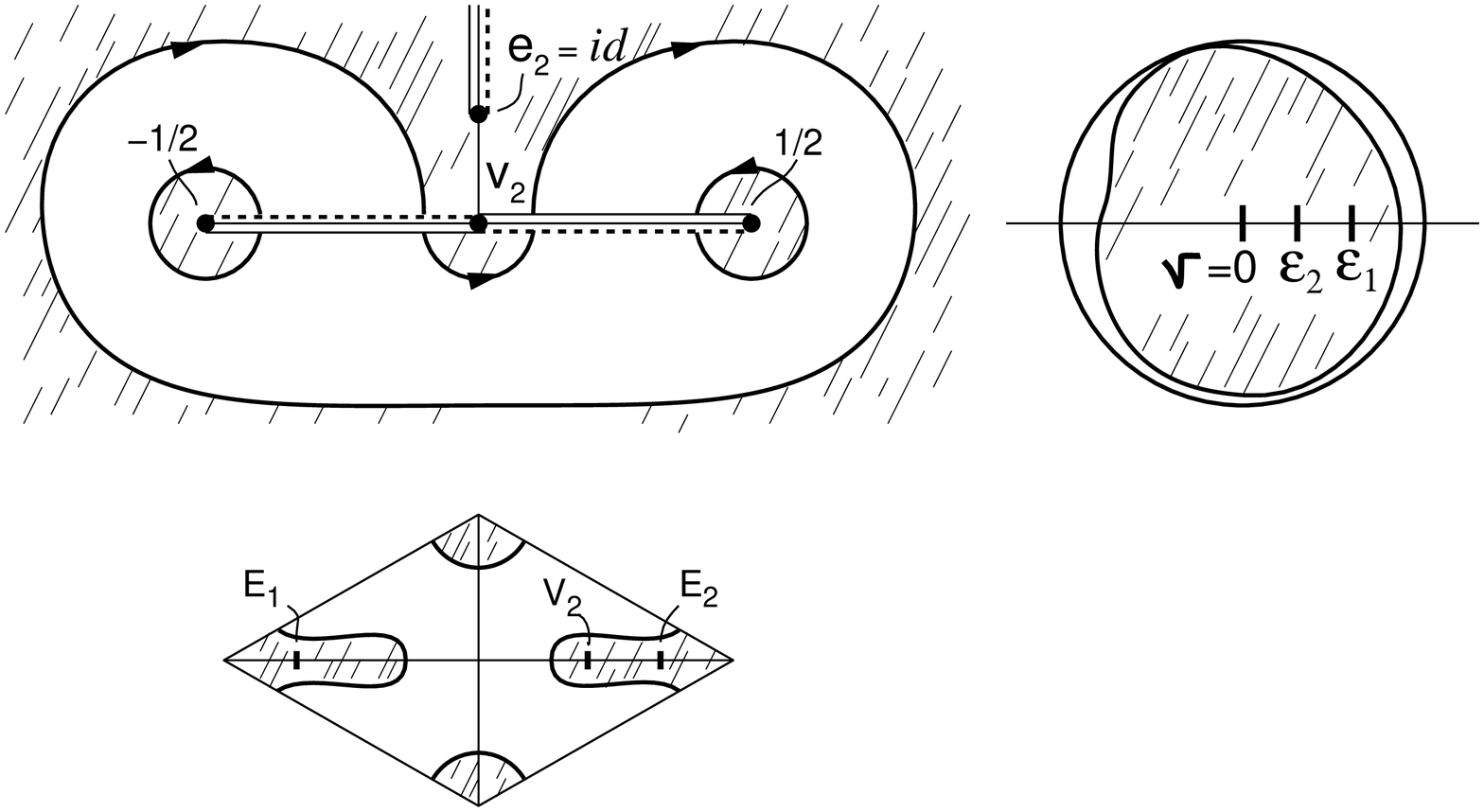}
\hfill
\begin{center}
\vspace{0.3cm}
\parbox{4.0in}{
\caption{
\label{figure:thecurveC}
 Top: The curve $\SC$ in the slit model is drawn on the left as is described 
in the text. It separates the torus into a simply connected  domain $D$  
(shaded) containing $v_2$, $e_1$ and $e_2$, and the complement of $D$, a 
domain of genus one. The 
curve $\SC$ is symmetric with respect to the imaginary 
axis. It consists of two circular arcs of radius $\f1{10}$, a semi-circular 
arc of radius $\f1{10}$ and a curve, symmetric with respect to the imaginary 
axis, that crosses that the imaginary axis once in the bottom halfplane. On 
the right is a represention of the $\Delta$-model normalized as in the 
discussion preceding Lemma~\ref{peskyestimates}. 
The point $\nu_2=0$ corresponds to the point 
$v_2$, and the points $\e_2<\e_1$ on the positive real axis
correspond to the points $e_2=id$ and $e_1 =\infty$, repsectively, in the
slit model. The $ds_{\Delta}$ metric on $\Delta$ is isometric to the 
$|d\z|$-metric on $D$ in the slit model. 
Bottom: Illustrated here is the shaded region in the rhombic model that 
corresponds to the region $D$ bounded by $\SC$ in the slit model.
 }
}
\end{center}
\endfig
}

\def\FigDegenerateOne{ 
\begfig
\hspace{1.6in}
\epsfxsize = 1.3in
\epsffile{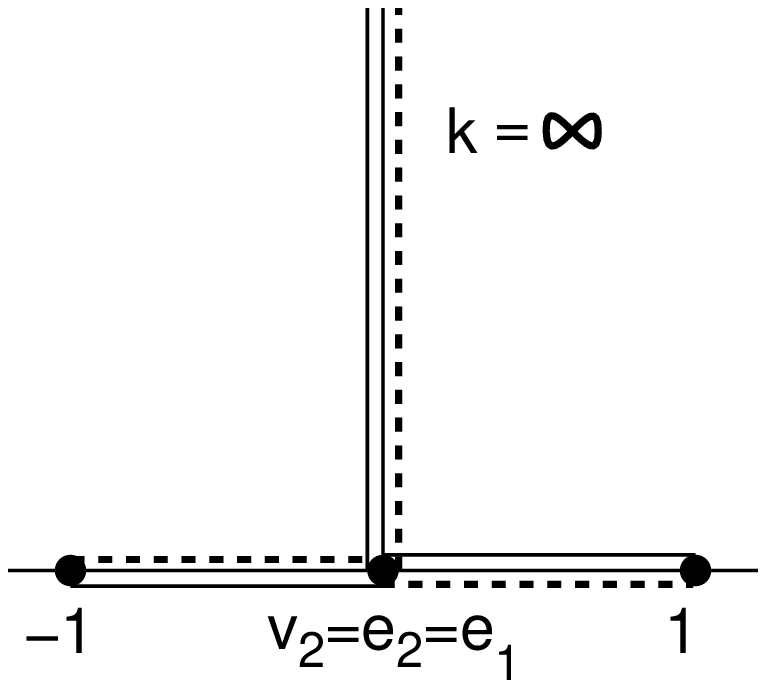}
\hfill
\begin{center}
\vspace{0.3cm}
\parbox{4.0in}{
\caption{
\label{figure:degenerate1}
The degenerate structure  $T_\infty(0)$.
 }
}
\end{center}
\endfig
}

\def\FigWeenydisk{
\begfig
\hspace{1.2in}
\epsfxsize = 2.3in
\epsffile{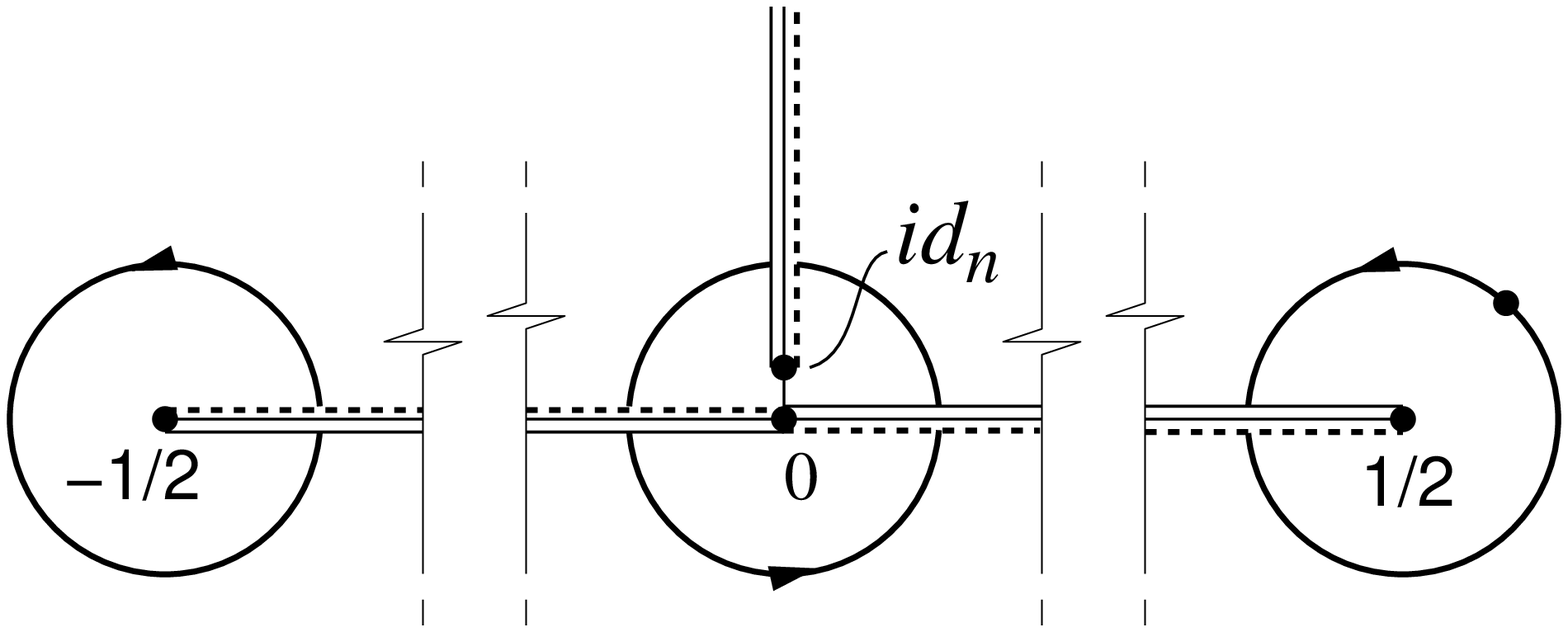}
\hfill
\begin{center}
\vspace{0.3cm}
\parbox{4.0in}{
\caption{
\label{figure:weenydisk}
The core curve of radius $2d_n$, bounding a small disk around $v_2 =0$.
 }
}
\end{center}
\endfig
}

\figon
\begin{document}
\ifx\href\undefined\else\hypersetup{linktocpage=true}\fi

\def\a{\alpha}
\def\d{\delta}
\def\e{\epsilon}
\def\g{\gamma}
\def\G{\Gamma}
\def\la{\lambda}
\def\La{\Lambda}
\def\om{\omega}
\def\Om{\Omega}
\def\th{\theta}
\def\z{\zeta}
\def\f{\frac}
\def\sbq{\subseteq}
\def\cd{\cdot}
\def\tl{\tilde}
\def\vp{\varphi}
\def\wt{\widetilde}
\def\wh{\widehat}
\def\x{\times}
\def\SC{\mathcal C}
\def\SE{\mathcal E}
\def\SF{\mathcal F}
\def\SO{\mathcal O}
\def\SP{\mathcal P}
\def\SR{\mathcal R}
\def\ST{\mathcal T}
\def\SU{\mathcal U}
\def\lb\{{\left\{}
\def\rb\}{\right\}}
\def\ovt{\overline{\mathcal T}}
\def\BC{\ensuremath{\mathbf C}}
\def\BE{\ensuremath{\mathbf E}}
\newcommand{\ms}{minimal surface}
\newcommand{\mss}{minimal surfaces}
\newcommand{\BR}{\ensuremath{\mathbf R}}
\newcommand{\area}{\operatorname{Area}}
\newcommand{\ext}{\operatorname{Ext}}
\newcommand{\genus}{\operatorname{genus}}
\newcommand{\ree}{\operatorname{Re}}
\newcommand{\dist}{\operatorname{dist}}
\def\ov{\overline}
\newcommand{\halo}[1]{\Int(#1)}
\def\Int{\mathop{\rm Int}}
\def\Re{\mathop{\rm Re}}
\def\Im{\mathop{\rm Im}}
\def\bfR{\mathop{\bf R}}
\def\bfC{\mathop{\bf C}}
\newcommand{\union}{\cup}
\newcommand{\goesto}{\rightarrow}
\newcommand{\hra}{\hookrightarrow}
\newcommand{\bdy}{\partial}
\newcommand{\n}{\noindent}
\newcommand{\p}{\hspace*{\parindent}}
\newcommand{\lap}{\bigtriangleup}
\def\({\left(}
\def\){\right)}
\def\[{\left[}
\def\]{\right]}
\def\S2/{${\bf S^2}$}
\def\R3/{${\bf R^3}$}
\def\Hel/{${\cal H}$}
\def\HEone/{${\cal H}e_{_1}$}
\def\Hone/{${\cal H}_{_1}$}
\def\Hk/{${\cal H}_{_k}$}
\def\He/{${\cal H}$}
\def\T/{${\cal T}_{_1}$}
\def\t/{$\cal T$}
\def\tkd/{${\cal T}_k(d)$}
\def\tec{Teichm\"uller\ }
\def\hksk{{\cal H}_{_k}/\sigma_k}
\def\s_k{\sigma_k}
\def\O/{${\cal O}$}
\def\aR/{${\cal R}$}
\newtheorem{theorem}{Theorem}
\newtheorem{assertion}{Assertion}
\newtheorem{proposition}{Proposition}
\newtheorem{remark}{Remark}
\newtheorem{lemma}{Lemma}
\newtheorem{definition}{Definition}
\newtheorem*{claim}{Claim}
\newtheorem{corollary}{Corollary}
\newtheorem{observation}{Observation}
\newtheorem*{convention}{Convention}
\newtheorem{conjecture}{Conjecture}
\newtheorem{question}{Question}
\newtheorem{example}{Example}
\newbox\qedbox
\setbox\qedbox=\hbox{$\Box$}
\newenvironment{proofspec}[1]%
{\smallskip\noindent{\bf Proof of #1.}\hskip \labelsep}%
{\nobreak\hfill\hfill\nobreak\copy\qedbox\par\medskip}
\newenvironment{acknowledgements}{\smallskip\noindent{\bf 
Acknowledgements.}%
\hskip\labelsep}{}

\setlength{\baselineskip}{1.2\baselineskip}

\title{An embedded genus-one helicoid}
\author{Matthias Weber\footnote{partially supported by NSF grant DMS-0139476}
\\
Department of Mathematics\\Indiana University \\
Bloomington, IN 47405\\
\and
David Hoffman\footnote{Partially supported by research grant
DE-FG03-95ER25250/A007 of the Applied Mathematical Science subprogram of the 
Office of Energy
Research, U.S. Department of Energy, and National Science Foundation,
Division of Mathematical Sciences research grant DMS-0139410.}
\\
 MSRI\\ 1000 Centennial Drive\\
 Berkeley, CA 94720\\
\and
Michael Wolf\footnote{partially supported by NSF grants
DMS9971563 and DMS-0139887}\\ 
Department of Mathematics\\   Rice University\\
Houston TX 77005
}
\maketitle

\begin{abstract}
There exists a properly embedded minimal surface of genus one with one end.
The end is asymptotic to the end of the helicoid. This genus one helicoid is
constructed as the limit of a continuous one-parameter family of screw-motion
invariant minimal surfaces---also asymptotic to the helicoid---that have 
genus equal to one in the quotient.
\end{abstract}

\newpage
\tableofcontents

\section{Introduction}

\subsection {An embedded genus-one helicoid}
We prove the existence of a properly embedded \ms\
in $\BR^3$ with finite topology and infinite total
curvature.\footnote {A surface is said to have {\em finite
topology} if it is  homeomorphic to a compact surface with a
finite number of points removed.}   It is the first such surface
to be found since 1776, when Meusnier showed that the helicoid was
a \ms\ \cite{Meu1}. Our surface has genus one and is asymptotic to
the helicoid.

We exhibit this minimal surface as a geometric limit of periodic embedded
minimal surfaces.  The periodic surfaces, \Hk/, indexed by $k\geq
1$, are invariant under a cyclic group of screw motions generated
by $\sigma_k$: rotation by $2\pi k$ about the vertical axis,
followed by a vertical  translation by $2\pi k$. Thus, 
for  fixed $k$, the quotient
surface has two topological ends and genus one. The limit is
taken as $k\rightarrow \infty$; a compact set in 
the limit surface \HEone/ is increasingly well-approximated (as 
$k\rightarrow \infty$) by corresponding pieces of fundamental domains 
of \Hk/$/\sigma_k$.

The requirement to prove embeddedness was the main motivation of
our work. We prove that embeddedness is inherited from the
embeddedness of the approximating simpler (periodic) surfaces, using 
that in this particular minimal surface setting,
the condition of being embedded is both open and closed on families.
This method of proving embeddedness for surfaces defined using the
Weierstrass representation contrasts with previous methods: here
the characteristic of being embedded follows naturally from the
property holding for simpler surfaces, while previously one 
proved embeddedness by ad hoc methods, for instance by
cutting the surface into graphs.
Recent work of Traizet-Weber \cite{TraizetWeber03} suggests
that, ultimately, the embeddedness of this genus-one helicoid
derives from the embeddedness of the helicoid itself.

A second important feature is that we approximate a surface of
finite topology (and finite symmetry group) by surfaces of
infinite topology (and infinite symmetry group). We believe this
is the first example of such a construction resulting in the
existence of a new surface.

The third  feature is that we  construct the Weierstrass data of
these approximating \mss\ in terms of flat singular structures on
the tori corresponding to the quotients. The salient feature to
note is that the defining flat structures have singularities
corresponding to the two ends \footnote{There is an additional 
cone point (with cone angle $6\pi$) in these structures at a 
vertical point.} with cone angles of $\pm2\pi k$. 
Thus, as the size of the
twist tends to infinity, the cone angles also tend to infinity,
with the limit surface---our genus-one helicoid--- represented in
terms of flat cone metrics with an infinite cone angle. This
corresponds to the Weierstrass data for a helicoid, whose Gauss
map has an essential singularity at the end. This required the
development of a theory of singular flat structures that admits
infinite cone angles.

Recently, Meeks and Rosenberg \cite{MR01} have shown that the
helicoid is the unique  simply connected, properly embedded
(non-planar) \ms\ $\BR^3$ with one end. The method of proof uses
in an essential manner the work of Colding and Minicozzi \cite{cm1},
\cite{cm2}, \cite{cm3}, \cite{cm4}
concerning curvature estimates for embedded minimal disks 
and geometric limits of those disks. (Colding and Minicozzi have
recently shown that a complete and embedded minimal surface with
finite topology in $\BR^3$ must be proper. \cite{cm04}.) These results
together with the present work  and
numerical work of Traizet  \cite{trp} and of Bobenko \cite{bobp} (see section 1.3) suggest
that there may be a substantial theory of complete embedded \mss\
with one end, infinite total curvature and finite topology. For
complete embedded surfaces of finite total curvature, the theory
is surveyed in \cite{hk2}.

\subsection{The main theorem}
In 1993 Hoffman, Karcher, and Wei \cite{howe3} constructed a surface, \HEone/
$\subset R^3$, which they called the genus-one helicoid.

It has the following properties: \hfill (1)
\begin{enumerate}
\item [(i)] \HEone/ is a properly immersed minimal surface;
\item [(ii)]\HEone/ has genus one and one end asymptotic to the helicoid;
\item [(iii)]\HEone/ contains a single vertical line ({\em the axis})
and a single horizontal line.
\end {enumerate}

\FigHelicoidZeroOne
We will take the liberty of referring to {\em any} surface
with the properties  (1) as a genus-one helicoid, and of
denoting such as surface by \HEone/.

The \HEone/ in \cite{howe3} was constructed by solving the period
problem for a  Weierstrass Representation (see (\ref{weier}) ---
(\ref{vpc})) chosen to force the surface to satisfy conditions
(1).

Computer-generated images of  this \HEone/, and computational
estimates produced by first solving the period problem
numerically and then triangulating the approximate surface,
showed beyond reasonable doubt that this \HEone/ was
embedded; nevertheless, a non-computional proof has been elusive.
(It is known that any \HEone/ must be embedded outside of
a compact set. See Proposition \ref{embeddedend}.)

We prove that

\begin {theorem}
There exists an embedded \HEone/.
\end{theorem}

We believe that the surface we have found is the same one
constructed by Hoffman, Karcher and Wei. In fact, we
believe {\conjecture There is a unique embedded  \HEone/ . }

\vspace{0.2in}

Note that Conjecture 1 does not assert that there is a unique
\HEone/ and that it is an embedded surface. This does not appear to be
true, as Bobenko \cite{bobe2} has given strong
computational evidence for the existence of an \HEone/ that is immersed
but {\em not} embedded.

Condition (1.ii) implies that any \HEone/ has finite 
topology and infinite total
curvature. The helicoid---a surface swept out by a horizontal
line rotating at a constant rate as it moves up a vertical
axis at a constant rate---is clearly properly embedded and
has finite topology (in fact it is simply connected). Since
it is singly periodic and  evidently not flat, it has
infinite total curvature.

\subsubsection{The place of an embedded \HEone/ in the global
theory of  minimal surfaces}

That complete minimal surfaces in $R^3$ with {\em finite
total curvature} must have finite topology is a consequence
of Osserman's theorem.\footnote{Osserman's theorem states that
a complete minimal surface with finite total curvature in
$R^3$ is {\em conformally} diffeomorphic to a compact Riemann
surface from which a finite number of points have been
removed. Moreover, the Gauss map and coordinate one forms of
such a surface, extend meromorphically to the punctures.}
\cite{os1}\cite{hk2} Finite topology does not imply finite total curvature
for complete minimal surfaces in $R^3$---as the example of
the helicoid shows---but Collin's  solution of the
generalized Nitsche Conjecture \cite{coll} implies that a
properly embedded minimal surface in $R^3$ of finite topology
with {\em more than one end} must have finite total
curvature.  One is naturally led to the following questions.

Let $S$ be a  properly embedded minimal surface of finite
topology with  infinite total curvature and
{\em one end}: \hfill 
\begin{enumerate}
\item [(2.i)] In addition to the helicoid, what are the other
examples?
\item [(2.ii)]Is the end of every $S$ asymptotic to the end of
the helicoid?
\item [(2.iii)] Is the helicoid the unique simply
connected $S$?
\end{enumerate}

In \cite{homc1}, one of the authors (DH) and John McCuan considered
properly immersed minimal ends that are conformally
equivalent to a punctured disk, upon which the Weierstrass
data $dg/g$ and $dh$ both have a double pole.  (This condition
is satisfied by the helicoid) . They also
assume that $dh$ has no residue at the puncture. (This 
condition  must hold if the end actually appears on a
properly immersed minimal surface of finite topology with one
end such as an \HEone/, because $dh$ is holomorphic away from
the puncture) They show that if the end contains a vertical
and horizontal ray, then the end is embedded, and it is
asymptotic to the helicoid. Specifically,

\begin{proposition} \cite{homc1} \label{embeddedend}. Let $E$ be a complete 
minimal 
annular end that is conformally a punctured disk, upon which both $dg/g$ and 
$dh$  have a double pole and $dh$ has no residue. If $E$ contains a vertical 
ray and a horizontal ray then that end is asymptotic to a helicoid. In 
particular, a subend is embedded.
\end{proposition}

This proposition allows us to construct Weierstrass data for an \HEone/ 
that meets the conditions of the Proposition
and then be assured that the resulting surface,
if it exists (i.e. if the period conditions for the Weierstrass data are 
satisfied) must have an embedded helicoidal end.

In terms of the  context for question (2.i)-- (2.iii),
the techniques of \cite{homc1} were used and extended by
Hauswirth, Perez and Romon \cite{hapero1} to study embedded minimal
surfaces of {\em finite type}\footnote{A minimal surface has
{\it finite type} if it is conformally a compact Riemann surface
with a finite number of points removed, and the Weierstrass
data $dg/g$ and $dh$ extend  meromorphically to the
punctures.}, a strengthening of the condition of finite
topology. They prove that questions (1.3ii) and (1.3iii) have
affirmative answers if one makes the additional assumption
that the minimal surface $S$ has finite type. Meeks and
Rosenberg \cite{mr2} showed that the answer to (2.iii) is "yes" under
the assumption that $S$ is also singly periodic. 
They  assume neither bounded curvature nor finite total
curvature of the quotient surface. As mentioned in
Section~1.1, Meeks and  Rosenberg \cite{MR01} recently resolved (2.iii)
in the affirmative. They also proved along the way that the
answer to (2.ii) is also ''yes.'' There is no assumption
of finite type or even that the end of the surface is
conformally a punctured disk.

Concerning (2.i), there is evidence that there are higher
genus examples. Traizet (unpublished) devised a computer
program to generalize the Weierstrass representation in
\cite{howe3} to higher genus and compute and solve the period
problem numerically. This yielded convincing numerical
evidence of a genus two helicoid, analogous to the
surface described in Theorem~3. For genus
three and genus four, Bobenko also has produced examples
computationally. See also Figures~\ref{figure:perihel} and
\ref{figure:g2heli}.

\FigPerihel

\FigG2Heli

\subsection{\HEone/ as the limit of a family of
screw-motion-invariant, embedded minimal surfaces}

The starting point of our investigation is the singly
periodic genus-one helicoid.

\begin{theorem} \label{PG1H} \cite{howe2,howeka4}
There exists a properly immersed,  singly periodic minimal surface
\Hone/, whose quotient by vertical translations:

\hfill 
\begin {enumerate}
\item [(3.i)]has genus one and two ends;
\item [(3.ii)] is asymptotic to a full $2\pi$-turn of a helicoid;
\item [(3.iii)]contains a vertical axis and two horizontal parallel lines.
\end {enumerate}
Furthermore, any surface satisfying conditions $(i)-(iii)$ is embedded.
 \end{theorem}

\FigHelPOne

See Figure~\ref{figure:help1} and the left column of
Figure~\ref{figure:hk} for  images of \Hone/. The fundamental domain
of the surface can be imagined as a modification of one full
turn of the helicoid, bounded by two parallel horizontal
lines that are identified in the quotient. The modification
consists of sewing in  a handle  at mid-level. In fact,

\begin {proposition}\Hone/ is the {\em unique} singly periodic
minimal surface satisfying the conditions $(i)-(iii)$ of
Theorem~\ref{PG1H}. \end{proposition}

It was observed by Karcher that the proposition
follows from the proof of the existence of \Hone/ in
\cite{howeka4}, together with a fundamental result of Weber about
rhombic tori \cite{weber1}. \nocite{weber6} This is discussed in
Section 3. (See Proposition~\ref{rmk2}.)

In 1993, Hoffman, Karcher and Wei, realized that \Hone/ could
be conceived as the limit of a one-parameter family of
deformations of the Karcher's genus-one modification of
Scherk's doubly periodic minimal surface.(See \cite{howe3, howe2} for
details.) Soon after, Hoffman and Wei imagined that \Hone/
could be also be deformed in a manner suggested by the
symmetries of the helicoid. The helicoid is not only
invariant under a vertical translation by $2\pi$, it is also
invariant under {\em vertical screw motions} $\sigma _k$: for
any real number, $k$, 
the isometry $\sigma _k$ is defined to be rotation
by  $2\pi k$  around the vertical axis followed by a vertical
translation by $2\pi k$. For $k\geq 1$, imagine a periodic
minimal surface, \Hk/, invariant under a vertical screw
motion $\sigma_k$  and satisfying the following conditions.

\begin{enumerate} 
\item [] The quotient of \Hk/ by  $\sigma_k$: \hfill (4)
 \item [(i)]has genus one and two ends;
\item [(ii)]is asymptotic to a portion of the helicoid that has
twisted through an angle of $2\pi k$;
\item [(iii)]contains a vertical axis and two  parallel horizontal lines.
\label{hk}
\end {enumerate}

See Figure~\ref{figure:hk}.

Hoffman and Wei defined a Weierstrass
representation that was amenable to numerical solution of the
period problem. The twist angle $2\pi k$ was not specified in
advance and was a calculated function of parameters that
specified the conformal type of the rhombus and the location
of geometrically specified points.  After  normalization to
make the Gauss curvature equal to $-1$ at the intersection of
the axis and the middle horizontal line, they observed (as 
you can in Figure~\ref{figure:hk}) that the handle
rapidly stabilizes and the surface quickly approaches a
helicoid away from the handle.  \cite {howe5} (See 
also the animation found at
\texttt {.www.msri.org/publications/sgp/SGP/indexc.html}  .) This  
reinforced the hope that \HEone/ could be
produced as the limit of the \Hk/ and that embeddedness could
be proved in this manner. However, the form of the
Weierstrass representation was not well-suited for proving
existence of the \Hk/, or continuous dependence on $k$.

The limit of these surfaces as $k\rightarrow\infty$, if it
existed in a geometric sense, should be an \HEone/. More
important, if it could be shown that the family depended
continuously on $k$, then the embeddedness of \Hone/ would be
inherited by the \Hk/. It was hoped that, under controlled
circumstances, embeddedness could be shown to pass to the
limit \HEone/.

We show that this can be done.
\FigHk

\begin{theorem}
        \label{MainThm}
For every $k>1$, there exists a complete,
$\sigma_k$-invariant, properly embedded minimal surface,
\Hk/, whose quotient by $\sigma _k$ satisfies conditions (4).
As $k\rightarrow\infty$, a limit surface exists and is an {\em
embedded} \HEone/, i.e. a properly embedded minimal surface
satisfying conditions (1).
\end{theorem}

\subsection{The ideas behind the proof of Theorem
\ref{MainThm}}

Weber realized that the Weierstrass data $\{g,
dh\}$ for \Hone/ of Theorem \ref{PG1H}
defined one-forms $gdh$ and $(1/g)dh$
that differed by a scale factor and a translation. He
used this to show that the {\em horizontal-period problem}
\eqref{hpc2} completely specifies the
conformal structure of the quotient of \Hone/ modulo
translations. \cite {mweber2}. 

By sewing in a cone of
angle $2\pi(k-1)$, we modify one of the {\em singular flat geometric
structures} used to define \Hone/  to produce candidate flat
structures for \Hk/{}. The position of the vertex of the cone
gives a real parameter $d>0$. (For $k=1$ we do not sew in a
cone, but the choice of $d$ determines the placement of the
ends.) We know from Theorem~\ref{PG1H} and Proposition~\ref{rmk2} 
that, for $k=1$, there
is a unique solution to the vertical period problem and it is
the embedded example \Hone/.  For each $(k,d)\in
[1,\infty)\times [0,\infty)$, the construction gives a
candidate structure on each of which $dh$ is determined up to
a scale factor and the horizontal period problem is solved.
These structures depend smoothly on $(k,d)$. 

For each fixed $k$, we have a single free parameter $d$.
 This 
 free parameter is used to satisfy the vertical period condition as follows:
First, we realize that the flat structure $|dh|$ of $dh$ can be understood qualitatively 
as a planar domain (see figure \ref{figure:DhonRone} right).
Secondly, for the parameter values at their limits, we are able to determine $|dh|$ explicitly
which allows us to apply the intermediate value theorem (see the images on the far right of figures
\ref{figure:DhonRtwo} and
\ref{figure:DhonRthree}).
This solves the vertical period problem for each fixed $k$.

Of course we need more than simply a single solution for each $k$:
we need a \emph{continuous} family $\SC$ of solutions that begins at the point $(1, d_1)$ corresponding to
\Hone/, and crosses each $k=const.$ line segment in the
$(k,d)$ rectangle.
Each point, $(k,d)$, on this curve will then define a
properly immersed minimal surface satisfing the conditions
(4).  

To do this, we note that the period of the height function for each $(k,d)$-structure
naturally defines a real analytic function, $h=h(k,d)$ whose
zeros are what we seek: if $h=0$ for some $(k,d)$, then  the
corresponding
$(k,d)$-flat geometric structure defines a properly immersed
minimal surface \Hk/   satisfying conditions (4).
Then, to find this curve $\SC$,
we first compactify the ``moduli space'' of $(k,d)$-structures,
adding in  degenerate structures for the loci $k=\infty, d=0$
and $d=\infty$ in a manner compatible with the topology of 
$[1,\infty] \times [0, \infty]$.  
Moreover, we show that the height
function $h$ is {\em continuous} on the full compact rectangle
$[1,\infty]\times[0,\infty]$ (i.e., $h$ extends continuously
to   to $d=0, \infty$ and $k=\infty$).
This has the advantage that
the signs of $h$ on the degenerate surfaces $(k,0)$ and
$(k,\infty)$ are evident and opposite, and so the intermediate
value theorem provides for a  solution $(k_0,d)\in \SC$ for
each choice of $k_0 \in (1,\infty)$.
More precisely, there must be a curve $\SC$ on which 
$h=0$ because this curve
separates the neighborhoods of the boundary components
$\{d=0\}$ and $\{d=\infty\}$ on which $h$ has opposite signs.
 Finally, a
maximum-principle argument then shows  that the embeddedness
of \Hone/ implies that all the \Hk/, $k>1$, on this curve
$\SC$ are embedded.

The $(k,d)$-structures are defined for $k=\infty$ by
sewing in an infinite cone. Thus, any structure
$(\infty, d)$ corresponds to a potential \HEone/, which will
exist provided $h(\infty, d)=0$; the endpoint
of $\SC$ on the locus $\{k=\infty\}$ is then such a point.
Thus we obtain
a limit flat structure that defines an \HEone/ and we argue
that, as a limit of a family $\{{\cal H}_k\}$ of embedded
surfaces,  the surface \HEone/
is also embedded.

\subsection{An outline of the paper}

In Section~\ref{background} we review the Weierstrass
representation and the associated period problem for minimal 
surfaces invariant under a screw motion. We introduce cone
metrics and establish for them existence and uniqueness results
necessary for our  work. A short review of extremal
length is presented.  The helicoid is  presented from both
the point of view of the Weierstrass representation and  of
singular flat structures (cone metrics). 

Section~3 is devoted to the singly periodic genus-one
helicoid, \Hone/,  of Theorem~2. In Section~3.1, we give a
derivation of the Weierstrass data for this surface under the
geometric assumptions of Theorem~2. In  Section~3.2, we state
the results of \cite{howeka4} about the existence and \Hone/:
the solution of the period problem.  The rest of Section~3 is
devoted to an  alternate construction of candidate data and
solution of the period problem for \Hone/ using singular flat
structures. 

Section~4 is devoted to a generalization of the cone-metric
construction  Weierstrass data  for \Hone/ in Section~3 to
the construction of analogous Weierstrass data for  the
surfaces \Hk/ of Theorem~3. This involves three singular
flat  structures corresponding to the one forms $gdh$, $\f1g
dh$ and $dh$. These  structures are naturally indexed 
in a rectangle \aR/ 
by two
real variables which control the twist
angle and the conformal type of the underlying  punctured
tori. These structures are generalized to the boundary of
\aR/ by  introduction of cone metrics that are natural limits
of the candidate cone  metrics in  the interior. One of the
limits involves letting the cone angle  tend to infinity;
this is what
we expect for the genus-one helicoid \HEone/.

Section~5 is devoted to the proof that the structures behave
continuously  on the extended rectangle \aR/. In particular,
the vertical period  ---a well-defined real number for each
structure in the interior---actually  extends to a continuous
function on the closed extended rectangle. We refer to this 
function as {\em the height function}. The surfaces \Hk/
correspond to the zeros of this function on the interior of
the rectangle.

In Section~6, we prove the first part of Theorem~2 by showing  that there 
is a continuous family of \Hk/ that begins with \Hone/ and is
defined for all $k\geq1$. We do this by an
intermediate-value-theorem argument using the height function
on \aR/. We then show that the limit structure as 
$k\rightarrow \infty$ produces an \HEone/. Embeddedness of
this \HEone/ is then established by using the fact
(Section~3.2) that \Hone/ is embedded  and an argument that
embeddedness propagates along the curve of \Hk/  structures.

In the Appendices, we give an alternate Weierstrass
representation of the \Hk/ using theta functions, and we
present a proof of the existence and  uniqueness results for
cone metrics that come up in Section~2.

\subsection{Acknowledgements}

The authors wish to thank Hermann Karcher for contributing to
this paper in many ways. As a collaborator, colleague,
advisor and friend, his critical intelligence and spirit have
been fundamentally important to our work.

We also wish to acknowledge useful conversations with Harold
Rosenberg, Bill Meeks and Fusheng Wei over the past years.
Wei's ingenious ideas about how to represent the family \Hk/
in a computationally tractable manner was key to the
establishing the belief that that the family not only existed
but also converged to a genus-one helicoid. The idea that
embeddedness was inherited from \Hone/ by the family \Hk/ came
from discussions that one of us had with Meeks about the same
phenomena for finite-total-curvature surfaces. Rosenberg's
pivotal idea to sew a helicoid into a genus-one surface to
create a genus-one helicoid, conveyed to Karcher in a
conversation reported in (\cite{howe3}), led to the discovery of
\HEone/. However,  the construction of \HEone/ in  \cite{howeka4}
was not that direct. In this paper, we are able to do this in
a concrete manner, realizing directly Rosenberg's idea,
but in a way he did not imagine. (See
Section 4.) 

\setcounter{equation}{4}

\section{Minimal surfaces, cone metrics and the geometry of
\Hk/ }\label{background}

\subsection{General background}

We begin with the Weierstrass representation of a minimal
surface in \R3/. Details can be found in \cite{hk2} or  \cite{os1}.
                                                                            
Let {\cal S} be an oriented minimal surface in \R3/. The
metric on {\cal S} induced by the immersion is analytic and
allows us to consider {\cal S} to be a Riemann surface. We
write $M$ for this Riemann surface. The immersion
$X:M\rightarrow$\R3/ is by definition conformal. Minimality
of $S$ is equivalent to the anticonformality of the Gauss map
$N:M\rightarrow  S^2$, which in turn is equivalent to the
conformality of $g:=\sigma\circ N$  where $\sigma$ is
stereographic projection to the extended complex plane. 

Minimality is also equivalent to the harmonicity of the
conformal immersion $X$. In particular, if $\hat{f}$ is a
linear function on $R^3$ (a coordinate function for example) 
then $\hat{f}$ is harmonic on $M$. Let $\hat{f^*}$ be the
(locally well-defined) harmonic conjugate of $\hat{f}$. Then
$f:=\hat{f}+i\hat{f*}$ is a locally  defined holomorphic
function and $df$ is a globally well-defined one-form on
$M$.  In particular, we define the one-form $dh$ to be the
exterior derivative of  the $h=x_3 +i{x_3}^*$, where $X=(x_1,
x_2, x_3)$. We will refer to $dh$ as  the ''height
differential.''

The Weierstrass representation allows one to write the
conformal parametrization $X$ in terms of $g$ and $dh$:
\begin{eqnarray}
X(p) =  Re\int_{p_0}^p \Phi, \mbox{ where   }
   \Phi = (\frac{1}{2}(g-g^{-1}), \frac{i}{2}(g +g^{-1}), 1)dh, \mbox{
  } p,p_0 \in M.    \label{weier}
\end{eqnarray}
The immersion will be regular provided the induced metric
\begin{equation}
ds=|gdh| +|\frac{1}{g}dh|\label{metric}
\end{equation}
is nowhere zero.
This requires the zeros of $dh$ to coincide with (and have
the same order as) the zeros and poles of $g$ (points where
the Gauss map is vertical).

The integral formula (\ref{weier}) can be used to construct
minimal surfaces. Given a Riemann surface $M$, a meromorphic
function $g$, and  a holomorphic one-form $\eta$,  on $M$,
the integral (\ref{weier}) -- with $\eta$  substituted for 
$dh$ -- defines a conformal and harmonic  mapping of $M$ into \R3/
whose image is a minimal surface. The mapping will be regular
provided $|g\eta|
+|\frac{1}{g}\eta|\neq 0$. The  stereographic projection of
the Gauss map of this surface will be $g$ and the 
one-form $\eta$
will be equal to the holomorphic one-form $dh$ that is 
constructed above from $h=x_3$. When constructing
minimal surfaces by  specifying ''Weierstrass data,'' that
is, when specifying $M$, $g$ and $\eta$, we will use the
notation $dh$ for $\eta$.

The Weierstrass representation is, in general, multivalued.
In order for (\ref{weier}) to be single-valued on $M$, it is
necessary and sufficient that $$Re \int_{\alpha} \Phi = 0$$ for
all closed cycles $\alpha$ on $M$. Using (\ref{weier}), this
can be rewritten as 

\begin{eqnarray}
 \overline{\int_{\alpha} gdh} - \int_{\alpha} (1/g)dh&=&0
\mbox{ \em (Horizontal Period Condition)} \label{hpc}
 \\
Re \int_{\alpha}dh&=&0  \mbox{ \em (Vertical Period Condition)}\label{vpc}
\end{eqnarray}

In dealing with {\em periodic} minimal surfaces $N$, it is
often useful---and sometimes necessary--- to work with $g$
and $dh$ on the Riemann surface of the quotient of $N$
by translations or screw motions. We will be dealing with
singly periodic surfaces invariant under screw motions;
without loss of generality, we may assume that the
translational part of the screw motion is vertical,  and that the axis of 
the screw motion is the
$x_3$-axis. Let $M$ denote the Riemann surface of the quotient
surface under the screw motion $\sigma_k$. The Weierstrass
representation \eqref{weier} defines an immersion
of $\tilde{M}$, the universal cover of $M$, into $\BR^3$.
Screw-motion invariance means that there is a basis $[\a_i]$
for the homology of $M$ such that if $A_i$ is the deck
transformation associated to $[\a_i]$, then 
$X \circ A_i =0$ or $X \circ A_i = \sigma_k X$.
\subsection{The helicoid}
\label{helicoiddef}
The helicoid, \He /, is a singly periodic minimal surface
swept out by horizontal lines moving at a contstant speed up
the $x_3$-axis while rotating at constant speed. (See 
Figure~\ref{figure:helicoidzeroone}.) It is invariant under
any vertical screw motion around the $x_3$-axis, in particular
vertical translation by $2\pi$. In the quotient of \He/  by
$\sigma_k$, the screw motion with twist angle $2\pi k$, the
Weierstrass data on $M={\bf C}-\{0\}$ can be chosen to be 
\begin{equation}
\label{helicgdh}
g=iz^k \mbox{\,and\,}
dh=\frac{kidz}{z}.
\end{equation}
When $k$ is not an integer, the Gauss map on the quotient is
multivalued, and so is $g$. From the Weierstrass
representation (\ref{weier}) we have
$$
2(x_1+ix_2) = \overline{\int{gdh}}- \int{(1/g)dh}=
\overline {z}^{k}+z^{-k}$$ 
$$x_3 =-k\cdot arg(z).$$

The only relevant cycle is represented by a circle
$\alpha$ around the origin. 
Hence $X \circ A_{[\a]} = \sigma_k X$, assuming $\a$ is
oriented in a clockwise direction.
The Riemann
surface upon which $dh$ is well-defined is the Riemann
surface of $w= ln(z)$, i.e. ${\bf C}$. 
\FigDivisorsZero 
The globally defined function $z=e^w$ on ${\bf C}$ allows the
expression of  the Weierstrass data for the helicoid in a
univalent manner: $g=ie^w, dh = -idw$. This gives a global
representation on ${\bf C}$ of \He/. Note that in the
representation on the quotient surface, ${\bf C}-\{0\}$,
there are two ends (one at $0$, the other at $\infty$) at
which $g$ has a simple zero and a simple pole, respectively, 
while $dh$ has a simple pole at both ends. In the global
representation, there is a single end at infinity where $g$
has an essential singularity and $dh$ has a double pole.

We will be using the helicoid as a model and have need to
restate the well-known facts above in terms of cone metrics.

\subsection{Flat cone metrics}\label{conemetric}
 
Consider, for $k>0$, the set $A_k$, described as an
identification space
$$
A_k=\{(r,\th)|0\le r<1,\mbox{\,\,}
|\th|\le\pi k\}/\sim,
$$
where $\sim$ identifies the top and bottom edges and
collapses the left  hand edge: 
$(r,-\pi k)=(r,\pi k)$ and $(0,\theta)= (0,0)$. 
Via the identification $(r,\th)\rightarrow z=re^{i\th}$, we
may regard $A_k$ as a possibly multisheeted ``sector'' with
vertex at the origin $0\in\bfC$ and edges identified. Away
from the origin, the sector $A_k$ inherits the flat metric
$|dz|$ on the plane. This metric is the metric of a flat cone
with vertex at the  origin.  We observe that a neighborhood
of the vertex is isometric to a  neighborhood of the vertex
of another identification space $A_{k'}$ if and only if
$k=k'$: this is because the circumference of a circle of
radius $\e$ linking the distinguished point of $A_k$ has length
$2\pi k\e$.  

We extend this definition of $A_k$ to $k=0$ by taking $A_0$
to be the identification space of 
$\{z\in \mathbf{C}| |Imz|\le 1 , Re z \le 0\}$, 
where we identify  the boundary rays by a
vertical translation and  consider the distinguished point
(the vertex) to be a point that compactifies the left end of
the cylinder.

A neighborhood of the vertex of $A_k$ is topologically (and
conformally) a disk.  Consider the map
$w:\mathbf{D} \to A_k$ from the  disk
$\mathbf{D}$ to a neighborhood of the distinguished
point given by $z=w^k$ for $k \ne 0$, and $z=\log w$, when $k=0$. We can 
pullback
the metric on $A_k$ to $\mathbf{D}$: 
$$
|dz| = |kw^{k-1}dw|.\label{ck}
$$
In particular, we
see that a metric $ds = |w^{\a}dw|$ on the disk 
$\mathbf{D}$ defines a metric isometric 
(up to a scale factor, which is unimportant in the 
present discussion)
to one in a
neighborhood of the distinguished point on $A_k$, with
$k=\a+1$. In particular, we may define for $k=0$ the metric
$|\frac{dw}{w}|$ on the   unit disk, and consider it to be a
flat cone with cone angle zero at the  vertex.

We now extend the range of definitions of these local
neighborhoods to all $k\in \mathbf{R}$ by working with these
flat singular metrics defined on the disk $\mathbf{D}$.

\begin{definition} {\em Cones and cone points with cone angle
$2\pi k$}. The cone $C_k$ is defined to be the  disk
$D$ with  the metric given up to scaling by $|w^{k-1}dw|$
on $\{w ||w|<1\}$. The origin is the cone point of the cone
$C_k$.
\end{definition}

Note that when $k=1$, this formula is the standard  regular
metric $|dw|$ on the  $w-$disk. We will adopt the
convention of not referring to the origin in $C_1$ as a cone
point but  as a regular point.  Also note that the metric
$|w^{k-1}dw|$ on the (compactified) exterior of $\mathbf{D}$ defines a
metric with vertex at infinity that is isometric  to
$C_{-k}$. This follows immediately from pulling the exterior 
domain back to $\mathbf{D}$ by $w\rightarrow\frac{1}{w}$ and
pulling back the  metric to the disk.

We will need to consider infinite cones of a specific type.

\begin{definition} {\em An exponential cone of simple type}.
The cone $C_e$ is the  disk $D$ with the
metric $|e^{\f1w}\f{dw}{w^2}|$. 
 The origin is the cone point of $C_e$.
\end{definition}

Note that the metric in the definition of $C_e$ is the
pullback of the metric $|e^{w}dw|$ on the (compactified) exterior of the
unit disk. Therefore we may (equivalently) consider $C_e$ to be
the (compactification of the) 
exterior of the unit disk with  metric $|e^{w}dw|$ and cone point at
$\infty$.

The following definition is nearly standard (see \cite{Tr86}):
our extension  allows the presence of cone points with negative cone angles 
and
exponential cone points of simple type. 

\begin{definition}{\em Cone metric.} Let $M$ be an oriented
surface and let $\{p_1,\dots,p_n\}$ be a discrete set
of points of $M$. A {\em flat cone metric on $M$}  is a
metric on $M$ so that every point on $M$ has a  neighborhood
as follows:

\noindent(1) {\em Regular points.} If 
$q\notin\{p_1,\dots,p_n \}$, then $q$ has a
neighborhood that is isometric to a neighborhood in the
Euclidean plane.

\noindent(2) {\em Cone points.} Every point $p_i$ in the 
distinguished set has a neighborhood that is isometric
to a neighborhood of the vertex in either $C_k$ or $C_e$.  In
the first case, the point $p_i$ is  a cone point with finite
cone angle and in the second case, it is a cone point with
cone angle of simple exponential type. 
\end{definition}

The definition naturally requires transition maps between
neighborhoods to be Euclidean isometries. This defines a
conformal, hence complex, structure on 
$M - \{p_1,\dots,p_n\}$.  Further, deleted
neighborhoods of the cone points are clearly conformally
punctured disks, so the Riemann surface structure on $M -
\{p_1,\dots,p_n\}$ then extends naturally to a
Riemann surface structure on $M$.

The simplest example of a cone metric is the extended complex
plane with metric $|dz|$. There is a single cone point at
infinity with cone angle $-2\pi$. More generally, consider a
compact Riemann surface $M$ and  $\omega$ a  meromorphic one
form on $M$.  Define the metric $ds=|\omega|$. Since 
$log|\omega|$ is harmonic, the Gauss curvature of the metric
is zero: $K=\frac{-\Delta\log|\omega|}{|\omega|^2} =0$, valid
away from the poles  and zeros of $\omega$. A zero (resp.
pole) of order $k$ of $\omega$ represents a cone point with
cone angle $2\pi (k+1)$ (resp. $2\pi(-k+1$)) of the metric
$|\omega|$. This is evident by  looking at local expansions
but it is useful to show this by use of the  developing map
$F$ given by
$$
F(p) = \int^p_a \omega: M - \{p_1,\dots,p_n\}\to\bfC,
$$
where $\{p_1,\dots,p_n\}$ are the zeros and poles of $\om$.
Values of $F$ at a point $p$ differ by periods of $\om$ so we
can pull back the flat metric on $\bfC$ to a well-defined
flat metric on $M-\{p_1,\dots,p_n\}$. (This gives another way
to show that the metric is flat away from the poles and zeros
of $\om$.) Near $p_i$, the surface $M$ has a local chart in which $\om$
takes the form $\om=z^{k_i}dz$. If $k_i\neq-1$, then we can
explicitly integrate to obtain that
$$
F(z)=\f{z^{k_i+1}}{k_i+1}
$$ near $p$. Thus, $p_i$ is a cone
point with cone angle $2\pi (k_i+1)$. 

In the previous paragraph, a cone metric was defined via
metric expressions. As we shall see in Examples 2 and 3
below, they can  also be pieced together from pieces of flat
cones like $C_k$ or $C_e$. The two methods of construction are 
related via the developing map.  It is important to our approach
to be able to pass freely between the two descriptions.

There is a natural version of the Gauss-Bonnet formula for
cone metrics on a surface $M$ with finite cone angles. 
Let $\g_i$ be a small circle around $p_i$,
$i=1,\dots,n$. Then $\bigcup^n_{i=1}\g_i$ bounds a connected
flat surface, $M'$, whose Euler characteristic is $2-2g-n$,
where $g=\genus(M)$. Each $\g_i$ has total geodesic curvature
on $M'$ equal to $-2\pi k_i$. The Gauss-Bonnet formula for
$M'$ gives
\begin{equation}
\sum^n_{i=1}k_i = n + 2(g - 1).\label{GB}
\end{equation}
We will see in Proposition~3 below that there is an extension
of this necessary condition to the case of cone metrics with
exponential cone points of simple type.

\begin{example} {\em The cone metrics $S_k$}.\label{Sk}
 
First let $k>0$, and 
let $M$ be the extended complex plane and $\om=z^{k-1}dz$. The
cone points of $M$ with cone metric $|\om|$ are $0$ and $\infty$ with cone 
angles $2\pi k$
and $-2\pi k$, respectively. From a constructive point of
view, for $k \in (0,1)$, 
the surface $S_k$ is the infinite sector of angle $2\pi k$, with
edges identified. When $k=1$, the surface 
$S_1$ is the extended complex
plane, equipped with the metric $dz$.
This surface has a single cone point of cone angle $-2\pi$ at
$\infty$. When $k>1$, consider that sector to be a multiple
covering of $\bfC-\{0\}$, with metric $|dz|$. If $z$ is the
variable in that plane, then $w=z^k$ and $\om=dw$. 
Note that $S_k$ and $S_{-k}$ are isometric with $z\to1/z$
producing the isometry.  Finally, for $k=0$, we set the metric
$\om=\frac{dz}{z}$, and $S_0$ is an infinite cylinder. 
\end{example}

Given two flat cone metrics, we can perform surgery to produce a
third one. Let $M_1$ and $M_2$ be cone metrics and let
$L_i\subset M_i, i=1,2$, be geodesics (straight lines) of the
same length that do not pass through (but may terminate at)
cone points. Join $M_1$ to $M_2$ along the $L_i$ by
identifying opposite edges of $M_1-L_1$ to $M_2-L_2$ in a
manner that produces a surface with orientation consistent
with the orientations of $M_1$ and $M_2$.

The 
end points of
the lines on the joined surface will, in general, be cone
points with cone angles equal to the sums of the angles at
the corresponding cone points of the $M_i$.

\begin{example} {\em Sewing an $S_k$ into a cone metric, $k>0$.}

Let $M$ be a cone metric and $L$ a straight line in $M$ of
infinite length, beginning at a point $p$ of positive cone
angle $2\pi k_p$ and terminating at a cone point, $q$, of
nonpositive cone angle $2\pi k_q$. Sew in $S_k$ to $M$ along
$L$ by matching the positive real axis in $S_k$ --- joining
$0$ to $\infty$ --- to $L$ with $0$ matched to $p$. The
resulting cone metric will have cone points at $p$ and $q$ of
cone angles $2\pi(k_p+k)$ and $2\pi(k_q-k)$, respectively (in
addition to any other cone points of $M$).
\end{example}

This is an example of the process of grafting of projective
structures (see e.g. \cite{Kl33} and \cite{Go87}.)

\begin{example} {\em Removing an $S_k$ from a cone metric. }

Let $M$ be a cone metric, $p\in M$ a point with cone
angle $2\pi k_p$ with $k_p\ge k>0$, and $L_1$ and $L'_1$ two
rays of infinite length in $M$ that satisfy the following
properties: $L_1$ makes an angle of $2\pi k$ with $L'_1$ at
$p$; the lines $L_1$ and $L'_1$ terminate at the same point
$q\in M$; the union $L_1\cup L'_1$ bounds a simply connected
region. Remove that region from $M$, with $L_1$ identifying
$L'_1$. The region removed is an $S_k$ (a fact that can be
seen easily or deduced from Proposition~\ref{uniquenessinf}
below). The resulting cone metric has cone points at $p$ and
$q$ of cone angles $2\pi(k_p-k)$ and $2\pi(k_q+k)$
respectively.
\end{example}
\subsubsection{Existence and uniqueness of cone metrics}

The next two propositions show that cone metrics are
essentially determined  by their cone points, and that the the
Gauss-Bonnet condition (\ref{GB}) and  its natural extension
(\ref{GBe}) are the only obstructions to existence. 
Here we aim to extend work of Troyanov \cite{Tr86} 
on cone metrics with positive and finite cone angles to 
the cases where the cone angles may be negative or
of simple exponential type.

\begin{proposition} \label{existenceinf} Let $M$ be a compact Riemann 
surface,
$\{p_1\dots p_r\dots,p_{r+\ell}\}$ a collection of distinct
points, $r>0$, $\ell\ge0$. Suppose $\{a_1\dots a_r\}$ is a
collection of real numbers satisfying (\ref{GBe})
\begin{equation}\label{GBe}
\sum^r_{j=1}a_j = -(2 - 2genus(M)) + r + 2\ell.
\end{equation}
Then there exists a cone metric on $M$ with finite cone
points $p_j$ with cone angles $a_j$, $j=1\dots r$ and
exponential cone points $p_k$, $0\le k\le\ell$ of simple type.
\end{proposition}

In order to state the uniqueness theorem for cone metrics we must introduce 
the following definition.
\begin{definition}
Two exponential cone points of simple type with local representations
$|e^{\f1w}\f{dw}{w^2}|$ and $|e^{\f1z}\f{dz}{z^2}|$ are {\em asymptotically
isometric} provided $\frac{dw}{dz}(0)=1$
\end{definition}

\begin{proposition} \label{uniquenessinf} A cone metric on a compact Riemann
surface with cone points with finite cone angles is determined
up to scaling by the location of the cone points and their
cone angles. The same result is true if one or more of the
cone points is an exponential cone point of simple type, provided that the
corresponding cone points are asymptotically isometric.

\end{proposition}

The proofs of these propositions are given in Appendix~B.

 The hypothesis of "asymptotically isometric" cone points
in Proposition~\ref{uniquenessinf} is necessary as the following example 
shows.
\begin{example} {\em Cone metrics with the
same cone points and cone angles are not necessarily scalar multiples
 of one another.}

Consider on $S^2=\BC\cup\{\infty\}$ the family of cone metrics given by
$$
\mu_\beta = \bigm|\f{z+1}{z-1}\f{z+2}{z-2}e^{\beta z}dz\bigm|
$$
$\beta>0$. All of the cone metrics $\mu_\beta$
have cone points at $\pm1$, $\pm2$,
and $\infty$ with cone angles $-4\pi$ at $z=1$ and  $z=2$, 
cone angles $+4\pi$ at $z=-1$
and $z=-2$, and an exponential cone point of simple type at $\infty$. From 
the definition of $\mu_\beta$ it is evident that the cone point at $\infty$
of $\mu_{\beta_1}$ is asymptotically isometric to the cone point at  
$\infty$ of $\mu_{\beta_2}$ if and only if $\beta_1=\beta_2$.
For any choice of
$\beta$, since $\mu_\beta|_z=\mu_\beta|_{\bar z}$, 
any segment of the real axis
is a geodesic. (In
fact, the real axis is a length-minimizing geodesic between the cone
points $\infty$ and $-2$ (and between 
$-2$ and $-1$) but we will not
need to use this observation). Even though the cone points and angles
are the same for all $\beta>0$, we will show that these metrics are not
all scalar multiples one of the other. Define
$$
f(\beta) = \int^{-2}_{-\infty}\mu_\beta\quad\text{and}\quad g(\beta) =
\int^{-1}_{-2}\mu_\beta.
$$
It is straightforward to show that
$$
0 < g(\beta)\le\f{e^{-\beta}-e^{-2\beta}}\beta,
$$
which implies that $\lim_{\beta\to0}g(\beta)\le1$. Since
$\lim_{t\to-\infty}\f{t+1}{t-1}\f{t+2}{t-2}=1$ and $e^{\beta t}>\f12$ for
$t>\f{-\log2}\beta$, it follows that $\lim_{\beta\to0}f(\beta)=\infty$.
Now suppose that the conclusion of Proposition 4 were true for the metrics
$\mu_\beta$. Then $$\mu_\beta=c(\beta)\mu_1$$ for some positive,
real-valued function $c(\beta)$. Moreover, since $f(\beta)$ is the
length of $\ov{(-\infty,-2)}$ in the $\mu_\beta$ metric and $g(\beta)$
is the length of $\ov{(-2,-1)}$ in the $\mu_\beta$ metric,
we would have
$$
\f{f(\beta)}{f(1)} = c_\beta = \f{g(\beta)}{g(1)}.
$$
But we have shown that $f(\beta)$ diverges as $\beta\to0$, and that
$g(\beta)$ is bounded as $\beta\to0$. Hence, it is not possible that all
the $\mu_\beta$ metrics agree up to a scalar stretch factor.

\end{example}
\subsection{The exponential cone as the limit of $S_k$ as
$k\goesto\infty$} \label{expconelimit}

In our construction of the surfaces \Hk/ in Section~4, we
will sew the cone  metrics $S_k$ into a torus and let
$k\rightarrow \infty$ with the expectation that the limit
corresponds to the creation of an exponential cone point of
simple type.  We will show here that the limit of $S_k$ as
${k\rightarrow\infty}$  (in an appropriate sense of limit) is
a cone metric on the sphere with  one exponential cone point
of simple type.

Let $z$ to be the variable on $S_k$ (considered as a
multisheeted sector with $|\th|\le k\pi$ in the
$z=re^{i\th}$-plane), and let $z=z(w)=(1+\f w k)^k$. The
metric $|dz|$  on $S_k$ with cone points at $0$ and $\infty$
is  isometric to the cone metric
$|(1+\f w k)^{k-1}dw|$ on $\bfC-\{-k\}$, whose cone points
at $w=-k$ and $w=\infty$ have cone angles $2\pi k$ and
$-2\pi k$, respectively. As $k\to\infty$, these metrics
on $\bfC-\{-k\}$ tend to
$|e^w dw|$ uniformly on compact subsets. 

We understand convergence of metric spaces here as relative
to a fixed base point; in this case, we take the origin as
the fixed point for each $S_k$.  Note that this point 
corresponds to $z=1$. The point $w=1$ correponds to a point
in  $S_k$ that is converging in the $|dz|$ metric to $e$.
Then the uniform convergence of the metrics on compacta, and the
choice of the origin as fixed for all $k$, implies that we
may regard the point $w=1\in \bfC$ as the limit of a bounded
set of points in $S_k$. 

Consider the annulus in $\bfC$ that is centered at the origin
and has inner radius $1+\f 12$ and outer radius $k-\f 1 2$.
This  annulus separates the pair of points $\{0,1\}$ from the
pair of points $\{-k,\infty\}$. The modulus of this annulus
(see Definition~\ref{modulus}) is equal to
$\log\f23({k-\f12})$, which goes to infinity with $k$.   From
this it  follows  from  Proposition~\ref{equivdef} in
Section~\ref{extlen} that the extremal length  of the class
of curves that separate these pairs of points goes to zero as
$k\rightarrow\infty$.  This shows that the points $-k$ and
$\infty$ coalesce as $k\rightarrow\infty$. What we mean by
this is explained in the next subsection in  
Remark~\ref{coalesce}. There are two types of limits of the
spaces $\{S_k\}$: the conformal limit of the punctured
Riemann surfaces $\bfC-\{-k,\infty\}$ and the (metric) limit
of the metric spaces $S_k$. The analysis above shows that the
conformal limit is $\bfC\union\{\infty\}$ with a distinguished
point at $\infty$, and the metric limit is $|e^wdw|$ uniformly on
compacta. This implies that, metrically, the spaces $S_k$ limit
on the sphere $\bfC\union\{\infty\}$ with a single cone point of
simple exponential type at $\infty$.

\begin{remark}
(i) Pulling back the metric on $S_k$ to the strip $|Im \zeta|<k$
by the map $z = \log\zeta$, one can consider the metric
$|e^{\zeta}d\zeta|$ on the strip to be a representation of
$S_k$ with cone points at $Re \zeta =-\infty$ and  
$Re\zeta =+\infty$ corresponding to the cone points $0$ and
$\infty$ in  the ``$z$ model'' of $S_k$. As
$k\rightarrow\infty$, the metric converges  to
$|e^{\zeta}d\zeta|$ on the entire complex plane. The argument using extremal 
length can be repeated here to show that
the strips converge to a cone metric with one cone
point at infinity.

(ii) We note here that there is a difference between the
exponential cone points and the cone points with finite cone
angle. In any cone metric, cone points  with finite {\em
positive} cone angles have neighborhoods where the metric  is 
precisely equivalent to the metric on a Euclidean cone, so
every curve from a regular  point to such a cone point has
finite length.  However,  every curve that goes  from a
regular point to a cone point with a finite, {\em
nonpositive} cone  angle  has infinite length. In particular,
a cone metric with all cone angles  finite defines a complete
metric space on the underlying Riemann surface  with the
nonpositive  cone points removed; each such point corresponds
to an end and the metric-space topology is identical to the
topology of the  underlying punctured Riemann surface. 

The
situation is not the same in the presence of cone points of 
simple exponential type. Consider the extended complex plane
with metric $|e^{\zeta}d\zeta|$, a cone metric with one cone
point of simple  exponential type at $\infty$. Horizontal
curves of the form $\{t+iy_0|\,  t<t_0\}$ have finite length
while those of the form $\{t+iy_0|\, t>t_0\}$ have infinite 
length. If the point at infinity is removed, the metric is
not complete. If it is left on the surface, then the metric
is complete but defines a topology that is not the  same as
the  topology of the Riemann sphere: for example, the sequence
of positive  integers  eventually leaves any neigborhood of
$\infty$.
\end{remark}

\subsection{Preliminaries on Extremal Length}\label{extlen}
We will make use of arguments using extremal
lengths, so we record the basics of this subject in this
subsection.

Extremal length assigns a conformal invariant
$\ext_\SR(\G)$ to a set of curves $\G$ on a Riemann surface
$\SR$. A flexible tool for distinguishing conformal
structures, extremal  length is especially useful when the
set of curves $\G$ is the free homotopy class of a simple
closed curve. For such a class of curves, there are two
equivalent definitions of the  extremal length $\ext(\G)$, with
one definition naturally suggesting lower bounds and the
other definition naturally suggesting upper bounds. 
We can often use this principle to obtain
good estimates (see e.g. \cite{Oht70})
for the asymptotics of $\ext_\SR(\G)$ in many situations under
which $\SR$ degenerates, sending $\ext_\SR(\G)$ to zero or
infinity.

We begin with the general definition.

\begin{definition} (Analytic) Let $\G$ be a set of curves on
a Riemann surface $\SR$. Then
$$
\ext_\SR(\G) =
\sup_\rho\f{\inf_{\g\in\G}[\ell_\rho(\g)]^2}{\area(\rho)},
$$
where the supremium is taken over measurable conformal
metrics $\rho|dz|^2$ on $\SR$, the notation
$\ell_\rho(\g)=\int_\g\sqrt\rho$ refers to the $\rho$-length
of a curve $\g\in\G$ on $\SR$, and
$\area(\rho)=\iint_\SR\rho$ is the $\rho$-area of $\SR$.
\end{definition}

\begin{example} Let $\SR$ be the annulus
$A(r_1,r_2)=\{r_1<|z|<r_2\}$ in the plane, and let $\G$
consist of all curves freely homotopic to the core curve
$\{|z|=\f{r_1+r_2}2\}\subset A(r_1,r_2)$. Then a simple
length-area argument \cite{Ah73} shows that
$\ext_\SR(\G)=\frac {1} {\log (r_2/r_1)}$.
\end{example}

\begin{definition}\label{modulus} The number $\log (r_2/r_1)$
is known as the modulus $\bmod A(r_1,r_2)$ of the annulus
$A(r_1,r_2)$.
\end{definition}

This leads to the second definition of extremal length in the
case that $\G$ is a free homotopy class of curves, all of
whose members are freely homotopic to a simple closed curve
on $\SR$.

\begin{definition} (Geometric) The extremal length
$\ext_\SR(\G)$ of a curve system $\G\subset\SR$ is defined to
be
$$
\ext_\SR(\G) = \inf_{A\subset\SR}\f1{\bmod A}
$$
where the infimum is taken over all conformal embeddings of
annuli $A$ into $\SR$ which take the core curve of $A$ into
$\G$.
\end{definition}

It is an important result (see \cite{St84}) that 

\begin {proposition}\label{equivdef} The geometric and
analytic definitions  of extremal length coincide.
\end{proposition}

Naturally, if we are interested in lower bounds for extremal
length, we compute
$\inf_{\g\in\G}[\ell_\rho(\g)]^2/\area(\rho)$ for a specific
conformal metric $\rho$ on $\SR$, and obtain a lower bound on
$\ext_\SR(\G)$.  On the other hand, if we are interested in
upper bounds, we compute the modulus of some specific annulus
embedded in $\SR$ with cone curve homotopic to $\G$, and
obtain an upper bound on $\ext_\SR(\G)$.

\begin {remark} \label {coalesce}In many of our applications,
we will wish  to show that a sequence of pairs of points, say
$\{p_n,p'_n\}$ ``coalesce'' to a single point $p_\infty$.
Conformally, this means that the set $\G$, of curves which
encircle $p_n$ and $p'_n$, have arbitrarily small extremal
length, i.e. the neck connecting a neighborhood of $p_n$ and
$p'_n$ is pinching off as $n\to\infty$. In this case, using
the geometric definition of extremal length, it is then
enough to show that there is a sequence of annuli $A_n$ with
$\bmod A_n\to\infty$ so that $A_n$ can be conformally mapped
into $\SR$ in a way that disconnects a disk containing $p_n$
and $p'_n$ from the rest of $\SR$.
\end{remark}

\subsection{The helicoid in terms of cone metrics}
\label{heliccone}
We conclude the background discussion by presenting the
helicoid from the point of view of cone metrics. From (9), we
have for the helicoid, \Hel/, modulo the screw motion $\s_k$:
$$
gdh = kz^{k-1}dz, \mbox{\,\,\,}\f1g dh=\f{-kdz}{z^{k+1}}\label{heli2}
$$
on $\bfC-\{0\}$ with ends at $0$ and $\infty$. We may consider
these forms to be defined on 
the cone mtrics $S_k$ and $S_{-k}$ respectively.
They are related by the  inversion  $z\rightarrow \f1z$, so
we may consider them both to be defined on the same domain.
This is  precisely the local expression of the form that 
produces the metrics on $S_k$ and $S_{-k}$ defined in
Example~\ref{Sk} in Section~\ref{conemetric}.

We may run this discussion backwards to construct
\He/$/\sigma_k$ from cone metrics. 
Both $S_k$ and $S_{-k}$ are defined on the extended plane.
We may develop (isometrically) both of these metrics onto the
Euclidean plane. If we pull back the naturally defined
one-form, say $dw$,  from that Euclidean plane to the
original extended plane, we obtain  
two one-forms we may use to define $gdh$ and
$\f1gdh$, respectively. Straightforward integration of (\ref{weier})
gives, as in Section~2.2, $x_1+ix_2=\bar z^{k}-z^{-k}$.
Moreover,
$$
dh^2=gdh\cd\f1gdh=-\f{k^2dz^2}{z^2},
$$
so $dh=\f{\pm ikdz}z$; thus we may recover the third component
of \He/$/\sigma_k$ from this data. Also
$$
g^2=gdh/(\f1gdh)=-z^{2k},
$$
so the Gauss map $g=\pm iz^k$  may also be recovered from this
data.

We showed in the previous section that $S_k$
converges, as $k\to\infty$, to a cone metric with a single exponential cone
point of simple type. Following the procedure there gives a
global representation of the limit of the cone metric
representation for \Hk//$\s_k$, yielding a representation of the
defining Weierstrass data of the helicoid by $gdh=e^zdz$ and
$\f1gdh=-e^{-z}dz$.

\subsection{Symmetries of the \Hk/}\label{Hk}
\FigSPGONEhalf

We conclude this section with a derivation of the symmetry
properties and the  conformal structure of the surfaces \Hk/
modulo $\sigma_k$. We assume   properties (i)-(iii) of (4),
which are the defining properties of the \Hk/,  whose
existence is asserted by Theorem~3. The results of this
subsection  are  collected in the lemma at the end of the
discussion.

From the assumptions (4),  we know that the periodic surface is
invariant under a vertical screw motion, $\sigma_k$, of angle $2\pi k$,
and that it contains a vertical axis.  By the Schwarz Reflection
Principle, which states that if a minimal surface contains  a
straight line then it is invariant under $180^\circ-$degree 
rotation  about that line, the surface is invariant under
$180^\circ-$degree rotation  about the vertical axis. In each
fundamental domain (a region of the surface that  generates the
whole surface by the action of $\sigma_k$), there are, by assumption 
two  parallel horizontal parallel lines.  It is easy to see
that---under the assumption  that the the  surface is singly
periodic but not doubly or triply periodic---the  horizontal
lines meet the vertical axis. 
To see this, recall that successive reflection
in two distinct  lines in $R^3$  results in a Euclidean motion
with the following properties: its  translational component is
in the direction of the segment of shortest length  between the
lines; its rotational component has this direction as axis
and rotation angle equal  to  twice the angle between the
lines. If either line does not intersect a third, say vertical,
axis, then the surface is invariant under a Euclidean motion
with a nonzero  horizontal  translational component. Since we
have assumed the existence of a vertical  translation $T$,
this contradicts our assumption of single-periodicity. 
In addition, if there were another non-horizontal line
on the surface, then the surface would be invariant 
under a screw motion along the line connecting the vertical
axis to that fourth line: thus the surface would 
be invariant under a non-vertical screw motion, which
is also a contradiction.

Let $d\neq 0$ be the distance between two consecutive
horizontal lines in  the surface. Clearly, $d<2\pi k$.
Consecutive reflection in these lines  results in a 
vertical screw motion $\sigma_d$. This screw motion is a
symmetry of the  surface and so must take horizontal lines
into horizontal lines. If $d<\pi k$, the image of one of the
lines under $\sigma_d$ is a third line not equivalent  to
either of the original two under the screw motion
$\sigma_k$, contradicting the assumption that the there are
only two lines in the quotient. If, on the other hand, we
have for some integer $m$, $2\pi k m>d\geq\pi k$, then either
$2\pi k m=2d<4\pi k$ or $2\pi km +d  =2d$, depending on
whether the $\sigma_d$ leaves the lines  invariant or
interchanges them, mod $\sigma_k$. In the first case, $m=1$ 
and $d=\pi k$, while the second case is not possible
because $2\pi km +d  =2d$ means that $d=2\pi km$, while
we required $d <2\pi k m$.
We
conclude that the two horizontal lines are separated by a
distance of $\pi k$, and that reflection in one of them is
the same as reflection in  the  other, mod $\sigma_k$.

Successive reflection in the vertical axis and a horizontal
line on the  surface that meets the axis at a point $p$,  
produces a reflection in the horizontal line through $p$ that
is   normal to the surface at $p$. Call this line $L$.
Reflection in $L$ is referred to as a {\em  normal symmetry}
of the surface. Modulo $\sigma _k$, this is the same 
symmetry no matter which horizontal line is chosen. 

Reflection in the vertical axis is an involution of the
minimal surface,   which leaves height unchanged and fixes
pointwise the vertical axis.  Its fixed-point set on the
quotient surface consists precisely  of those points that get
mapped to the vertical axis, a connected set. The  normal
symmetry around $L$ fixes a point on the vertical axis and
leaves  the vertical axis invariant. In the quotient surface,
the normal symmetry  fixes {\em two} points on the vertical
axis, the points where the  horizontal lines meet the
vertical axis: this follows because the normal
symmetry acts
as an isometry on the segment of the vertical axis 
in one fundamental domain of the surface.  Hence
as it inverts this segment about the point where the line
$L$ meets the axis, it fixes exactly that intersection
point and the (end)point at distance $2\pi k$ along
the segment. 

By hypothesis, the quotient surface is a torus, and we may
model it as the  region bounded by a parallelogram in the
complex plane, with opposite edges identified. (When we refer
to a ''parallelogram,'' we will  mean the closed region
bounded by a quadrilateral with opposite sides parallel, or, 
depending on the context, the Riemann surface of genus one
produced by identifying the opposite edges of  its boundary.)
The normal symmetry fixes two points on the vertical axis. 
Without loss of generality,  we will  assume that one of those
points correponds to the center of the parallelogram.  We
label the center by \O/.

On the torus, there is a holomorphic one-form with no zeros.
Up to scaling,  this one-form is equal to the one-form that
descends to the  parallelogram  from the one-form $dz$ on the
complex plane. We will also refer to this one  form as $dz$.
Let $\rho$ be the involution on the quotient surface that is 
induced by the normal symmetry. Let  $r$ be 
$180^\circ-$rotation about the  center of the parallelogram, also an
involution of the torus. Both $\rho^*(dz)$ and $r^*(dz)$ are 
zero-free holomorphic one-forms  and  so must agree up to a
scalar multiplicative factor. Since they are  involutions
that fix the  center point, $r^*(dz)=-dz=\rho^*(dz)$ at that
point. Hence $r^*(dz)= \rho^*(dz)$ everywhere, which  implies
that $\rho=r$. 

Let $\mu_v$ be the (anticonformal) involution of the
parallelogram corresponding to the involution of the
quotient surface produced by  reflection in the vertical
axis. We know that the fixed-point set of $\mu_v$ is a connected
curve that passes  through the center of the parallelogram. We
look at the action of $\mu_v$  near  the center point of the
parallelogram.  If $W$ is the line through the center point of
the parallelogram that is  tangent to the fixed-point set of
$\mu_v$, let $r_W$ be reflection across $W$  in the complex
plane.  At the center, 
$\mu_v^*(dz)=\lambda(d\overline{z})= r_W^*(dz)$, for some
complex number $\lambda$, $|\lambda|=1$. As in  the previous
paragraph, we can conclude $\mu_v^*(dz)= r_W^*(dz)$ and hence
that $\mu_v=r_W$ near the center,  hence everywhere. In
particular, we have shown that reflection in the line $W$ is
an involution of the parallelogram. This implies that the 
parallelogram is either a rhombus with a diagonal on $W$, or 
a rectangle with a side parallel to $W$. But we know that the
fixed-point  set of $\mu_v$ is connected, which implies that
the parallelogram is a  rhombus.

Rotate the rhombus if necessary so that the diagonal on $W$
is  a vertical line segment. We will refer to this
diagonal as the vertical  diagonal. Reflection in the other
diagonal of the rhombus is equal to 
$\mu_h:=\mu_v\circ\rho$, and so must correspond to the
reflection in the  horizontal lines of the quotient surface.
Since this diagonal represents  both lines and since the
lines diverge, the punctures at the ends must  appear on this
diagonal. We will denote the end punctures by $E_1$ and $E_2$. Because $r$ 
leaves $\{E_1, E_2 \}$ invariant---in fact, it interchanges
$E_1$ and $E_2$---they  must be symmetrically placed with respect to
\O/. 

The careful reader will note that we have specified  \O/ to
correspond to  one of the two points where a horizontal line 
on the surface meets the  vertical axis. The rotation $r$ by
$180^\circ-$degrees about \O/ fixes four points; the center,
the vertex,  and  the two half-periods. Since $r$ corresponds
to the normal  symmetry, and $r\circ \mu$ corresponds to
reflection in the horizontal  lines, the two
off-vertical-axis  fixed points of $r$ must lie on the same
horizontal line. {\em Without loss  of generality, we may 
assume that it is this horizontal line that crosses the axis
at  \O/. } We collect the above discussion as

\begin{lemma} \label{symmetry} The defining properties (4)(i)-(iii) of the 
surfaces \Hk/
of Theorem~3 imply that 
\begin{enumerate}
\item[(iv)] The horizontal lines on \Hk/  meet the vertical axis. Two
successive horizontal lines are separated by a vertical distance of 
$\pi k$. Composition of rotation  about  two successive lines is a
vertical screw motion $\sigma_k$,  which descends to the  identity 
transformation on \Hk// $\sigma_k$.  Rotation by $180^\circ$ about one
of the horizontal lines is a symmetry of the surface that descends to
an involution of quotient surface, and that involution does not depend
on the choice of horizontal line.;
\item[(v)] Rotation by $180^\circ$ about the vertical axis is a symmetry
of \Hk/.  The composition of this rotation with rotation about a 
horizontal line on the surface is rotation by $180^\circ$ around a
line  orthogonal to the axis and the horizontal line.  These order-two 
symmetries---referred to as a normal symmetries---induce the same
involution on the quotient surface of \Hk//$\sigma_k$.
\item [(vi)] The quotient surface of \Hk//$\sigma_k$ has the conformal
structure of a rhombic torus with two punctures. Without loss of
generality, we may assume this rhombus: is conformally modelled by the
domain bounded  by a rhombus in the plane with opposite edges
identified; and  is oriented so that the one of the diagonals is
vertical, the other horizontal; the  vertical diagonal is mapped into
the vertical axis; the horizontal diagonal  is mapped onto the
horizontal lines of \Hone/. In particular, the two 
punctures---corresponding to ends---occur on the horizontal
line, and they  are symmetrically placed with respect to the origin.
They separate it into  segments mapped to the two different horizontal
lines.
\item [(vii)] On the rhombus, the reflection $\mu_v$ in the vertical
axis  corresponds to $180^\circ$-rotation about the axis of \Hk/. The 
reflection $\mu_h$ in the horizontal diagonal of corresponds to
$180^\circ-$rotation around a horizontal line of \Hk/. The
$180^\circ$-rotation, $r=\mu_v\mu_h$, about the center, \O/, of the
rhombus corresponds to the  normal symmetry of the quotient surface.
Two of the fixed points of $r$, namely \O/ and the vertex of this
rhombus, correspond to the two points in  the quotient of
\Hk//$\sigma_k$ where the horizontal lines cross the vertical axis. The
two other fixed points of $r$ lie at the half-period points and
correspond to two off-axis fixed points of the normal symmetry. These 
points lie at the same height as one of the on-axis fixed points:
without loss of generality, we may assume that they lie at the same
height as the  point correponding to \O/.
\end{enumerate}
\end {lemma}

\section{The singly periodic genus one helicoid}

We give a new proof of the existence of the surface \Hone/ of 
Theorem~\ref{PG1H} in Section 1.3.  The alternative
construction  presented here is due to Weber, and it is
key the  proof of the existence of the  \Hk/ of
Theorem~\ref{MainThm} in  Sections 4 and 5. It is in this
section that the transformation of the  problem from an
analysis of forms and functions on Riemann surfaces to 
geometric manipulation of singular flat structures is most
clearly  displayed. The existence and embeddedness  of
\Hone/  was  originally  proved  in \cite{howeka4}.

 In 3.1, we  derive 
necessary conditions for the Weierstrass data of \Hone/ that are 
sufficient to specify them uniquely. Together  with a real
parameter that gives us the underlying rhombic torus, this 
determines not only the Riemann surface structure of the
quotient surface  but also the Weierstrass data $g$ (up to a unitary 
multiplicative factor) and $dh$ (up to a positive real multiplicative factor).
(Geometrically, the surface is determined up to a rotation about the 
vertical 
axis and a scaling in $R^3$. If we  make
the reasonable assumption that the horizontal lines on the
surface are  parallel to the $x_2-$axis, and that the
translational symmetry of the  surface is generated by a
vertical  translation of length 
$2\pi$, then the Weierstrass data is completely determined.)

In Section 3.2,  we formulate the period problem for \Hone/ and
state the main results of \cite{howeka4}, namely that there
exists a choice of parameters that solve the period problem. Any such choice
produces a minimal surface that has all the required properties
of  Theorem~\ref{PG1H}  in Section 1.3, including the property
of embeddedness.
It is necessary to prove the existence of \Hone/ using the methods here 
because of the relevance of the methods to the construction of the surfaces 
\Hk/ in Section~4.

In Sections 3.3--3.4 we present an alternative proof of the
results of  Theorem~\ref{PG1H} (excluding embeddedness). We
begin in Section~3.3  by  showing that the properties of the
Weierstrass data derived in Section 3.1  uniquely determine
the rhombic torus. This means that there is in fact only one
parameter in the Weierstrass data: the position of the ends.
Next we construct this special torus, which we call \T/,
by    cone-metric methods. This method actually determines
not only the torus but  also the one-form $gdh$. For any
placement of the ends on \T/, we use a  symmetry construction
to produce a  candidate for $(1/g)dh$ for which the 
horizontal period problem is {\em automatically solved}. In
Section 3.4, we  then show that there is a placement of the
ends for which the vertical  period problem is also solved.

\subsection {The Weierstrass data: derivation from geometric
assumptions} We want a singly periodic, properly immersed,
minimal surface with the properties i)-iii) of
Theorem~\ref{PG1H} in Section 1.3. Namely, \Hone/ modulo
translations has the properties stated in (3):

\begin {enumerate} 
\item [(i)] \Hone/ has genus one and two ends,
\item [(ii)] \Hone/ is asymptotic to a full $2\pi$-turn of a
helicoid, and
\item [(iii)] \Hone/$\sigma_1$ contains a vertical axis and two horizontal
parallel lines.
\end {enumerate}
See Figure~\ref{figure:help1} for an image of this surface.

We will now assume that such a surface exists and derive its
Weierstrass  data.

\subsubsection{The Gauss map and the placement of the ends}
\label{abplacement}
The total curvature of the quotient surface, $M =$\Hone/$/sigma_1$, is
$2\pi(\chi(M)-W(M))$, where $\chi(M)$ is the Euler
characteristic of $M$ and $W(M)$ is the total winding number
at the punctures. \cite{jm1} We know that $M$ is a
twice-punctured torus so $\chi(M)=-2$; further, each end is,
by  assumption, asymptotic to a single full turn of the
helicoid, so $W(M)=2$. Hence, the total curvature of $M$ is
equal to $-8\pi$. This implies that the  the degree of the
Gauss map is equal to two. 

The assumption (iii) that each end is asymptotic in the
quotient to a  single full turn of the helicoid forces us to
assume, according to (\ref{helicgdh}),  that the Gauss map is
vertical at the ends. Let $g$ be the stereographic 
projection of the Gauss map, defined on the underlying
Riemann surface.  Since the degree of the  Gauss map is two,
there must be one other point where $g=0$ and one other 
point where $g=\infty$. According to Lemma~\ref{symmetry},(vi)and (vii),
the underlying  Riemann surface may be modelled by a rhombic
domain, and reflection in the diagonals of that rhombus 
correspond to
rotation by $180^\circ$ about the lines on \Hone/.  These 
rotations are orientation-reversing and preserve verticality.
Therefore, the reflections in the diagonal preserve
verticality and so must leave the collection of poles
and zeros of $g$ invariant. If $g$ has a zero  or pole that
is not on  the diagonals of the rhombus, then it has {\em at
least} two zeros and  two poles not on the diagonals, in
addition to the two ends that are located on the horizontal
diagonal.  This implies that the degree of $g$ is at  least
$3$, contradicting the fact that the degree of $g$ is two.
Hence the  other points  where $g=0$ or $g=\infty$ lie on the
diagonals. Since the vertical diagonal  corresponds to the
vertical axis, we require the Gauss map to be horizontal
there: $|g|=1$ on the vertical diagonal.  We conclude that 
the other two vertical points of the Gauss map lie on the
horizontal  diagonal. 

Label the ends $E_1$ and $E_2$ and the vertical points $V_1$
and $V_2$. Since reflection in the vertical diagonal
corresponds to the  orientation-reversing symmetry of
rotation about the  vertical axis, a symmetry that preserves
verticality of the Gauss map, the  pair of ends and the pair
of vertical points must be symmetrically placed  with respect
to the center, \O/, of the rhombus.

For the remainder of Section~3, we will assume that the the
length of the horizontal diagonal is two, that  \O/ is placed
at the origin of ${\bf C}$, and that $E_1$ lies to the left
of \O/. Then we may write
$$
E_1=-b, \mbox{     \,\,} E_2=b,
$$ 
$$
\mbox{     \,}V_1=-a,  \mbox{    \, \,\,} V_2=a,
$$
for some real numbers, $0<a$, $b<1$,  The strict inequality is
required in order to have two ends and two vertical points.
Also, we require $a\neq b$ to prevent the ends from 
coinciding with the finite points.  Without loss of
generality, we may also  assume that the surface is oriented
so that $g(E_1)=\infty$ (and therefore $g(E_2)=0$.)  

Either $g(V_1)=\infty$ and $g(V_2)=0$ or $g(V_1)=0$ and
$g(V_2)=\infty$. Applying Abel's theorem to $g$ yields
$|a-b|=1$ in the first case and $a+b=1$ in the second case.
The first case is impossible because $0<a,b<1$.  Therefore
$g(V_1)=\infty$ and $g(V_2)=0$ and 
$$
a+b=1.
$$
We have now determined the divisor of $g$. (See Figure~8).
The points $E_i$, $V_i$ are symmetrically placed with respect
to the quarter  points of the horizontal diagonal. (For  {\em
any} degree-two elliptic function, the branch points are
symmetrically  placed with respect to the zeros and poles of
the function.) Also, the branch points must be symmetric with
respect to the symmetries of the surface.  We conclude from
this that the quarter points of the horizontal diagonal  are 
branch points of $g$  and that there are two symmetrically
placed branch points of $g$ on the vertical diagonal: in fact
they are also the quarter points.

We will assume that the surface is rotated so that its normal
vector at \O/ is $(1,0,0)$: that is, $g({\cal O})=1$. This
rotation makes the  horizontal lines parallel to the
$x_2-$axis.

\begin{remark}
\label{a<b?}
At this point, our specification of Weierstrass data depends
on the choice  of rhombic torus and a choice of $b$ between
$0$ and $1$ to place the ends $E_i$. We will see in
Section~3.1.2 that $dh$ is determined by these  choices. 

We note that we have not specified whether or not $a<b$, i.e.
whether or  not the vertical points lie on the line that
passes through the image of \O/. We do not at this point have
the freedom to assume one way or the other. It turns out
that, in fact, $a<b$,  and the $V_i$ lie closer to \O/
than  do the $E_i $. This is the result of a  computation
(see Proposition~\ref{oldsinglyperiodic}, Statement~3) that 
shows that the period problem {\em cannot be solved if
$b>a$}. 
\end{remark}

\FigDivisorsOne

\subsubsection{The one-form dh}

The expression (\ref{metric}), for the induced metric on $M$
show that $dh$  has simple zeros at the vertical points $V_i$.
Because we want the  ends to be helicoidal, we require $dh$
to have simple poles at the the ends $E_i$. (See
Section~\ref{helicoiddef} and Figure~\ref{figure:divisors0}.) The one-form 
$dh$ can have no other poles or 
zeros. Hence we know the divisor of $dh$. Because we require
the horizontal  diagonal to be mapped into horizontal lines,
$dh$ must be purely imaginary on the horizontal diagonal. This
determines $dh$ up to a real scalar factor,  which
corresponds to scaling the surface in $\BR^3$.

\FigPieceHeOne 

\subsection{The period problem}
\FigHOneRhombus
In Section~3.1, we  specified a  Weierstrass representation
for \Hone/ which was forced by the geometric conditions 
4(i)--(iii) of Theorem 2 in Section 1.3. The underlying torus
is rhombic and the ends are placed on one of the  diagonals.
The divisors of $g$ and $dh$ are determined. The situation is
encapsulated in  Figure~\ref{figure:divisors1}. The function $g$ is then
completely determined by the condition that $g=1$  at the
center of the rhombus, and the one-form $dh$ is determined up
to a real scaling  (We could of course determine that factor
by the condition imposed by Theorem 2 and (\ref{weier}), that
$Re\int{dh}=2\pi$ on  the vertical diagonal, but we do not do
that at this time.)

We still have to impose the period conditions (\ref{hpc}) and
(\ref{vpc}).  In Figure~10, the indicated cycles $B$ and
$\beta$, together with their  reflections in the vertical
axis, generate a homology basis for $M$. The  nonzero
translational period we require to produce a singly periodic 
surface  will be evident on the cycle $\beta$ that surrounds
an end. By the symmetry we have imposed upon the Weierstrass
representation, we need only consider  these two cycles. The
horizontal and vertical period condition for \Hone/  can be
written as follows:
\begin{equation}
\int_B gdh = \overline{\int_B {\frac{1}{g}dh}} \label{hpc2}
\end{equation}
\begin{equation}
 \Re \int_B dh = 0. \label{vpc2}
\end {equation}

\subsubsection{Existence of \Hone/}

The following result, translated into the terminology used
here, is proved in \cite{howeka4}.
 
\begin{proposition} \cite{howeka4} \label{oldsinglyperiodic}
\begin{enumerate}
\item[] Every regular, complete, periodic minimal surface
satisfying the conditions of Theorem \ref{PG1H} (stated in 
(3)) may by
represented by Weierstrass data that satisfies all of the
conditions described in the table of Figure~\ref{figure:divisors1}
in Section~3.1.
Conversely,
\item Given any rhombic torus and any $b\in (0,1)$, there
exists Weierstrass data $\{g, dh\}$ on a rhombic torus with
ends and vertical points determined by the choice of $b$ as
in Section~\ref{abplacement}, whose divisors are specified in the table
of Figure~\ref{figure:divisors1}. The Weierstrass integral
(\ref{weier}) produces a {\em multivalued}, regular, minimal
and complete immersion of this punctured torus into $R^3$
with all the required symmetry properties of 
Lemma~\ref{symmetry};

\item The immersion described in 1) above will be singly
periodic if and  only if the Weierstrass data satisfy
(\ref{hpc2}) and (\ref{vpc2}). The translational period $T$
is determined the vertical vector $(0,0,c)$ where
$$
\pm c = \int_\beta dh =  2\pi i Res_{E_i}dh;
$$
\item The period conditions (\ref{hpc2}) and (\ref{vpc2})
cannot be  satisfied unless $\frac{1}{2}< b <1$. In
particular, $1-b=a<b$. Thus it is necessary that the vertical
points, $V_i$, be located (as illustrated in Figure~9) closer
to the center, \O/, than the end points $E_i$. 
\end{enumerate}
\end{proposition}

This proposition is proved in Section~1 of \cite{howeka4}: see page
259. How  Theorem \ref{PG1H} is proved from
Proposition~\ref{oldsinglyperiodic} we outline in the
following  remark. 

\begin{remark}\label{verticalperiodfunction} One can
establish the  existence of \Hone/ in two stages.
First, it is shown that there is an open interval  of 
rhombic tori (parametrized by branch values---say $\rho$---of
a geometrically-normalized {\cal P}-function) for which the
vertical period problem (\ref{vpc2}) can be  solved by an
appropriate choice of $b$. That choice of $b$ is unique, 
satisfies $b>1/2$ and depends smoothly on $\rho$. For  values
of $\rho$ outside that interval, (\ref{vpc2}) cannot be
solved.   Second, the horizontal period condition 
(\ref{hpc2}) changes sign on the curve $(\rho, b(\rho))$, and
hence is zero for at least one value of $\rho$. This is
proved in Section~2 of \cite{howeka4}. Embeddedness of \Hone/ is
proved by a separate argument, as is almost always the case
in these matters.
\end{remark}
 
In the next section, it will be shown that the form of the
Weierstrass data  and the horizontal period condition {\em
determine the conformal structure}. Putting that  together 
with the Remark above will show that there the singly periodic
genus-one helicoid is  unique. See Proposition~\ref{rmk2}.

\subsubsection{The rhombic torus \T/ and the uniqueness of \Hone/}

\vspace{.2in}
We assume that the Weierstrass data for \Hone/ has divisors
described in  Figure~8 and that the zeros and poles are
symmetrically placed along a  diagonal of a rhombus. Note
that the divisors of $gdh$ and $\frac{1}{g}dh$ each have one
double zero and one double pole (with no residue) on a
diagonal,  and no other zeros or poles. Since the torus is a
group, it follows that $gdh$ and $\frac{1}{g}dh$ differ by a
translation and scaling:
$$gdh = ct^*(\frac{1}{g}dh),$$
for some translation, $t$, of the torus and nonzero
constant $c$. Then, for any closed curve $\a$ on the 
torus, we have
\begin{equation}
\int_{\alpha}{gdh}=
c\int_{\alpha}t^{*}({\frac{1}{g}dh})
=c\int_{t^{-1}(\alpha)}{\frac{1}{g}dh}
=c\int_{\alpha}{\frac{1}{g}dh}, \label{tstareqn}
\end{equation}
the last equality following because the form $\frac{1}{g}dh$ 
has no residue. If
$\{\alpha_1,\alpha_2\}$ is a basis for the homology of the
torus, then
$$
r:=\frac{\int_{\alpha_1}{gdh}}{\int_{\alpha_2}{gdh}}
=\frac{\int_{\alpha_1}{\frac{1}{g}dh}}{\int_{\alpha_2}{\frac{1}{g}dh}}.
$$
The horizontal period condition (\ref{hpc})
$$
\int_{\alpha_i}{gdh}=\overline{\int_{\alpha_i}{\frac{1}{g}dh}}
$$
now implies that $r=\bar r$, i.e. $r$ is real. If we modify
$gdh$ by multiplication by a constant---if necessary---to make
$\int_{\alpha_1}{gdh}$ real, then $\int_{\alpha_2}gdh$ must
also be real.

Weber \cite{weber1} notes that this simple condition 
characterizes the underlying rhombic torus.

\begin{proposition} \label{uniquetorus}\cite{weber1} There
exists a unique rhombic torus carrying a one-form $\eta$ with
the following properties: the form $\eta$ has a double zero and a
double pole (with no residue) on a diagonal, no other poles
or zeros, and all periods real.
\end{proposition}

\begin{definition}We will use the symbol \T/ to refer to the rhombic torus of
Proposition~\ref{uniquetorus}. \label{T1}
\end{definition}

Proposition~\ref{uniquetorus} can be used to prove
(see also the recent preprint \cite{} by Mart\'in)

\begin{proposition}\label{rmk2} \cite{Ka03} 
\Hone/ is unique.
\end{proposition}

\begin{proof}  
In Lemma~\ref{symmetry}, we showed that the geometric
conditions of Theorem 2 in Section 1.3 implied that the
quotient surface of any \Hone/  had to be a rhombic torus 
upon which $gdh$ satisfied the conditions on $\eta$ in 
Proposition~\ref{uniquetorus}. That is, the geometric data
that solve the horizontal period condition (\ref{hpc}) or
(\ref{hpc2}) specify the rhombus uniquely. Yet, as noted in  
Remark~\ref{verticalperiodfunction}, it is shown in 
\cite{howeka4} that
for each admissible
rhombus there is a  unique placement of ends (i.e. a unique
$gdh$) satisfying the vertical  period condition (\ref{vpc})
or (\ref{vpc2})).  As
Proposition~\ref{uniquetorus} implies that
there is a unique rhombic torus, \T/, on which the horizontal
period condition (\ref{hpc2}) can be  satisfied,
there is then at most one surface satisfying the  conditions of
Theorem~2 in Section~1.3: the surface \Hone/ is unique.
\end{proof}

\begin{remark} \label{thetavalue}
The rhombic torus \T/ is the torus ${\bf C}$ modulo 
the lattice generated by $\{1,
e^{i\theta_1}\}$, where $\theta_ 1\approx 1.7205$.
\cite{howe5} \end{remark}

\subsection{The cone metric construction of \Hone/}

In Section~3.1, and in Lemma~\ref{symmetry}, various
conditions for the  Weierstrass data of \Hone/  were derived
using geometric and analytic  arguments. In Section~3.2, we
stated the period problem for \Hone/. In this  section, we
 give a cone-metric derivation of the Weierstasss data, the most important 
feature of which is that it solves the horizontal period problem 
(\ref{hpc2}) 
by construction.  It  does not use
Proposition~\ref{oldsinglyperiodic} of Section~3.2, which 
depends on the analysis and estimates of \cite{howeka4}. 

\subsubsection {The cone-metric construction of \T/ and
Weierstrass data  for \Hone/}
\FigslitTone
We will construct the torus \T/ in a manner that produces,
at the same time, a candidate for the one-form $gdh$. 

A torus may be constructed by identification of opposite
edges of a parallelogram. Consider the region bounded by a
parallelogram with  vertices $0$, $1$, $\tau$ and $1+\tau$ in
${\bf C}$. The one-form $dz$ on $\bfC$ induces a  holomorphic
one-form on the torus. For the  cycles on this torus that
correspond to the edges $\overline{0,1}$ and 
$\overline{0,\tau}$, the periods of this one-form are clearly
visible in the construction; they are the complex numbers $1$
and $\tau$.

Instead of the region bounded by the parallelogram, consider its
complement in  the extended $z$-plane. Identify opposite
boundary edges; again, this is  topologically a torus. We
have flat charts given by $\int dz$ away from the vertex
point. On this torus, with the  induced flat metric from the
plane, the vertex is a cone point with cone  angle $6\pi$.
(See Section~\ref{conemetric} for a discussion of cone points and cone 
metrics.) Allowing 
a  slight abuse of notation,  we will write $dz$ for the
induced  one-form  and $|dz|$ for the  associated metric on
the torus. In this language then, the form $dz$ must have
a double zero at the vertex point and, of  course, it has a
double pole at infinity. The periods of $dz$ along the  cycles
corresponding to the edge vectors of the parallelogram are
$1$ and $\tau$. 

Recall that we desire to produce \T/, a rhombic torus 
carrying a one-form with a double pole, a double zero and
real periods as in Proposition~\ref{uniquetorus}. If we construct the torus 
from the point-of-view of the 
previous paragraph, we are forced to choose both the edge
vectors of our parallelogram to be real. Even though such a 
``parallelogram" is degenerate, the construction produces a
topological torus with a flat  structure given by  $|dz|$,
which has isolated conical singularities.  Therefore it has a
regular conformal structure. 

Essentially, we are slitting the plane (without loss of
generality, along $[-1,1]$), choosing $c\in[0,1)$, and
connecting the region above (below) $[-1,-c]$ with the region
below (above) $[c,1]$ by identifying these boundary segments.
Only for the  choice of $c=0$ will the torus be rhombic;
reflection in the imaginary axis provides an involution with
a  connected fixed point set. The rhombic torus so
constructed carries a one-form descended from $dz$, with one
double pole and one double zero on a diagonal, no  other
zeros or poles, and all periods of this one-form are real.
Because  the sum of the residues of any one-form is zero and
there is only one pole, this form 
$dz$ has no residue at its pole. This torus 
satisfies 
the requirements of Proposition~\ref{uniquetorus}, 
and therefore must be the unique \T/. 

We will take $dz$ as a candidate for $gdh$ on \T/, and we will refer
to this  presentation of \T/ as the {\em slit model}. 

\begin{remark}  \label{specificrhombus}In the construction of
tori by removal of the interior of a parallelogram, $P$, from
the plane, it is important to note that the torus constructed is
not, in general, conformal to the torus produced by taking the
interior of $P$ and identifying opposite sides. This is clear in
the limit case of \T/, constructed by removal of a slit from the
extended plane; the torus \T/ is a nondegenerate rhombic torus
conformal to the one produced by $\bfC/\{1,e^{i\theta_1}\}$ with, according 
to 
Remark~\ref{thetavalue}, 
$\theta_1\sim 1.7205$. It is true, however, that removing the
square from the extended plane produces the square torus; the
torus produced is both rhombic and rectangular and such a torus
must be the square torus.
\end{remark}
\Figslitandrhombus

We now relate the slit model of \T/  to the rhombic model,
which we will take, as in Section~3.1, to be  the region bounded
by a rhombus whose diagonals are parallel to the  coordinate
axes.  (See Figure~\ref{figure:slitandrhombus}, left.) We will 
scale the rhombus so that the length of the horizontal diagonal
is equal to  two.  The slit model and the rhombic model are
conformally diffeomorphic.  We choose a conformal map from the
slit model to the rhombus---in terms of 
Figure~\ref{figure:slitandrhombus}, a mapping from the domain on
the right  to the domain on the left---that takes the imaginary
axis of the slit model onto the horizontal  diagonal of the
rhombus: both of these lines are  symmetry lines of the conformal
structure. We will do this in such a manner  as to have the
standard vertical orientation of the imaginary axis in the  slit
model correspond to the standard left-to-right orientation of
the  horizontal diagonal. The conformal diffeomorphism between
the slit model  and the rhombic model is now completely
determined up to the action of a  translation on the rhombus,
which with our normalization must be a  horizontal translation.
Therefore, the conformal diffeomorphism is determined up to
composition on  the left by a horizontal translation.

We label by
$v_2$ the vertex point in  the slit model and refer to the point
at infinity in the slit model as $e_1$.
We want the pullback of $d\z$ from the slit model to correspond to
the one-form $gdh$ on the rhombic model. This  pullback will
have a double pole at the inverse image of the point at infinity
and a double zero at the inverse image of the vertex point in
the  slit model. Since the double pole of $gdh$ occurs at the at the
end where $g=\infty$ (that is at $E_1$) and the double zero of $gdh$
must  occur at $V_2$, we see that
the inverse image of $v_2$ is $V_2$ and the inverse 
image of $e_1$ is $E_1$.  In Section~3.1, we saw
that we could assume without loss of generality that $E_1$ lies
to the left of the center of the rhombus.  Without loss of
generality we restrict ourselves to conformal diffeomorphisms
with this property. 

As in Section~3.1, we define $E_2$ and $V_1$ to be the points
on the  horizontal diagonal symmetric---with respect to
\O/--- to $E_1$ and $V_2$,  respectively. We may again write
$$
E_1=-b,\,\, E_2=b,$$ 
$$
V_1=-a,\,\,  V_2=a
$$
for some real numbers $a$, $b$ with $0<a$, $b<1$, $a\neq b$.
Abel's theorem  applied to $gdh$ requires $a+b=1$. We do {\em
not} know at this point  whether or not we may assume that
$a<b$. This does follow from the estimate of \cite{howeka4} discussed in
Remark~\ref{verticalperiodfunction} but as we shall see, we can achieve our
result without this external reference. 

We will adopt the convention that points on the slit model that are labelled 
by
 letters in lower case have corresponding points on the rhombus labelled by 
the same letter in  upper case.

\begin{remark} \label{degeneration}
As discussed in Section~3.1, it does not make geometric sense
to allow $a$  or $b$ to take on the values $0$, $1/2$ or $1$.
There are  three  possibilities: $a=0$ and $b=1$; $a=1$ and
$b=0$; $a=b=1/2$. In the first two cases, the  ends coincide,
while in the third case each end coincides with a  vertical
point on the surface; in all three cases, the data are 
incompatible  with our geometric assumptions.  {\em However, it
does make analytic sense to  allow this to happen.} The pullback
of $dz$ will define a one-form on the  rhombus with one double
pole and one double zero. In the cases where $a=0$ or $a=1$, the
pole occurs at the center or at the vertex point. In the case 
where $a=1/2$, the pole and zero occur at the half-period points
on the  horizontal diagonal.
\end{remark}

\subsubsection{The definition of $\frac1gdh$ and the solution
of the  horizontal period \label{1overgdh}
problem}
\label{horizontalperiod}
Let ${\bf t}$ be the translation of the torus in the rhombic
model that is induced by the translation in the plane
satisfying ${\bf t}(E_1)=E_2$  (and therefore ${\bf
t}(V_1)=V_2$). In terms of our normalization, ${\bf t}$ is  the
translation on the torus induced by the translation in ${\bf
C}$ by $2b$, where $E_1=-b$. (Note that by Abel's Theorem,
this translation $2b$ can be written as $2b=2(1-a)=2-2a$,
and so an equivalent translation is by $-2a$.  We will
use this observation in the proof of Proposition~\ref{findgooda}.)
We define
$$
\frac1gdh:={\bf t}^*(gdh).
$$
The one-forms $gdh$ and $(1/g) dh$ have the divisors as
specified in  Figure~\ref{figure:divisors1}, and {\em they
automatically satisfy the  horizontal period condition
(\ref{hpc2}).} To see this, recall that the periods of $gdh$
are  real by construction and, by \eqref{tstareqn},
$$
\int_{\a_i}\frac1g dh = \int_{\a_i}{\bf t}^*(gdh) =
\int_{{\bf t^{-1}}(\a_i)}gdh =
\int_{\a_i}gdh.
$$

It is important to observe that while we have one free parameter
to define $gdh$---essentially the parameter $a$ defined in
Section~3.1 that places the point $V_2$---the
choice of $\frac1gdh$ is  determined by the horizontal period
condition (\ref{hpc2}). We now seek a  choice of $a$ for which
the vertical period condition (\ref{vpc2}) is satisfied. This
condition requires us to integrate $dh$ over a specified  cycle
on the rhombus. To find a value of $a$ that satisfies
(\ref{vpc2}) we first need to define $dh$ in terms of $gdh$ and
$\frac1gdh$ and to verify that it has the desired symmetries.

\subsubsection{The Weierstrass data \{g,dh\} and the
symmetries of $dh$} 

Since $dh^2=gdh\cdot\frac1gdh$ and $g^2=\frac{gdh}{\frac1gdh}$,
we can use these forms to determine $g$ and $dh$ ---up to
sign--- by taking a  square root. The one-form $dh$ defined in
this manner will have simple  zeros at the $V_i$ and simple
poles at the $E_i$, precisely what is required by the divisor diagram in 
Figure~\ref{figure:divisors1}.

To prove the existence of \Hone/ we must find a choice of the
value $a$ so  that the one-forms $gdh$ and $\frac1gdh$
determined by this choice yield a candidate
$dh=\pm\sqrt{gdh\cdot\frac1gdh}$ that satisfies the vertical
period condition (\ref{vpc2}). (We discuss below how the sign
is determined.) We solve this problem in the
Section~\ref{solvevpc}. First, we establish some expected but
important properties of any $dh$ defined in the manner just
described.

\begin{lemma} \label{dhlemma} (i) The one-form  $dh$, considered
as a one-form  on the rhombus, is imaginary on the horizontal
diagonal and real on the  vertical diagonal.

(ii) Let $\mu_v$ and $\mu_h$ be reflection in the vertical diagonal
and the horizontal diagonal, respectively, and let $\rho$ be
$180^\circ-$rotation about \O/, the center of the rhombus.
Then
$$
-\rho^*dh=-\overline{\mu_h^*dh}=\overline{\mu_v^*dh}=dh.
$$
\end{lemma}

By saying that a one-form is real or imaginary on a curve we
mean that  evaluation of the one-form on tangent vectors to
that curve produces  real or imaginary values. The Lemma does
not depend on our choice of sign of $dh$. We will make the
choice of sign according to the following  geometric 
consideration: {\em according to the Lemma, the form $dh$ is real
along the vertical  diagonal, where it is never zero. We
choose the sign of $dh$ so that $dh$  is positive on
upward-pointing vectors tangent to the vertical diagonal.} 

It follows immediately from the Lemma above and our choice of
sign of $dh$ that 

\begin{corollary}\label{dhcorollary} The involutions $\mu_v$,
$\mu_h$, and $\rho$ are isometries of the cone metric $|dh|$.
In particular, the fixed  point sets of the reflections (the
vertical and the horizontal diagonals) develop into  straight
lines under the developing map $$p\rightarrow \int_{\cal O }^p
dh.$$  The developed image of the horizontal diagonal is a
vertical line.  The developed image of the vertical diagonal
is a horizontal line, and the upper half of the vertical
diagonal develops to the positive $x$-axis (i.e., the image
of $p$ develops to the right on a  horizontal line as one
moves up the horizontal diagonal from \O/).
\end{corollary} 

\begin{proof}[Proof of Lemma~\ref{dhlemma}]

In the slit domain, the one-form $gdh$ is given by $d\z$:
here we set $\z=\xi+i\eta$.  On the rhombus,  the one-form
$gdh$ is produced by pulling back $d\z$ from the slit domain.
The imaginary axis in the slit domain corresponds  to the
horizontal diagonal in the rhombus model: it is the developed
image under $gdh$ of the horizontal diagonal. (See
Figure~\ref{figure:slitandrhombus}.) Along the imaginary
axis in the slit domain, $d\z$ is imaginary:
$d\z(\frac{\bdy}{\bdy\eta})=i$.

If we parametrize the horizontal diagonal from left to right
by $s$, then the point $s$ on the rhombus develops to a
point $i\eta=if(s)$ on the  imaginary axis in the slit model, for
some real-valued function $f(s)$ with $f'(s)>0$. Then at $s$:
$$
gdh(\frac{d}{ds}\bigm|_s) = f'(s)
d\z|_{if(s)}(\frac{\partial}{\partial \eta})\in i{\bfR}.
$$
Since on the rhombus 
$(1/g)dh={\bf t}^*(gdh)$, where ${\bf t}$ is a horizontal
translation, we can write 
\begin{equation*}
\begin{split}
\f1gdh(\frac{d}{ds}\bigm|_s)  &= {\bf t}^*(gdh)
(\frac{d}{ds}\bigm|_s)\\
&= (gdh)({\bf t}_*\frac{d}{ds}\bigm|_s)\\
&= (gdh)(\frac{d}{ds}\bigm|_{{\bf t}(s)})\in i{\bf R}_+.
\end{split}
\end{equation*}
But then 
$$
dh^2(\frac{d}{dt},\frac{d}{dt})=
gdh(\frac{d}{dt})(1/g)dh(\frac{d}{dt})
$$
is negative. Therefore $dh(\frac{d}{ds})$ is imaginary along
the horizontal  diagonal. 

Recall that $dh=\pm\sqrt{gdh\cdot (1/g)dh}$ has simple zeros
at the points $V_1$ and $V_2$, and simple poles at the points
$E_1$ and $E_2$. Also, the pair of points $E_1$ and $E_2$ and 
the pair of points $V_1$ and
$V_2$ are each interchanged by the reflection $\mu_v$ in the
vertical diagonal. Therefore, $\overline{\mu_v^*dh}$ is a
meromorphic one-form whose poles and zeros  match those of
$dh$. Hence, up to a nonzero scale factor, the form
$\overline{\mu_v^*dh}$ is equal to $dh$. To determine the
scale factor we evaluate the form at \O/, the center point of
the rhombus. Since \O/ is on the horizontal diagonal where 
$dh$ is imaginary, since \O/ is fixed by $\mu_v$,  and since
$\mu_{v*}$ changes the sign of horizontal  vectors, we have
that 
$$
\overline{\mu_v^*dh(\frac{d}{ds})}
=\overline{dh(-\frac{d}{ds})}  =dh(\frac{d}{ds})
$$
at \O/. Hence $\overline{\mu_v^*dh}= dh$ on the torus. On the 
vertical diagonal, which is fixed by $\mu_v$, vertical vectors are also 
fixed by $\mu_v*$. Hence $\overline{dh}=\overline{\mu_v^*dh}=dh$ along 
the vertical diagonal, 
which implies that $dh$ is real along the vertical diagonal. 

In an analogous manner, observe that $\rho^*dh$ has the same
poles and zeros as $dh$, and that, at \O/, we have
$\rho^*dh=-dh$.  Hence, $\rho^*dh=-dh$ everywhere. Since
$\mu_h=\mu_v\rho$, we have $\overline{\mu_h^*dh} = -dh$. This
completes the proof of the Lemma.
\end{proof}

\subsection{Solving the vertical period problem geometrically 
\label{solvevpc}}
\Fig_a_and_b

In order to try to satisfy the vertical period condition
(\ref{vpc2}) we  will vary the conformal diffeomorphisms with
which we pull back $d\zeta$ in the slit  model to produce $gdh$
on the rhombus. Any such conformal  diffeomorphism is
determined by the position of the inverse image, $V_2$, of the 
vertex point $v_2$ of the slit model. We may write $V_2=a$,
$0<a<1$, $a\neq 1/2$. Our goal is to find a value of $a$
between $0$ and $1/2$ for which  the vertical period
problem (\ref{vpc2}) is solved by using the intermediate 
value theorem.

\begin{remark} Because $a+b=1$, the restriction of $a$ to lie
in the  interval $(0,\f12)$ is equivalent to the requirement
that $a<b$. As noted in  Remark~\ref{a<b?}, there are no
values of $a>1/2$ that satisfy the  horizontal period
condition (\ref{hpc2}).
\end{remark}

Even though they do not produce admissible Weierstrass
data---see  Remark~\ref{degeneration}---the values $a=0$, and
$a=\frac 12$ do define conformal diffeomorphisms from the
rhombus to the slit model, and hence they produce well-defined
one-forms $gdh$. The extreme value $a=\frac12$ represents the 
case where $E_i=V_i$, $i=1$, $2$; the extreme value $a=0$ is
case where $E_1=E_2$ and $V_1=V_2$. In each of these  cases we
may define $dh$ by translating $gdh$ to produce $(1/g)dh$ and
then taking the square root of the product of these two forms.

\FigRectangle_R

For each $a\in[0,\frac12]$, define the map $F_a
:$\T/$\rightarrow {\bf C}$ by
$$
F_a(z)=\int^z_{\cal O}dh.
$$
The vertical period condition (\ref{vpc2}) is
$$
Re\int_B dh=0,
$$
where $B$ is the curve in the rhombus illustrated in 
Figure~\ref{figure:Rectangle_R}.

We may also consider $F_a$ to be a map from  \aR/ in
Figure~\ref{figure:Rectangle_R}---a  rectangle that 
comprises half  of \T/---
to the complex plane. According to Lemma~\ref{dhlemma}, we
have $\rho^*dh=-dh$. Therefore, we may reexpress the period
condition  (\ref{vpc2}) as
\begin{equation}
\Re\int_B dh=0, \label{p2} 
\end{equation}
where $B$ is 
now considered to be the diagonal in \aR/ from \O/ to $P$.
(See Figures~\ref{figure:Rectangle_R} and \ref{figure:DhonRone}.)

Considering $F_a$ as a map from \aR/ to ${\bf C}$, we define
\begin{equation}\label{f}
f(a) := \Re \{F_a(P)-F_a({\cal O})\},
\end{equation}
in terms of which we may restate the vertical 
period condition (\ref{p2}) as asserting that
our choice of $a$ must satisfy
\begin{equation}
\label{horizontalB}
f(a)=0
\end{equation}

\FigDhonRone

\begin{proposition} There exists a value of $a$, with $0<a<\f12$
for  which the vertical period condition (\ref{p2}) is 
satisfied. \label{findgooda}
\end{proposition}

\begin{proof}
We will prove the proposition by showing the existence of a
value of $a$  satisfying (\ref{horizontalB}). We begin with a
discussion of the continuity of $F_a$ and $f$. Recall that
$dh$ is defined as a square root of $(\f1g)dh\cdot gdh$. The
one-form $gdh$ is a pullback to the
rhombus of $d\z$ in the slit model. 
We will write $\eta _a$ for the pullback of
$d\z$ that corresponds to the choice of $gdh$ with a double
zero at $z=a$ and a double pole at $z=-(1-a)$. In particular,
$\eta_{_0}$ has a double zero at $z=0$, (which is the center 
point \O/) and a double pole at the vertex at $z=1$. Let 
${\bf t}_a$ be the conformal diffeomorphism 
${\bf t}_a(z)= z-a$ of the rhombus induced by
horizontal translation by $-a$, with $a\in{\bf R}$. Then
$$
\eta_{_a}= {\bf t}^*_{a}\eta_{_0}.
$$
The one-form $\eta_{_a}$ has a double zero at $a$ and a double
pole at $-b=-(1-a)$ as required. For $(gdh)_a=\eta_{_a}$, the
corresponding one-form $(\f1gdh)_a$ is
${\bf t}^*_{-2a}gdh={\bf t}^*_{-2a}{\bf t}^*_{+a}\eta_{_0}={\bf
t}^*_{-a}\eta_{_0}$. Then,
$$
dh_a^2 = {\bf t}^*_{-a}\eta_{_0}\cdot{\bf 
t}^*_{a}\eta_{_0}=\eta_{-a}\cdot\eta_{a}.
$$
In particular, 
\begin{equation} \label{dh_0}
dh_0^2 =\eta_{_0}^2,
\end{equation}
while $dh_{\f12}^2$, having no poles or zeros, is a constant multiple of
$dz^2$:
\begin{equation} \label{dh_half}
dh_{\f12}=cdz
\end{equation}

It is clear that $dh_a$ depends continuously on $a$ on the open interval
$(0,\f12)$. 

{\em Claim.} $(dh_a)$  depends continuously on $a$ on
the closed interval $[0,\f12]$. 

The claim has two  immediate consequences. According to 
Corollary~\ref{dhcorollary}, for $a\in (0,\f12)$,  the developed image under 
$dh$ of the upper half of the vertical diagonal lies on the positive 
$x-$axis 
(assuming \O/ develops to the origin in the plane). By  the claim, the 
same must be true for $dh_0$ and $dh_{\f12}^2$. It follows immediately from
(\ref{dh_half}) that 
\begin{equation} \label{dh_half1}
dh_{\f12}=-c_1idz
\end{equation} 
for some positive real constant $c_1$.
Turning our attention to $dh_0$, it follows from
(\ref{dh_0})  that $dh_0$ must have a double zero at \O/, 
and the developed 
image 
of the horizontal diagonal to the right of \O/ must lie on the negative 
imaginary axis. Using the fact that $\eta_0$ is the pullback to a rhombus
of $d\z$ in the slit model by an orientation preserving diffeomorphism, it 
follows from (\ref{dh_0}) that

\begin{equation} \label{dh_01}
dh_0 =\eta_{_0}.
\end{equation}

{\em Proof of claim.}  There is a potential problem
at the endpoints where either a zero coalesces with a pole
(at $a=1/2$) or the zeroes coalesce and the poles coalesce (at
$a=0$). We will address this by using sigma functions to
represent the one-forms $\eta_{_a}$.

We recall from the theory of the sigma function $\sigma(z)$ on a
lattice that any meromorphic form $\vp(z)dz$ on a torus may be
expressed as a ratio
$$\vp(z)dz=C\prod^n_{k=1}\sigma(z-s_k)\sigma(z-t_k)^{-1}dz,$$ where each 
$s_k$
represents an orbit of zeros and each $t_k$ represents an orbit
of poles. In this representation, it is crucial that
$\sum^n_{k=1}s_k=\sum^n_{k=1}t_k$, and zeros of multiplicity are
considered as separate entries in the list in the customary way.

With this representation, we may write
$$
\eta_{_0} = C\frac{\sigma^2(z)}{\sigma(z-1)\sigma(z+1)}dz
$$
for some complex constant $C$ (independent of $a$). 
It then follows, 
using ${\bf t}_a(z)= z-a$, that
$$
(gdh)_{a} = \eta_{_a} = t^*_{a}\eta_{_0} = 
C\frac{\sigma^2(z-a)}{\sigma(z-a-1)\sigma(z-a+1)}dz
$$
and 
$$
(\f1gdh)_{a}=t^*_{-a}\eta_{_0}= 
C\frac{\sigma^2(z+a)}{\sigma(z+a-1)\sigma(z+a+1)}dz.
$$
Therefore 
$$
\frac{(dh)_a^2}{dz^2} = C^2
\frac{\sigma^2(z-a)\sigma^2(z+a)}
{\sigma(z-a-1)\sigma(z-a+1)\sigma(z+a-1)\sigma(z+a+1)}.
$$
As $a\rightarrow 1/2$, this function limits on
$C^2\frac{\sigma(z-1/2)\sigma(z+1/2)}{\sigma(z-3/2)\sigma(z+3/2)}=C^2$,
since $\sigma(z+2)=\sigma(z+\om_1+\om_2)$ differs from
$-\sigma(z)$ by a factor of the form $e^{\a z+\beta}$ where $\a$
and $\beta$ depend only on the lattice. Also, note that as
$a\rightarrow 0$, the function $(dh_a/dz)^2$ limits on
$\eta_{_0}^2/dz^2$.  This completes the proof of the claim.
\vspace{.2in}

In the three-dimensional product
$[0,\f12]\times$\aR/ of the interval $[0,\f12]$ with
the rectangle \aR/, the set 
$$
\{(a,p) | dh_a \text{ has a pole at } p\}
$$
consists of two line segments  on the boundary of this 
(three-dimensional) box. One 
line 
segment connects the bottom left corner of \aR/$\times \{0\}$  
to the midpoint of 
the top of \aR/$\times \f12$.  
The other line segment connects the midpoint 
of the bottom of \aR/$\times \{0\}$
to the top right corner of \aR/$\times\{\f12\}$. 
After removing from $[0,\f12]\times$\aR/
a small tubular neighborhood ${\cal N}$ of these line 
segments, we may assert that 
$F:[0,\f12]\times$\aR/$\setminus{\cal N}\rightarrow{\bfC}\union\infty$
defined by 
$$
F(a,p)=F_a(p)=\int_0^z dh_a
$$
is continuous and bounded. Since $[0,\f12]\times B$,
where $B$ is the path of integration from ${\cal O}$ to $P$, 
lies in the domain of 
$F$, we may assert that
$$f(a):= \Re \{F_a(P)-F_a({\cal O})\}$$ is continuous on $[0,\f12]$. 
(The generic image of $F_a$, $0<a<\f12$, is illustrated in
Figure~\ref{figure:DhonRone}. In this case, $a\neq b$ and
$E_i\neq V_i$, $i=1,2$. ) 
According to (\ref{dh_half1}), $$F_{\f12}(p)=\int_0^p 
dh_{\f12}= -c_1 i\int_0^pdz =-c_1 i p,$$
where $c_1$ is a positive constant, i.e. $F_{\f12}$ is clockwise rotation by 
$-\pi/2$ followed by a scaling. (See Figure~\ref{figure:DhonRtwo}.) In 
particular,
\begin{equation}
f(\f12)= \Re \{F_{\f12}(P)-F_a({\cal O})\}>0. \label{apositive}
\end{equation}

\FigDhonRtwo

We now consider the other extreme case: $a=0$.
\newline
{\em Claim:}
\begin{equation}
\label{anegative}
f(0)<0.
\end{equation}
The proposition follows from (\ref{apositive}),
(\ref{anegative}) and the intermediate value theorem.
\newline

{\em Proof of Claim (\ref{anegative})} According to (\ref{dh_01}),  
$dh_0=\eta_{_0}$. The symmetry lines of 
$|dh|$ on the
rhombus---the horizontal and vertical diagonals according
to  Corollary~\ref{dhcorollary}---correspond to the 
imaginary and real axes,  respectively, in the slit model.

\FigTjigsaw  

If we remove the symmetry lines from the rhombus---also the symmetry lines
of $|dh|$ according to Corollary~\ref{dhcorollary}-- what is
left consists of  two rectangles, one of which is \aR/.  If
we remove the symmetry  lines from the slit model, 
 what is left consists of two  conformal
rectangles, the ones  illustrated in 
Figure~\ref{figure:Tjigsaw}. Each rectangle has two vertices
at infinity. Since $dh_0=\eta_0$ is the pullback of $d\z$ on the slit model, 
the developing map $F_0$ must send \aR/ onto one of these rectangles. Again, 
as with
all $F_a$, the map $F_0$ takes vertical (resp horizontal)
boundary segments to horizontal (resp vertical) lines.  
Using  Lemma~\ref{dhlemma},  we know that the image
under $F_a$ of the right-hand-side of \aR/ is a horizontal
line going to the right as one ascends from \O/. Therefore,  
the image $F_0({\cal R})$ must be the conformal rectangle on
the right of Figure~\ref{figure:Tjigsaw}  and the left of
Figure~\ref{figure:DhonRthree}, with
$F_0({\cal O})$ equal to the  right-hand  vertex of this
rectangle  and $F_0(P)$ equal to the left-hand vertex. Hence
$F_0(P)$ lies to the left of $F_0({\cal O})$ which means that
$$f(0)= \Re \{F_0(P)-F_0({\cal O})\}<0,$$
which is (\ref{anegative}). This 
completes the proof of the claim, which completes the proof of the 
proposition.
\end{proof}

\FigDhonRthree

\newpage

\section{The construction of screw-motion-invariant \Hk/}
  
Our plan is to construct a family of surfaces $\{$\Hk/$\}$
that includes the surface \Hone/; in this family, there should
exist at least one surface of the type \Hk/ for each $k\ge1$.
We begin by specifying  necessary conditions for the
Weierstrass data of $\hksk$. 
A fundamental domain of \Hk/ is the image of $\hksk$ under
the  Weierstrass mapping. It is a genus-one surface with two
ends, possessing a vertical axis and containing two
horizontal lines making an angle of $\pi k$, one with the
other. We have established  in  Lemma~\ref{symmetry}
of Section~2 that $\hksk$ is a rhombic torus and that it has
the same  symmetries as \Hone/$/{\bf T}$, where ${\bf T}:=\sigma_1$ is 
vertical translation by $2\pi$.

\FigDivisorsFive

Because we want the \Hk/ family to  give a deformation of
\Hone/, we are  justified in assuming that the placement of
vertical points of \Hk/  conforms qualitatively to what
happens on \Hone/. Specifically, the vertical points lie
on the same lines in \Hk/ as they do on \Hone/,  and in the
same relative position. As we did in Section~3, label the
ends $E_1$, $E_2$ and the vertical points $V_1$, $V_2$. All
four of these points lie on a symmetry line of $\hksk$, which
we expect to be mapped into the horizontal lines of \Hk/.
Another symmetry line of the closed Riemann surface \T/ that will 
be mapped into 
the
vertical axis of \Hk/
must meet the first symmetry line
orthogonally in two points that we may assume, without loss
of generality, lie between $E_1$ and $E_2$ and $V_1$ and $V_2$
respectively. Also, we will assume that $g(E_1)=\infty$ and
$g(E_2)=0$. We expect the Gauss map to behave locally---near
$E_1$ and $E_2$---like $z^k$ near infinity and the origin,
respectively. At one end we expect a pole of order $k$, and
at the other a zero of order $k$. Since we will consider $k$ 
to  take on all real values greater than $1/2$, we will be
considering multivalued Gauss maps and one-forms $gdh$ and
$(1/g)dh$.

Because we want the family \Hk/ to be continuous in $k$, we
require---as is  the case on \Hone/---that $g(V_1)=\infty$
and  $g(V_2)=0$. This determines the divisors of $gdh$ and
associated forms that are  presented in
Figure~\ref{figure:divisors5}.

\FigRectangleRTwo

\subsection{The rhombic ($|dz|$) model  \label{rhombmodel}}     
\FigFundDomainhK
In Section~3, we chose to represent \T/, the Riemann surface
of \Hone/  modulo translations, as a rhombus in the following manner: a
 diagonal is a horizontal line segment of
length two, whose center point  is the origin of ${\bf C}$.
In dealing with  family \Hk/, we do not know in  advance 
the underlying rhombic structure. Also, we will choose a
different  normalization for rhombi underlying $\hksk$. Here is our 
convention:

For each $\hksk$, the associated rhombus in the complex plane
will be chosen to have the following  properties: the top
vertex of the rhombus sits at $0\in\bfC$; the point $1\in\bfC$
is the right-most vertex; the left-most vertex is at some
unitary value $\tau$; the points $V_1$ and $V_2$ lie 
symmetrically placed on the diagonal from $1$ to $\tau$. The
last statement  is a consequence of Abel's theorem applied to
$dh$ (or a consequence of the  symmetry of the surface
imposed by Lemma~\ref{symmetry}). In this  setup, the center
point of the rhombus lies at ${\cal O}=\f{1+\tau}2$, and a
point on the horizontal diagonal must be of the form
${\cal O}\pm t(\f{1-\tau}2)$, $-1\le t\le1$. In particular,
the  points $V_i$, $E_i$, $i=1$, $2$ are of the form
\begin{equation}
\begin{split}\label{ab}
E_1   &= {\cal O} - b\frac{1-\tau}2\\
E_2   &= {\cal O} + b\frac{1-\tau}2\\
V_1   &= {\cal O} - a\frac{1-\tau}2\\
V_2   &= {\cal O}  + a\frac{1-\tau}2
\end{split}
\end{equation}
for some $0<a<b<1$.

Because we want the family \Hk/ to
include \Hone/, we may assume that 
$a<b$, as we know this to be the case for that surface.
We will refer to this model of \t/ as the rhombic or
the rhombic ($dz$) model, the latter when we 
wish to emphasize that we are understanding the
rhombus as coming equipped with a (non-singular)
cone metric. In subsequent sections, we will develop 
two other models of \t/ using the forms $gdh$ and $dh$.

The positions of the ends $E_i$ and the vertical points
$V_i$  are not independent. For $k=1$, we observed in Section~3
that $a+b=1$, a consequence of Abel's theorem. In general, we have the 
following result.
 
\begin{proposition}
\label{a+kb}
Suppose Weierstrass data is given that satisfies the
geometric  conditions (4) for  $\hksk$; in particular, $g$
possesses the divisor specified in
Figure~\ref{figure:divisors5}. If $k$ is an integer, then
$a+kb$ is also an integer and
$$
1\leq a+kb\leq k.
$$
In particular if $k=1$, $a+b=1$ . 

If such Weierstrass data exists for a continuously varying
family of $\hksk$ that contains \Hone//$\sigma_1$, then
\begin{equation}
a + kb = k \label{abk}
\end{equation} 
\end{proposition}

\begin{proof} If $k$ is an integer, then the function $g$ is
single-valued. We will apply Abel's Theorem to $g$, whose 
divisor is given in Figure~\ref{figure:divisors5}. From
(\ref{ab}), we have
$$
-kE_1-V_1+V_2+kE_2=(a+kb)(\tau-1).
$$
Therefore, by Abel's theorem, $(a+kb)(\tau-1)$ is in the
lattice, which implies that $a+kb$ is an integer. Since
$0<a<b<1$, we see that $a+kb$ is an integer between $1$ and
$k$. In particular, when $k=1$, we have $a+b=1$.

For any (possibly non-integral) $k>0$, consider the
meromorphic one-form $\f{dg}g$, whose divisor and residues are
given in Figure~\ref{figure:divisors5}.
\FigGammaPaths
Let $\g_1$ and $\g_2$ be the cycles on ${\cal
T}=\bfC/\{1,\tau\}$ given by the vectors $1$ and $\tau$,
respectively. For any closed one-form, $\mu$, with simple
poles at $P_1,\dots,P_r$, and residues $a_1,\dots,a_r$, the 
bilinear relation \cite{FK91} gives
\begin{equation} \label{bilinearrelation}
2\pi i\sum^r_{j=1} P_j a_j = \om_2\a_1 - \om_1\a_2,
\end{equation}
where $\a_i=\int_{\g_i}\mu$ and $\om_i=\int_{\g_i}dz$,
$i=1$, $2$. Applying (\ref{bilinearrelation}) to $\mu=\f{dg}g$ we have
$$2\pi i(-kE_1 - V_1 + V_2 + kE_2) = \a_1\tau - \a_2,$$
and using (\ref{ab}), we obtain
\begin{equation}
2\pi i(a + kb)(1 - \tau) =   \a_1\tau-\a_2.  \label{int}
\end{equation}

It is left to evaluate $\a_1$ and $\a_2$. 

In Figure~\ref{figure:gammapaths}, we have drawn curves
$\tilde{\g_1}$ and $\tilde{\g_2}$ that are homotopic to
${\g_1}$ and ${\g_2}$, respectively. On the semi-circular
arcs  near the $E_i$,  the integral of $\f{dg}g$ will give
one half  of $2\pi i$ times the residue with a sign change
due to the orientation of the semicircles  associated to
$\tilde{\gamma_1}$. From the residues given in
Figure~\ref{figure:divisors5}, we know that the
contributions to  the integral  of $\f{dg}{g}$ along the
semicircles is the same on $\tilde{\gamma_1}$ as on
$\tilde{\gamma_2}$ and is equal to $-\pi i(k+1).$
 
The horizontal line segments of each $\tilde{\g_i}$ are mapped by the
immersion into horizontal lines, which implies that on each of the
segments, the Gauss map $g$ takes values on  radial lines in $\bfC$.
Since $\f{dg}g=d\log g$, the integral of $\f{dg}g$ along these lines
depends only on the change  in $|g|$ along these segments. We can
assume without loss of generality  that  each of the semicircular arcs
begins and ends at points where $|g|$ has the  same value. Since $|g|=1$
at the center and at the vertex of the  rhombus---both points on the
vertical diagonal which is mapped into a  vertical line---we  can
conclude that the total contribution to $\alpha_i$  of integrating
$\f{dg}g$ along these segments is zero.

It remains to compute $\int_A\f{dg}g$,  where $A$ is the top
half of the vertical diagonal. See Figure~\ref{figure:gammapaths}. Before 
doing this, we observe that since this segment is common to
$\tilde{\g_1}$  and $\tilde{\g_2}$ we now know that 
\begin{equation}
\alpha_1=\alpha_2=-\pi i(k+1) +\int_A \f{dg}g; \label {int4}
\end{equation}
and we denote this common value as $\alpha$. From (\ref{int}) we
have
$$2\pi i(a + kb)(1 - \tau) = -\alpha(1 - \tau)$$  
or
\begin{equation}
a+kb= \f{-\a}{2\pi i} \label{int1}.
\end{equation}
As noted above, the vertical diagonal is mapped by the
immersion into the  vertical axis of \Hk/ along which $g$ is
unitary.  By our assumption that \Hk/ has the geometric
properties outlined in (4) and in Lemma~1, the tangent planes at the
endpoints of $A$ are vertical planes making an angle of $\pi
k$ with one another.  Hence
\begin{equation}
\int_A\f{dg}g = -\pi i (k + N(k)) \label{k+N(k)} 
\end{equation}
for some integer $N(k)$. Therefore from (\ref{int4})
\begin{equation}
\a= -\pi i(k+1)-\pi i (k + N(k))= -\pi i (2k +N(k) +1)\label{int2} 
\end{equation}
and from (\ref{int1})
\begin{equation}
a+kb=-\frac{\a}{2\pi i} = k+ \frac{(N(k)+1)}{2}\label{int3}
\end{equation}

Equation (\ref{int3}) is valid for any value of $k>0$. If we
assume that we have Weierstrass data for a continuous family
of \Hk/$/\sigma_k$, then the integer-valued function  $N(k)$
is a constant. If that family includes \Hone/ then it follows
from (\ref{int3}) and the first part of the proposition that
$N(k)= -1$. Thus (\ref{int3}) states that $a+kb=k$, which is
(\ref{abk}).
\end{proof}

Since $0<a<b<1$ we have

\begin{corollary} Under the assumptions of Proposition
\ref{a+kb}, \label{corabk}
$$ 
a < \f k{k+1} < b,\quad\text{and}\quad\lim_{k\to\infty}b=1.
$$
\end{corollary}

\begin{remark}
The integral $\int_A\f{dg}g$ in Proposition~\ref{a+kb}
measures the turning  of the normal along the vertical axis
of \Hk/: the full turning in one period of \Hk/ will be
$2\pi(k-1)$ according to the equation (\ref{k+N(k)}) in the
proof of  Proposition~\ref{a+kb} and the fact that  $N(k)$ is
identically equal to $-1$. On the helicoid, \He/, the normal
turns by $2\pi k$ on each \He//$\sigma_k$. We interpret this
as saying that  the presence of a handle in each  fundamental 
domain of \Hk/ has the effect of costing one full turn of the
normal. In particular, the normals on \Hone/  do not
wind around the vertical axis at all, a  surprising 
geometric consequence of Abel's theorem.
\end{remark}
\label{dz}

\subsubsection{The period conditions}
The Weierstrass data we produce must satisfy the  horizontal and vertical  
period
conditions (\ref{hpc2})and (\ref{vpc2}): 
\begin{eqnarray}
 \int_B gdh &=& \overline{\int_B \f1gdh} \label{hpc3};\\
Re \int_Bdh&=&0. \label{vpc3}
\end{eqnarray}
The curve $B$ is defined in 
Figure~\ref{figure:rectanglertwo}.  

\subsection{The $|gdh|$ model and the solution of the
horizontal period  problem}\label{sec:gdhmodelsection} 
We will produce candidate Weierstass data that depend on the
conformal  type and placement of the  distinguished points.
They will satisfy the  divisor conditions of
Figure~\ref{figure:divisors5}, all the symmetry  conditions
of Lemma~\ref{symmetry}, the distinguished-point-placement
requirements of (\ref{abk}) and the  horizontal  period
condition (\ref{hpc3}), but not necessarily the vertical
period  condition (\ref{vpc3}).  We will do this by 
constructing a second model to which we will sometimes refer as the
``$|gdh|$ model'' but more often---for reasons that will become evident--- as
the {\em slit model} or the $|d\zeta|$- model.
 
The construction of \T/ in Section~3.3 involved slicing the
$\zeta$-plane  along $[-1,1]$ and making identifications. Into
this model of \T/  we are going  to sew in a copy of the cone
$S_{k-1}$. This 
cone metric construction procedure is described in  Section~2.3, 
Example~2. For each $d>0$ and $k>0$, we slice the
$\zeta$-plane from $di$ to $\infty$ along the imaginary axis,
and sew in $S_{k-1}$, placing the cone point of $S_{k-1}$ with
positive cone angle at the point of \T/  corresponding to
$\zeta = di$, and the cone  point with negative cone angle at
$\infty$.  

{\center We will refer to this torus as \tkd/.}

It has three
cone points, $\zeta=0$, $\zeta=di$, and $\zeta=\infty$, with
cone angles $6\pi$, $2\pi k$ and $-2\pi(k+1)$, respectively.
We label these points $v_2$, $e_2$ and $e_1$ respectively.
See Figure~\ref{figure:tkd}. The multivalued one-form $d\zeta$ on 
the $\zeta$ plane with an $S_{k-1}$ sewn in descends to a
multivalued one-form   on \tkd/ with a zero of order $k-1$ at
$e_2$, a pole of order $k+1$ at $e_1$ and a double zero at
$v_2$.  (Our choice of notation for these points comes from the desire to 
have the $e_i$ and the $v_i$ correspond to the $E_i$ and $V_i$ in the
rhombic model.)

{\center We will also refer to this one-form on \tkd/ as 
$d\zeta$.}
 
All these special points lie on a fixed point set
of an  involution of \tkd/, namely the imaginary axis (fixed
under $\z\rightarrow -\bar\z$).
The rhombic torus \tkd/
is conformally  diffeomorphic to one of the tori constructed in 
Section~\ref{dz}, which underly the
$|dz|$ model. It is clear that  we may choose a conformal 
diffeomorphism---from the rhombus model to the slit model of the 
torus we have
just  constructed---that  has the following properties: 
\begin{itemize}
\item[i.]the
imaginary axis of the $\z$-plane is the image of the
horizontal diagonal of a (rhombic) fundamental domain; 
\item [ii.]the
preimages $E_1$ and $E_2$ of $e_1=\infty$ and $e_2=di$ 
(respectively) are symmetric with respect to
the vertical diagonal of the rhombus;
\item [iii.] the point $E_1$, the preimage of the point $e_1=\infty$ in  the 
slit model, lies to the left of the center point \O/ in the rhombus.
\end {itemize} As pointed 
out  in Section~4.1, Abel's theorem applied to $dh$ forces us
to define $V_1$ as the point symmetric to $V_2$ on the
horizontal diagonal of the rhombus and therefore to identify
$v_1$ as a uniquely specified point on the imaginary axis of
the $\z$-plane. Because of the order of points on the
horizontal diagonal on the rhombus, we know that $v_1$ lies
on the negative imaginary axis of the $\z$-plane.

Note that for each choice of $k$ and $d$, our geometric
normalization  determines $a=a(k,d)$ and $b=b(k,d)$ uniquely.
Also, the lattice generated  by $\{1, \tau\}$  that is
associated to \tkd/ changes as $k$ and $d$ vary. That is
$\tau =\tau (k,d)$.

\FigTkd

\begin{definition}\label{gdh2}
The one-form {\em $gdh$} on the torus \tkd/ is given by $d\zeta$ in the slit 
model. On the 
rhombic model of \tkd/ it is given 
by the pullback of  $d\zeta$ under the conformal diffeomorphism  between the 
two models, and the conformal diffeomorphism
is determined by conditions (i)---(iii) above.
\end{definition}

\begin{remark}  Our construction in this section is ambiguous
when $k=1$ and we make it precise here in a manner that is
consistent with the discussion of \Hone/ in Section~3.  When
$k=1$ we are not sewing in a cone at all, which means that
$\tau(1,d)$ is constant: the rhombus we get is the rhombus of
\T/ constructed in Section~3. For any choice of  conformal
diffeomorphism that places the  inverse image of the point at 
infinity at $E_1$ (a point on the horizontal diagonal to the
left of \O/),  let $E_2$ be the  symmetric point on the
horizontal diagonal to the right  of \O/. Then the image of
$E_2$ under this conformal diffeomorphism is some point,
$\z=di$, on the positive imaginary axis.   
\end{remark}

The one-form $gdh$ is multivalued when $k$ it not an integer.
The periods of $gdh$ (and its companion one-form $\f1gdh$ to be
defined below) will be defined according to the following
convention.

\begin{convention} If $\a$ is a curve in the slit model \tkd/
or in the rhombic model, we specify a base point $p$ on $\a$ and
a choice of value for $gdh$ at $p$. There is a unique continuous
choice of values of $gdh$ along $\a$ that agrees with our choice
at $p$, and it is with this choice that we evaluate $\int_\a
gdh$. In practice we will choose a preferred branch of $gdh$. In
the slit model \tkd/, we choose a branch of $d\z$ that is real
on the part of the real axis where $|\z|>1$; this branch extends
in a unique manner to \tkd/ after removing a cut from $e_2$
to $e_1$. The cut we will choose is in the cone and is the image
of the segment of the horizontal diagonal in the rhombus model
that runs from $E_2$ to $E_1$, passing through the vertex of the
rhombus. It follows that, on the rhombus, we are choosing the
branch of $gdh$ that is well-defined
in the complement of this same slit, and that is also imaginary along the 
remaining 
segment of the
horizontal diagonal.

When dealing with specific curves that do not cross the branch
cuts, we will always choose this branch to evaluate and we will
not need to specify a point on the curve in question.
\end{convention}

With this convention in mind we evaluate $\int_Bgdh$, where $B$
is the cycle in the rhombic model in Figure~\ref{figure:thecurveB}. It is 
illustrated
there with its diffeomorphic image, also labelled $B$, in the
 slit model.
\FigTheCurveB
 Since
these curves do not cross the cut used to define the principal
branch of $gdh$, we have
\begin{equation}
\int_B gdh = \int_{B}d\z = \pm1\in\bfR.\label{gdhred}
\end{equation}

We note that this choice of branch is consistent with the
case $k=1$, as described in the proof of Lemma~\ref{dhlemma}. 
The one-form $gdh$ has the divisor specified in
Figure~\ref{figure:divisors5}. Also, as was the case for the
$gdh$ constructed for \Hone/ in Section 3, it follows
immediately from the definition of $gdh$ and the convention
above that all of the periods of $gdh$ are real.

\subsubsection{The $|\f1gdh|$ model and the horizontal period
condition
\label{section:hpp} }
We will now specify the mulivalued one-form $\f1gdh$. Recall
from Section~\ref{Hk} the involution $r$ of the rhombus
which is a 180$^\circ$-rotation about the center point \O/. Note
that $r(E_1)=E_2$ and  $r(V_1)=V_2$. We expect
geometrically that $g\circ r  =\f1g$ and $r^*dh=-dh$. 
This motivates the formal definition of $\f1g dh$ on the rhombic model: 
\begin{equation} \label{oneovergdh2}
\f1g dh := -r^*(gdh).
\end{equation}
The one-form $\f1g dh$ has the divisor specified in 
Figure~\ref{figure:divisors5}. Furthermore, the pair 
$\{\f1g dh, gdh\}$ satisfies the horizontal period condition
(\ref{hpc3}). To see this, note that the cycle $B$ in
Figure~\ref{figure:thecurveB} satisfies $r\circ B=-B$ and
therefore (using $r^2=id$ in the penultimate equality)
$$
\int_B\f1g dh = -\int_{B}r^*gdh = -\int_{r^{-1}(B)}gdh
= -\int_{r\circ B} gdh = \int_B gdh.
$$

Since we know from (\ref{gdhred}) that $gdh$ has real periods,
we see that
\begin{equation}
\int_B \f1g dh
= \overline{\int_B gdh} \label{hpc4}
\end{equation}
is satisfied. This is the horizontal period condition (\ref{hpc3}). For 
future 
reference, 
we state this result as a
proposition.

\begin{proposition} \label{prop:hpc}
For every $k>1/2$ and every $d>0$, the one-forms $gdh$ and
$\f1gdh$,  defined in Definition~\ref{gdh2} and
(\ref{oneovergdh2}), respectively,  satisfy the  horizontal period condition 
(\ref{hpc3}) in the sense of the Convention of the previous section.
\end{proposition}

\FigOneOverG

\begin{remark} 
Let $\phi$ be the diffeomorphism from the rhombus to the slit
model \tkd/ that is used to define $gdh$ in Definition~\ref{gdh2}:
$$
gdh = \phi^*(d\z).
$$
Using (\ref{oneovergdh2}), we have $-\f1gdh=(\phi\circ r)^*(d\z)$,
where $r$ is the 180$^\circ$-rotation about the center of the
rhombus.  Let $\tl r$
be the conformal involution of the slit model \tkd/ induced by $r$: 
that is, $\tl r=\phi\circ
r\circ\phi^{-1}$. Then
$$
-\f1gdh = (\tl r\circ\phi)^*(d\z).
$$
Clearly $\tl r$ interchanges the $e_i$, interchanges the $v_i$,
fixes $\SO$ and leaves the imaginary axis invariant. If $f$ is
the involution of the \tkd/ model induced by $\z \rightarrow -\z$, then
$$
\f1gdh = (f\circ\tl r\circ\phi)^*(d\z) = \phi^*(f\circ\tl r)^*(d\z).
$$
This equation above tells us that
$\f1gdh$ is the one-form on the rhombus produced by pulling back
the one-form $d\z$ on the modified slit model in 
Figure~\ref{figure:oneoverg}.
\end{remark}

\subsubsection{Symmetries inherent in the $|gdh|$  and
$|\f1gdh|$  models
\label{section:symmetries}}
We show in this section that the
Weierstrass data we have just defined satisfies the geometric
conditions required of $\hksk$ in (4) and in 
Lemma~\ref{symmetry}. 

The metric induced by the Weierstrass
immersion (\ref{weier})  is given by 
(\ref{metric}):
\begin{equation}
ds = |gdh| + |\f1g dh|. \label{metric2}
\end{equation}
Two different branches of $gdh$ (or $\f1g dh$) differ by a
unitary constant. Thus the metric expression in
(\ref{metric2}) does not depend on our choice of branch.

We may define $g$ and $dh$ to be square roots of the ratio
and the product, respectively, of $gdh$ and $\f1g dh$, and we
choose the sign of the square root as we did in Section~3 for
\Hone/. The function $g$ is then likely multi-valued, and our
convention extends to choose the branch determined by our
choice of branch of $gdh$ and $\f1g dh$. Let
$r$ be the involution of the rhombic model defined
by 180$^\circ$-rotation about the center point \O/, 
and let and $\mu_h$ and
$\mu_v$ be reflection in the horizontal and vertical
diagonals, respectively.

\begin{lemma} \label{symmetriesgdh} 
{\rm i)} The involutions $r$, $\mu_v$ and $\mu_h$ are isometries
of the metric $ds$ in (\ref{metric2}). The involution $\mu_h$ is
an isometry of the  rhombus equipped with the
$|gdh|$ metric. The involutions $r$ and $\mu_v$ are isometries
between the rhombus equipped with the metric $|gdh|$ and the
rhombus equipped with the metric $|\f1gdh|$.

{\rm ii)} The forms $gdh$ and $\f1gdh=-r^*(gdh)$ on the
rhombus are fixed by $-\bar\mu^*_h$ and interchanged by
$-r^*$ and $\bar\mu^*_v$. 

{\rm iii)} The one-form $dh$ on the rhombus is imaginary on
the horizontal diagonal and real on the vertical diagonal.
There is a definition of $g$, consistent with our convention
so that the product of $g$ and $dh$ agrees with the form
$gdh$, and the quotient of $dh$ by $g$ agrees with the form
$\f1gdh$. Moreover, $g$ is unitary on the
vertical diagonal, and $g^2$ is real on the horizontal diagonal
and satisfies $g^2\circ r=\f1{g^2}$. In particular, any branch of
the multivalued function $g$ takes segments of the horizontal
diagonal disjoint from $E_1$ and $E_2$ to a straight line.
\end{lemma}

We note that the principal branch of $gdh$ and $\f1gdh$ are
mapped into one another by these automorphisms (see our
``Convention'' in the previous section.)

\begin{corollary} \label{cor:symmetry} The Weierstrass immersion 
(\ref{weier})
defined by the forms $gdh$ and $\f1gdh$ maps the vertical
diagonal into a vertical line and the horizontal diagonal
into horizontal lines that project onto lines in the
$(x_1,x_2)$-plane making an angle of $\pm\pi k$ with one another.
\end{corollary}

\begin{proof}[Proof of Lemma] Statement i) follows immediately
from statement ii) and the form of the metric (\ref{metric2}).

Statement ii) for $-r^*$ follows immediately from the
definition of $\f1gdh$ and the fact that $r^2=id$. The
meromorphic one-form $\bar\mu^*_h(gdh)$ has the same divisor
as $gdh$. Along the portion of the horizontal diagonal where
(our chosen branch of) $gdh$ is imaginary, we therefore have
$\mu^*_h(gdh)=gdh$. Hence
\begin{equation}
\bar\mu^*_h(gdh) = -gdh,\label{mh1}
\end{equation}
(and, analogously $\bar\mu^*_h(\f1gdh)=-\f1gdh$). Similarly,
$\bar\mu^*_v(gdh)$ has the same divisor as
$\f1gdh=-r^*(gdh)$. Hence
$$
\bar\mu^*_v(gdh) = c\f1gdh = -cr^*(gdh)
$$
for sone nonzero constant $c$. Since $r\circ\mu_v=\mu_h$,
applying $r^*$ to the equation above and using (\ref{mh1})
gives $c=1$. Hence,
\begin{equation}
\bar\mu^*_v(gdh) = \f1gdh.\label{mv1}
\end{equation}
Since $\mu_h\circ r=r\circ\mu_h$ and $\mu^2_v=id$, the
equalities (\ref{mh1}) and (\ref{mv1}) imply statement (ii)
for $\mu _h$.

For statement iii), observe that since
\begin{equation}
dh^2 = gdh\cdot\f1gdh,\label{dh2}
\end{equation}
the identities (\ref{mh1}) and (\ref{mv1}) imply that
$\mu^*_v(dh^2)=\ov{dh}^2=\mu^*_h(dh^2)$. Therefore
when applied to tangent vactors along both diagonals, we have
\begin{equation}
dh^2 = \ov{dh}^2,\label{dh22}
\end{equation}
implying that $dh$ is either real or imaginary on the
diagonals. Along the horizontal diagonal, the form $gdh$ is by
definition imaginary when applied to tangent vectors, and
the same is true for $\f1gdh$, either by our convention, by
ii), or by construction of $\f1gdh$. Hence, $dh^2$ is real and
negative along the segment of horizontal diagonal from $E_1$ to
$E_2$ through $\SO$, which implies that $dh$ is imaginary. We
also know that $dh$ has no zeros on the vertical diagonal, as
neither $gdh$ nor $\f1gdh$ vanish there. Since the form $dh$ must
be real for vertical directions at \O/ (because $dh$ is
imaginary on tangent vectors to the horizontal diagonal), the
form $dh$ is real on the entire  vertical diagonal.

Our next goal is to provide a definition of $g$ that is
consistent with our conventions and in accord with statement
(iii).

To begin, define a multi-valued function $G$ by
\begin{equation}
G = gdh/\f1gdh.\label{gsq}
\end{equation}
Since both $gdh$ and $\f1gdh$ are imaginary when applied to
tangent vectors to the horizontal diagonal, the function $G$ is
real there. Since we are expecting a multivalued $g$---defined
only up to integer powers of $e^{2\pi ik}$---we have shown what
is required in statement iii) for the horizontal diagonal. 
Applying (\ref{mv1}) and the fact that $\bar\mu_v$ is an
involution to (\ref{gsq}), we have
$$
\bar\mu^*_v(G) = \f1G.
$$
This implies that when applied to tangent vectors to the
vertical diagonal we have $\overline G=\f1G$, i.e. $|G|=1$ as
required. Also, $G\circ r=\f{r^*gdh}{r^*(\f1g
dh)}=\f{-\f1g dh}{-gdh}=1/G$. We define $g$ to be the
square root of $G$. (This is permissible since we can do it on a
branch and then extend.) The statements in iii) for $g$ follow
directly from those for $G$. This completes the proof of the
lemma.
\end{proof}

\FigGammaPathsTwo
\begin{proof}[Proof of Corollary] The involutions $\mu_h$ and
$\mu_v$ are isometries of the metric $ds$ by statement i) of 
Lemma~\ref{symmetriesgdh}. Therefore,  their fixed-point sets
are geodesics. On a minimal surface, a geodesic is a straight
line if and only if it is an asymptotic curve. Moreover, a
curve on a minimal surface is asymptotic if and only if
$\f{dg}gdh$ is imaginary along it. \cite{hk2} Along the
vertical diagonal, the form $dh$ is real on tangent vectors
and the function $g$ is unitary according to statement iii) of
Lemma~\ref{symmetriesgdh}. Thus the form $\f{dg}gdh$ is
imaginary on tangent vectors to the vertical diagonal, which
means that the vertical diagonal corresponds to an asymptotic
curve, hence a straight line on the minimal surface. Since $g$
is unitary along the vertical diagonal, the corresponding straight
line on the minimal surface is vertical.

By statement iii) of the lemma, the form $dh$ is imaginary on
the horizontal diagonal --- and on any interval there not
containing $E_1$ or $E_2$. Hence $x_3=\Re\int dh=0$ on any
interval of the horizontal diagonal. It follows that the image
of segments of this diagonal bounded by $E_1$ and $E_2$ are
mapped to curves in horizontal planes in $\BR^3$. Also by
statement iii) of Lemma 3, the function $g$ takes values on a
line through the origin. This implies that $dg/g dh$ is
imaginary along those segments, so their images are horizontal
straight lines. The values of $g$ along these lines increase or
decrease by integer multiples of $\pi k$ because $\int_\a
gdh=2\pi ki$ for any simple closed curve $\a$ surrounding $E_1$ or
$E_2$. 
Hence projection of successive lines make angles of
$\pm\pi k$ with one another.
\end{proof}

For future reference, we gather together the key results of 
Sections~\ref{section:hpp} and
\ref{section:symmetries}. It is a consequence of  as a 
Proposition~\ref{prop:hpc} and 
Corollary~\ref{cor:symmetry} that 
\begin{proposition}\label{prop:almostthere}
For any values of $k>1/2$ and $d>0$, the Weierstrass data defined on \tkd/  
produces an 
\Hk/ with  both the horizontal period conditon (\ref{hpc3}) and all the 
properties of (4) and Lemma~1, 
provided that the vertical period condition (\ref{vpc3}) is solved.
\end{proposition}

To prove Theorem~3, we must show (among other things) that we can not only 
choose for each $k>1$ a value of $d>0$ for which  the horizontal period 
condition 
(\ref{vpc3}) is solved, but also that we can choose a continuous path 
$(k(t),d(t))$ along which $k(t)$ takes on all values between $1$ 
and $\infty$. We will formulate our approach to this problem in Sections~4.3 
and 
4.4, then carry it out in Sections~5 and 6. Before doing so, we verify
in the next section that 
Proposition~\ref{a+kb} is satisfied by the our Weierstrass candidates.
 
\subsubsection{The placement of the ends and vertical points 
\label{placementofends}}

Proposition~\ref{a+kb} establishes a necessary condition on the Weierstrass 
data, specifically restrictions on the relative positions of the $E_i$
and $V_i$, in order for the continuous choice of $d$ to be possible. We 
prove in this section (see also the alternative
proof in Appendix A)
that the Weierstrass data we have defined satisfies this 
condition for all admissible pairs $(k,d)$.

\begin{proposition}\label{a+kb2} For all finite $k>1/2$ and all finite $d>0$,
the points $E_i$ and $V_i$ defined by the construction of $gdh$
satisfy (\ref{abk}):
$$a+kb=k$$
when written in the form (\ref{ab}):
$E_1=\SO-b(\f{1-\tau}2)$, $E_2=\SO+b(\f{1-\tau}2)$,
$V_2=\SO+a(\f{1-\tau}2)$.

\end{proposition}

\begin{remark} In Proposition~\ref{a+kb}, we showed that if
an \Hk/ family  exists, then the placement of the ends and
vertical points satisfies the condition $a+kb=k$. In
Proposition~\ref{a+kb2}, we verify that the data we are
constructing in the $gdh$ model does indeed satisfy this
condition. We  are not assuming here that the surfaces \Hk/
exist and we use only the  properties of the $gdh$ and $\f1g
dh$ models to prove the proposition. On the other
hand, the function $G$ in the  proof of
Lemma 3 will turn out to
be the function $g^2$ (the square of the stereographic projection
of the Gauss map on the desired surfaces).

We restate here the immediate consequence of
Proposition~\ref{a+kb}, namely  Corollary~\ref{corabk}, now
valid for our constructed models:
\begin{equation}\label{ineqabk}
a < \f k{k+1} < b,\quad\text{and}\quad\lim_{k\to\infty}b=1.
\end{equation}
\end{remark}
\FigPathsforaPluskb

\begin{proof}
(An alternate proof using theta functions may be found in the Appendix.)
We begin by recalling from the proof of 
Lemma~\ref{symmetriesgdh} the definition of the function 
$$
G=\frac{gdh}{\f1g dh},
$$
a multivalued function with poles at $E_1$ and $V_1$, and
zeros at $E_2$ and $V_2$, of orders given in the following
table:
\begin{equation}
\begin{array}{rlllll}\label{table1}
           & E_1      & V_1  & V_2 & E_2 \\
 \hline
 G       & \infty^{2k}   & \infty^{2}    & 0^2   & 0^{2k}\\
\f{dG}{G}        &\infty          &\infty        &\infty & \infty \\
residue\,\f{dG}{G} & -2k                  & -2           &2      & 2k
\end{array}
\end{equation}
The one-form $dG/G$ has simple poles at precisely these four
points, with  residues given by the values in the array
above. We apply the bilinear  relations to $dG/G$ and $dz$,
obtaining 
\begin{equation} \label{bilinear2}
2\pi i(-2k E_1 -2V_1+2V_2+2kE_2)=\alpha_1 \omega_2-\alpha_2
\omega_1.
\end{equation}
Here the $\alpha_i$ are the periods of $dG/G$ and the
$\omega_i$ are the  periods of $dz$ on the cycles given by
the sides of the rhombus: $\gamma_1=\overline{0,1}$ and
$\gamma_2=\overline{0,\tau}$. (See
Figure~\ref{figure:pathsforapluskb}.) 
Clearly, we have $\omega_1=1$
and $\omega_2=\tau$. Using the expressions for the four poles
of $dG/G$ given in the statement of the proposition, we may
write (\ref{bilinear2}) as
\begin{equation} \label{akb4}
(a+kb)(1-\tau)=\frac{-1}{4\pi i}(\alpha_2-\alpha_1 \tau).
\end{equation}
To prove the proposition we will now show that 
\begin{equation} \label{akb5}
\alpha_1=\alpha_2=-4\pi i k.
\end{equation}

We compute the $\alpha_i$ by replacing the paths $\gamma_i$
with the paths $\tilde{\gamma_i}$ in
Figure~\ref{figure:pathsforapluskb} to which they are 
homotopic. Each $\tilde{\gamma_i}$ consists of the top half
of the vertical  diagonal, three horizontal segments and two
semicircular paths. Since, by Lemma~\ref{symmetriesgdh}, 
we know $G\circ r=1/G$, we have
$$
\alpha_2 = \int_{\tilde{\gamma_2}}\f{dG}G =
\int_{r^{-1}(\tilde{\gamma_2})}r^*(\f{dG}G) =
\int_{r(\tilde{\gamma_2})}\f{d(G\circ r)}{G\circ r} =
\int_{-\tilde{\gamma_1}}\f{d(\f1G)}{\f1G} =
\int_{\tilde{\gamma_1}}\f{dG}G = \alpha_1.
$$
It remains to prove the second equality of (\ref{akb5}). The
contribution to each $\alpha_i$ from the two circular paths
on $\tilde{\gamma_i}$ is clearly equal to $-2\pi i(k+1)$. We
make two  assertions from which the second equality of
(\ref{akb5}) follows immediately. 
\begin{enumerate}
\item The integral of $dG/G$ along the top half of the
vertical diagonal (oriented downward) is $-2\pi i(k-1)$.
\item The integral of $dG/G$ along the union of the three horizontal 
segments of each $\tilde{\gamma_i}$ is equal to zero.
\end{enumerate}

{\em Proof of assertion 1.}  We label the top half of the
vertical diagonal by $\delta$ as in
Figure~\ref{figure:pathsforapluskb}. Using Statement ii)  of
Lemma~\ref{symmetriesgdh} (and its proof), we may assert that
$G=g^2$ is unitary along $\delta$. The integral of $dG/G$ 
along $\delta$ easily computed as follows. Since the form $gdh$
is defined via pullback of $d\z$ from the slit model \tkd/,
the image of $\delta$ under the conformal diffeomorphism $\phi$
from the rhombus model to the slit model is a curve whose
derivative at $\phi(\delta(t))$ is precisely
$gdh(\dot{\delta}(t))$. It follows from Lemma 3, part ii) and the
fact that $\delta$ is fixed by $\mu_h$ that
\begin{equation}
G\circ\delta = \f{gdh(\dot{\delta})} {\f{1}{g}dh(\dot{\delta})} =
\f{gdh(\dot{\delta})}{\overline{gdh(\dot{\delta})}}
\end{equation}
and therefore $G\circ\delta(t)=e^{2i\theta(t)}$, where
$\theta(t)$ is angle determined by the tangent to the image of
$\delta$ under the conformal map $\z$.  Therefore the integral of
$dG/G$ along $\delta$ is $2\pi i$ 
times the total turning of the
normal of the image (under $\phi$) of $\delta$. A curve homotopic
to the image of $\delta$ is drawn in
Figure~\ref{figure:pathsforapluskb} and is also labeled
$\delta$. It is clear that its normal turns through an angle of
$-\pi(k-1)$  as it leaves $P$ and passes through half of the
sewn-in cone.  It is also clear  that the normal does not turn
at all from the time the curve leaves the  cone until it reaches
\O/. This proves the first assertion.

{\em Proof of assertion 2.}  Since $gdh$ is the pullback of $d\z$
under $\phi$ from the rhombic to the slit model, it follows from 
Lemma~3 and our convention that any branch of $gdh$ on the
rhombus takes values of the form $i\la(t)e^{(2\pi ik)n}$ 
on an
interval of the horizontal diagonal bounded by $E_1$ and $E_2$;
here, $\la(t)$ is some nonzero real valued function, $k$ is
determined by our choice of \tkd/, and $n$ is an integer
depending on the branch of $gdh$ we choose.

The same is true for $\f1gdh$ with the same value of $k$, but 
with a
possibly different integer $n$ and 
with a possibly different real-valued function $\la(t)$.
Therefore, the function $G$ on the interval in question takes
values on an open ray in $\BC$, from which it follows that
the integral of $\f{dG}G$ on this interval is equal to the
difference of the values of $\log|G|$ at its endpoints. Since we are free
to choose the two ``semi-circular'' arcs of $\tl\gamma$, so
that $|G|$ is constant along them, it follows that the sum of
the integrals of $\f{dG}G$ along the three horizontal line
segments of $\tilde{\gamma_1}$ is equal to 
$\log(|G(P)|)-\log(|G(\SO)|$. But both $P$ and $\SO$ are on the
vertical diagonal, where $G$ is unitary, which establishes
assertion 2.
\end{proof}

\subsection{The $|dh|$ model and the global formulation of the
vertical period  problem}

We now address the second part of the period problem; that
is, the vertical period problem on \tkd/.  To do this, we
will describe what  amounts to a fourth model for the torus
in question, in terms of which the $dh$-period problem becomes
a question in flat (singular) geometry.

We begin with the  rhombic model described in 4.1. Consider
half of the rhombic torus, produced by identifying the
top-left and bottom-right rectangles as in 
Figure~\ref{figure:rectanglertwo}. 
The  vertical period problem (\ref{vpc3}) is still:
\begin{equation} \label{vpc5}
Re\int_B dh = 0,
\end{equation}
but now we consider, as in Section~3, 
the curve $B$ to be a path in the rectangle 
${\cal 
R}$.
It follows from Lemma~\ref{symmetriesgdh}(iii) that $dh$ is imaginary on the
horizontal sides of the rectangle (as that side corresponds to
the imaginary axis in the $|gdh|$ model) while it is real on
the vertical edges, corresponding to the vertical axis of
$\hksk$. The divisor of $dh$ is given in
Figure~\ref{figure:divisors6}.

\FigDivisorsSix

We define the map $W$ from the rectangle to the complex plane
by
\begin{equation}
W(p)=\int^p_\SO dh,\label{developmentW}
\end{equation}
where $\SO$ is the point on the rectangle ${\cal R}$ corresponding to
the center of the rhombus.

The vertical period problem \eqref{developmentW} is solved
when
\begin{equation}
\ree\int_Bdh = \ree W(P) = 0 \label{vpc6}
\end{equation}
where $P$ is the left vertex of the rectangle.

\subsection{The $(k,d)$ rectangle $\SR$}

Our goal in Section 5 will be to to prove the existence statement of
Theorem 3. This we will do by proving that for every $k>\f12$ there is
a finite positive value $d$ so that \eqref{vpc6} is satisfied, and that
there is a continuous $(k(t),d(t))$ for  which $k(t)$ assumes all
values of $k$ in the open interval $(\f12,\infty)$.  This we accomplish
by an intermediate-value theorem argument with boundary estimates
coming from $d=0$ and $d=\infty$ for fixed $k\in(\f12,\infty]$.

Formally, this requires first adjoining, to the parameter space of
tori $\ST_k(d)$  (with $(k,d)\in(\f12,\infty)\x(0,\infty)$, as defined in
section~\ref{sec:gdhmodelsection}), a set of punctured tori parameterized by
$(k,d)\in[(\f12,\infty]\x\{0,\infty\}]\cup[\{\infty\}\x[0,\infty]]$
and then showing that the height function $h=\int_B dh$ is
continuous on the union, $\SP$, of the parameter space with these
boundary points. These boundary points will be ``degenerate'' in the
following sense. For $1/2<k<\infty$ and $0<d<\infty$, the surface
$\ST_k(d)\setminus\{e_1,e_2,v_1,v_2\}$ is conformally a four-punctured
torus. When we set $k=\infty$, $d=0$, or $d=\infty$; we will be
describing a torus with fewer than four punctures. These punctured
tori are degenerate four-punctured tori because, informally, some of
the punctures have coalesced; formally, these tori are
noded surfaces.

We begin by carefully defining the total parameter space
$\SP$. In Section~5 we will study the height function, $h$ on $\SP$
and prove that it is continuous. In Section 6, we will estimate the
height function in order to apply an
intermediate-value argument to prove the existence part of Theorem~3.

\subsubsection{The definition of $\SP$}

To each $(k,d)\in(\f12,\infty]\x[0,\infty]$, we associate a
punctured Riemann surface. The formal setting for this is the
augmented \tec space $\ov\ST_{1,4}$ of four-times-punctured tori. This
space is the bordification of the \tec space $\ST_{1,4}$ of four-times
punctured tori obtained by attaching, to $\ST_{1,4}$, strata
representing surfaces with nodes; the nodes are obtained by pinching
simple closed curves on a (topological) four-times punctured torus. (We
will stick to our suggestive, but somewhat sloppy terminology of
discussing the ``coalescing'' of distinguished points
$\{p_1,\dots,p_k\}$ on the torus --- here we of course are referring to
the pinching off of a curve surrounding the points $\{p_1,\dots,p_k\}$.
See Remark 2.) In this sense, our definition of the rectangle $\SP$ is
simultaneously a map of the rectangle $(\f12,\infty]\x[0,\infty]$ into
$\ov\ST_{1,4}$.

We begin the definition of $\SP$ on its interior
$(\f12,\infty)\x(0,\infty)$. We have in fact done this in the
Sections~4.1 and 4.2 using the slit model. We remind the reader that in
the  slit model, the parameter $d$ determines the placement of $e_2$,
while $v_2$ is always located at the origin of the $\zeta$-plane. The
pair $(k,d)$ determine the the conformal structure of \tkd/, and a
specifically defined diffeomorphism with a rhombus  representation of
\tkd/ determines the location of the four distinguished points $E_1$,
$E_2$, $V_1$ and $V_2$.  Their positions (see (\ref{ab})) determine the
functions $a=a(k,d)$ and $b=b(k,d)$ satisfying $a+kb=k$ according to
(\ref{abk}) and Proposition~\ref{a+kb}. 

For $d=0$, and $\f12<k<\infty$, we are sewing in a cone point at $v_2$
with $v_2=e_2$. In particular, if we want the structures \tkd/ to be
continuous at $d=0$, we are required to have $a(k,0)=b(k,0)$, and by
\eqref{abk}, we have $a=b=k/(k+1)$. In particular, this will force
$e_1=v_1$. Hence $T_k(0)$ must be a rhombic torus with two
distinguished points. As $k\to\infty$, we have $a=b\to1$, which means
that these two distinguished points should converge at the vertex in the
rhombic model of $T_\infty(0)$. We define $T_\infty(0)$ to be the torus
produced by sewing in a cone of simple exponential type at $v_2$ as in 
Figure~\ref{figure:degenerate1}. 
All the points $V_1$, $V_2$, $E_1$ and $E_2$ then coincide with $a=b=1$.

\FigDegenerateOne

We now define $T_k(\infty)$ for $\f12<k\le\infty$. Here, we are sewing
in a $2\pi k$ cone at infinity. The underlying torus will be defined to
be \T/ (see Section~3) and the points $e_1$ and $e_2$ coincide. This
can only happen when $b=b(k,\infty)=1$ in (\ref{ab}). Hence by
(\ref{abk}), we have $a=0$. We require $E=E_1=E_2$ to be at the vertex
and $V=V_1=V_2$ to be at the center of the rhombus \T/.

Finally, we define $T_\infty(d)$. For $d=0$ or $d=\infty$, we have
already done so. For $0<d<\infty$, we are sewing in a cone of simple
exponential type at the point $di$ in the slit model. This means that
the points labelled $e_1$ and $e_2$ coincide {\em but it does not
necessarily mean that the $v_i$ have to coincide. In fact, as we will
prove in Proposition~\ref{anotzero} below, the points $v_1$ and
$v_2$ are distinct.} Thus, in this case, we will have {\em three}
distinct points: $E=E_1=E_2$, $V_1$ and $V_2$. Here $b=b(\infty,d)=1$
but $a=a(\infty,d)\neq0$, $1$.

\FigKdrectangle

\begin{proposition}\label{anotzero} For $k=\infty$ and any value of
$d\in(0,\infty)$, the points $V_1$ and $V_2$ defined by the
construction of $T_\infty(d)$ do not coincide. In particular, in the
form (\ref{ab}) we have the expression
\begin{align*}
V_1  &= \SO - a\(\f{1-\tau}2\)\\
V_2  &= \SO + a\(\f{1+\tau}2\),
\end{align*}
with $a\neq0$, $1$.
\end{proposition}

\begin{proof} We begin in the slit model. Since we are sewing in a cone of 
simple exponential type, we know that $e_2=di $ and $e_1=\infty$ coincide on 
the torus $T_\infty(d)$. We 
will label that point $e$. When $d>0$, we have $v_2=0\neq e$. Let us now 
look at what this says in the rhombic model of $T_\infty(d)$. First, the 
ends $E_i$ must coincide at a point we label $E$ and that point must be the 
vertex of the rhombus. In particular, $b=1$. Consider the points $V_1$ and 
$V_2$ corresponding to $v_1$ and $v_2$. Because  $v_2\neq e$ we know that 
$V_2\neq E$. In particular, $a\neq1$. We will now show, by contradiction, 
that $a\neq0$. This is equivalent to showing that $V_1\neq V_2$, which is 
in turn equivalent to showing that
 $v_1\neq v_2=0$. 

Suppose $V_1=V_2=\SO$ in the rhombic model. Then the
one-forms $\f1gdh=-r^*(gdh)$ and $gdh$  have the same divisors. This means 
that the cone metrics $|gdh|$ and $|\f1gdh|$ have the same cone points: a 
cone point at $\SO$ with cone angle $6\pi$, and an exponential cone point of 
simple type at $E$, the vertex of the rhombic model. The vertex of the 
rhombus is a fixed point of $r^*,$  
and therefore $r^*_{|E}=-id$. Therefore 
$$\f1gdh_{|E}=-r^*(gdh)_{|E}=gdh_{|E},$$ which implies that the exponential 
cone points of the $|gdh|$ and $|\f1gdh|$ metrics are asymptotically 
isometric.  It then follows from Proposition~\ref{uniquenessinf} of 
Section~\ref{conemetric}
that $|\f1gdh|$ and $|gdh|$ define the same cone metric up to a
constant scale factor.  Since the metrics agree at the fixed points of
$r$, they are in fact equal: $|\f1gdh|=|gdh|$.  

We now use statement ii) of Lemma~\ref{symmetriesgdh}. It says, among other
things, that  $|gdh|$ and $|\f1gdh|$ are interchanged by 
$\overline{\mu}^*_{v}$, where $\mu^*_{v}$ is reflection in the vertical 
diagonal of the rhombus. (Lemma~\ref{symmetriesgdh} 
is presented in a context where it is natural to assume that  that 
$k<\infty$, but its proof does not use this, allowing
us to use its statement ii) here.). Therefore, the involution $\mu^*_{v}$ is
an isometry of the metric $|gdh|$. The fixed points of $\mu^*_{v}$ consist 
entirely of the vertical diagonal, which meets the  horizontal diagonal at 
precisely the points
$\SO$  and $E$.  

We now look at what this implies in the slit model. Here, $gdh$ is the one 
form induced by $d\zeta$. This means that the fixed point set of 
 ${\mu}^*_{v}$, considered now as acting in the slit model,  must consist of 
straight lines segments and rays. Moreover, reflection in these lines and
rays must induce ${\mu}^*_{v}$ in the slit model. Since  ${\mu}^*_{v}$
in the rhombic model leaves the horizontal axis invariant, it follows that 
 ${\mu}^*_{v}$ in the slit model leaves the imaginary $\zeta$-axis invariant.
This is possible only if the fixed point set is either the imaginary 
$\zeta$-axis or the 
real $\zeta$-axis. However, $\mu^*_{v}$ can't be reflection in 
the imaginary axis because reflection in the imaginary axis is $\mu^*_{h}$ 
and $\mu^*_h \neq \mu^*_v$. Similarly,  $\mu^*_{v}$ can't be reflection in 
the real axis of the slit domain because it takes
$e=di$ to  $-di$,  and $-di$ is not a cone point. Therefore, the assumption 
that $V_1=V_2=\SO$ leads to a contradiction.
\end{proof}

\begin{remark} \label{kinfinity}
The definition of the Weierstrass data for $T_\infty(d)$, $d\neq0$,
$\infty$, was forced by the desire to have the data $T_k(d)$ depend
continuously on $(k,d)$ as $k\rightarrow\infty$, something we will prove
to be the case in the next section. In addition to the properties
proved in Proposition~\ref{anotzero}, 
these data have other critical properties.
First, they automatically solve the horizontal period problem. Second,
they define a (possibly multivalued) conformal  minimal immersion at
which there is no period at any end.  Third, they have the properties
of  statement~(iii) of Lemma~\ref{symmetriesgdh}. These assertions are
proved in  the first four paragraphs of the proof of Lemma~\ref{limits}
in  Section~\ref{mainproof}.
\end{remark}

\section{Continuity  and boundary estimates for the height function}

In the previous section we defined for each pair
$(k,d)\in\SP=(\f12,\infty]\x[0,\infty]$, a point \tkd/ in the closure of
the \tec space, $\ST_{1,4}$ of four-times punctured tori. This torus
$T_{k,d}$ carried a pair of one-forms $gdh$ and $\f1gdh$.

\begin{proposition}\label{kdcontinuity} The mapping $(k,d)\to T_k(d)$ is
continuous on $\SP$.
\end{proposition}

We will prove Proposition~\ref{kdcontinuity} in Section 5.1 and its 
subsections. 

\begin{remark} Because the one-forms $gdh$ and $\f1gdh$ are both varying
continuously (on, say, a family of continuously varying models), 
it follows from Proposition~\ref{kdcontinuity} that the
development map $W$ in (\ref{developmentW}) is also continuous as a
function of $(k,d)\in\SP$. 
\end{remark}

\begin{definition} The real-valued function $h:\SP\to\BR$ is
defined by the left-hand side of
\eqref{vpc6} (i.e. $h(k,d) = \ree\int_Bdh = \ree W(P)$) , 
where for $(k,d)\in\SP$ the path $B$ and the one-form $dh$ are taken in
\tkd/.
\end{definition}

It follows from Proposition~\ref{kdcontinuity} that

\begin{proposition}\label{hcontinuity} The function $h:\SP\to\BR$ is
continuous.
\end{proposition}

\begin{proof}[Proof of Proposition~\ref{hcontinuity}] The one-form $dh$
on \tkd/ is determined up to a factor by the conformal structure of
\tkd/ and the divisor of $dh$ on $T_{k,d}$. Of course,
Proposition~\ref{kdcontinuity} asserts that the conformal structures and
divisors of \tkd/ vary continuously in $\SP$. We know from Lemma 3(iii)
that for $(k,d)\in\SP^\circ$, the one-form $dh$ on $T_{k,d}$ is always
imaginary on the horizontal diagonal and real on the vertical diagonal.
Thus if we insist on the forms $dh$ defined on $(k,d)\in\partial\SP$
also satisfying Lemma 3(iii), then we need only  establish the
continuity in $(k,d)$ of the cone metric $|dh|$ to verify continuity of
the forms $dh$: the phase is already continuous. On the other hand, 
the line element 
$|dh|$ can be expressed as $|dh|=\{|gdh|\cd|\f1gdh|\}^{1/2}$, and thus
any scale of $|dh|$ is determined by the scales of $|gdh|$ and
$|\f1gdh|$. We know that these scales are determined and are clearly
continuous by the restriction that the slit $[-1,1]$ on the slit model
always has $|gdh|= |d\zeta|-$length and $|\f1gdh|$-length equal to two.
This  then corresponds to the continuity of the $|dh|$-length of a curve
on \tkd/ for $(k,d)\in\SP ^\circ$. The proposition then follows
once we note that as the underlying conformal structures of
(unpunctured) tori vary continuously over $\SP$, then so do
representatives of the cycle $B$, and hence $h=\int_Bdh$.
\end{proof}

\subsection{The proof of Proposition~\ref{kdcontinuity}}

We prove the continuity in successive steps, each step focusing on
continuity on a particular region on $\SP$. We begin with a preliminary
observation, separating out the important issues  of convergence of
punctured tori from the minor issue of convergence of underlying
(unpunctured) tori.

\begin{lemma} \label{compactnesslemma} For every closed subrectangle 
$[\frac{1}{2}+\epsilon, \infty]\times[0,\infty] \subset \SP$,
with $\epsilon >0$, the corresponding punctured tori form a 
compact set.
\end{lemma}
\FigCompactnessLemma

\begin{proof}
For $(k,d)\in\SP$, we will work with the slit model of the \tkd/.
Consider  two cycles $\g_1$ and $\g_2$ in this model which connect
the lower edge of the slit with its corresponding upper edge via a path
that avoids the positive imaginary axis. They may be chosen to be
symmetric with respect to the imaginary axis. (See 
Figure~\ref{figure:compactnesslemma}.) 
The space
of tori is compact if, and only if, those curve classes have extremal
lengths which are uniformly bounded above and uniformly bounded away
from zero.

Recall from section \ref{extlen} that there are two equivalent
definitions of extremal length, one (the geometric) which lends itself
to upper bounds, and one (the analytic) which lends itself to lower
bounds. Note that in each of the homotopy
classes of $\g_1$ and $\g_2$, because the cone angle
of $e_2$ is bounded away from $\pi$, there is a fixed annulus of positive
modulus which embeds in the torus disjointly from the imaginary axis,
hence independently of $(k,d)$. Using Definition 6 of Section 2.5, this
provides a uniform upper bound for the extremal lengths of those
curves. Next, to prove a lower bound for the extremal lengths of
$\g_1$ and $\g_2$, it is enough using Definition 4 to exhibit a
conformal metric $\rho$ for which
$\ell_\rho(\g_i)^2/\area(\rho)$ is bounded below. But such a
metric is evident: on the box $[-2,2]\x[-1,1]\sbq\BE^2$, set
$\rho\equiv1$ and, on the complement, set $\rho\equiv0$
(here, we implicitly set $\rho\equiv0$ on the conical region
sewn in along the imaginary axis above $(0,di)$). It is clear
that in this metric  $\ell_\rho(\g_i)\geq 1$ while $\area(\rho)=8$,
proving the positive lower bound.
\end{proof}

\subsubsection{Continuity in the interior of $P$}

Continuity on the interior of $\SP$ is almost self evident. Here we take
$\f12<k<\infty$ and $0<d<\infty$, and note, roughly, that small changes in
either $k$ or $d$ in those ranges only
change the $|gdh|$ structure slightly, and so the conformal structure of
the torus with the four distinguished points $\{V_1,V_2,E_1,E_2\}$
changes only slightly.

However, in preparation for the next two subsections, we will discuss 
continuity in the interior carefully 
and in the context of the methods used later.

In the rhombic or $|dz|$-model, there is a flat cone metric
that is isometric to the  pullback to the rhombus of the metric
$|d\zeta|$  metric on in the slit model. Recall that this is the
$|gdh|$ metric on the rhombus, and we note that it is determined up to
a scaling factor by its divisor (by Proposition 3); this factor, say
$C$, is a normalizing
constant. In particular, the location and type of the cone points (the
ends and vertical points) are uniquely determined by the metric $|gdh|$
or even enough of its geometric invariants such as $d_{|gdh|}$ ($V_2$,
$E_2$) or the lengths, $\ell_{|gdh|}(\g_i)$, of the shortest
representative of the free homotopy class of $\g_i$. It  is therefore 
evident that since under a perturbation of $(k,d)$, all the geometric
invariants on the slit model change but slightly, it follows that the 
positions of the ends and vertical points vary only  slightly in the 
rhombic model. Because  $V_1$ and $V_2$ are symmetrically placed with
respect to $\SO$, it follows that $V_1$ also  depends continuously on
$(k,d)$. By Lemma~\ref{compactnesslemma}, the underlying  unpunctured 
torus varies continuously. Therefore, $T_k(d):\SP\to\ov\ST$ is
continuous on the interior of $\SP$.

\subsubsection{Continuity on the top edge away from the right-hand corner.}

By Lemma~\ref{compactnesslemma},  we  know that the underlying rhombic 
tori subconverge, so our attention is focused on the positions of
the points $E_1$, $E_2$, $V_1$, and $V_2$ in the rhombic model.

Recall that the points in the rhombic model labeled by 
$E_i$ (resp. $V_i$) correspond to the points in the slit model labeled
by $e_i$  (resp. $v_i$).

Assuming that $k\le k_0 <\infty$ and we want to prove that $E_1$ and
$E_2$ coalesce as $d \to \infty$.   
Since $E_2-E_1 = 2b = 2b(k,d)$ according to 
\eqref{ab}, the coalescence of the $E_i$ is equivalent to the requirement
that $b(k,d)\to1$ as $d\to \infty $ with  $k\le k_0$.  Again we have,
by \eqref{ab}, that $V_2-V_1= 2a= 2a(k,d)$, and also that $a+kb=k$, a
consequence  of Proposition~\ref{a+kb}. It follows from the boundedness of $k$
that $b(k,d)\to 1$ as $d\to\infty$ if and only if 
$a(k,d)\rightarrow 1$ as $d\to\infty$.
Hence the coalescence of the $E_i$ is equivalent to coalescence of the
$V_i$.

The coalescing of $E_1$ and $E_2$ as $d\to\infty$ for $k\le k_0<\infty$
follows from a simple extremal length argument.  Consider, for small
fixed $\d>0$, an annulus $A(1+\delta,d-\delta)$ 
in the slit model with center at the
origin ($v_2$),  inner radius $1+\d$ and outer radius equal to $d-\d$
(with $d>1$).  This annulus has a large modulus and a core curve 
encircling $e_1$ and $e_2$.  It is clear that the core curve  of 
$A(1+\d,d-\d)$  separates the slit from  the topological disk in the 
complement of $A(1+\d,d-\d)$ that contains $e_1$ and $e_2$.  For $d$
large, The modulus of $A(1+\d,d-\d)$ is essentially $\log d$, and
therefore the  extremal length of a curve encircling $e_1$ and $e_2$
(which we know to be bounded above by $\f1{\log d}$ by Definition 5)
becomes arbitrarily small as $d\to\infty$. Therefore $e_1$ and $e_2$
coalesce on the compact set of rhombi under discussion and it follows 
immediately that the same is true for $E_1$ and $E_2$ in the rhombic
model.

\subsubsection{Uniform estimates near the right-hand edge away from the 
bottom vertex. 
\label{uniformestimates}}

\FigA1D

In order to prove continuity along the right-hand edge, we need to
have  good estimates for the positions of the $E_i$, as functions of $k$
and $d$, as $k$ and/or $d$ go to  $\infty$. To do this we concentrate
on a simply connected neighborhood of  these points. We begin in the
the  slit model. It is convenient for the exposition to do a single
homothety of the slit model, scaling it so that the slit on the real
axis is now along the segment $[-1/2, 1/2]$.  Consider a circle $\SC$
of, say, radius $1$ about the origin. This circle $\SC$ separates the
model into two components: one, say $T$, contains the slit and is
topologically a punctured torus, and the other, say $D$, contains $e_2$
and $e_1$ and is a topological disk. This disk with the metric
$|d\zeta|=|gdh|$ is isometric to a domain, say $\Delta$, which we expect to
be nearly a disk with the metric
$$
ds_\Delta: = C\left|\f{(\xi-\e_2)^{k-1}}{(\xi
-\e_1)^{k+1}}\right| |d\xi|.
$$
(See Appendix B, Lemma~\ref{metricform}.) 
Here, $\e_i$ represents the position in the domain
$\Delta$ of the point $e_i$ in the slit model, and $C$ is a constant
depending upon the geometries of $\Delta$ and of the slit model with the
metric $|gdh|=|d\zeta|$. We translate the domain $\Delta$ so that
$\e_1=0$; then the metric $ds_\Delta$ takes the form
\begin{eqnarray} \label{dsdelta}
ds_\Delta  &= C\left|\f{(\xi-\e_2)^{k-1}}{\xi^{k+1}}\right|
|d\xi|\\
&= C\left| 1 - \f{\e_2}\xi\right|^{k-1}\f{|d\xi|}{|\xi|^2}.
\end{eqnarray}

\FigDeltaModel

We see that in situations with $|\e_2|$ small, we can expect
$\bdy\Delta$ to be a nearly round circle of radius nearly
$1/C$; we will later show that $|\e_2|$ must be small along
any path in $\SP$ along which $(k,d)\to(\infty,\infty)$.

We are interested in determining the asymptotics of $C$ and $\e_2$ as
$k$ or $d$ tends to infinity; this involves matching the geometric
invariants of $(\Delta, ds_\Delta)$ with those of $(D,|gdh|)$. To do
this we insist that $(\Delta,ds_\Delta)$ should be bounded by a 
$ds_\Delta$-round circle of $ds_\Delta$-circumference $2\pi$ and that
the $ds_\Delta$-distance from this circle to $\e_2$ should equal the
$|gdh|$-distance from $\SC$ to $E_2$, which of course is $d-1$. These
are conditions that hold true for the isometric domain $(D,|gdh|)$.

We now introduce two additional normalizations. First, we rotate
$\Delta$ so that $\e_2$ lies on the positive real axis. Second, we
dilate $\Delta$---and consequently alter the constant $C$----so that
$\bdy\Delta\subset D_1(0)$ with $\bdy D\cap\bdy D_1(0)\neq\emptyset$,
i.e. we do a homothety to $\Delta$ so that it just fits within the unit
ball $B_1(0)$, touching $\bdy B_1(0)$ at at least one point.
(See Figure~\ref{figure:deltamodel}.)
\begin{lemma} $C=1+o(1)$, $\e_2=\f1{kd}+o(\f1{dk})$ as
$k\to\infty$, with estimates holding independently of $d$ for
$d$ sufficiently large. \label{candetwoestimates}
\end{lemma}

\begin{proof}
We first prove that $C$ is bounded and that $\e_2\to0$ as $k\to\infty$.
Continuing with our description of the domain $\Delta$, we next make
use of the symmetry of our situation. In the slit model, the disk $D$
is symmetric with respect to 
reflection in the imaginary-axis, and the portion of
the  imaginary axis between $e_2$ and $\SC$ is the unique geodesic
connecting $e_2$ to $\bdy D$. In particular, $\bdy\Delta$ meets the
positive real axis in one point, say $r_+$, and meets the negative real
axis in one point, say $r_-$. As the metric $ds_\Delta$ is invariant
under reflection in the real axis, we conclude that the portion of the
real axis joining $\e_2$ to $r_+\in\bdy\Delta$ is a geodesic. As there
is but one geodesic joining $e_2$ to $\bdy D=\SC$, we conclude that the
isometry taking $D$ to $\Delta$ must take the imaginary axis between
$\SC$ and $e_2$ to the positive real axis between $r_+$ and $\e_2$.

We can now get a crude bound on $C$. Consider the intersection of
$\bdy\Delta$ with $D_{\f{|r_-|}4}(r_-)$, the disk around $r_-$ of
radius  $\f{|r_-|}4$. 
The arc $\bdy\Delta\cap D_{\f{|r_-|}4}(r_-)$ must be properly embedded
in the ball $D_{\f{|r_-|}4}(r_-)$. Also, in that ball,  we have that
$|\xi-\e_2|>|\xi|$ as $\e_2\in\BR_+$. This gives the estimate
\begin{align*}
\f{|\xi-\e_2|^{k-1}}{|\xi|^{k+1}}  &> \f1{|\xi|^2}\\
&> \f{1}{|\f{5}{4}r_-|^2}\\
&= \f{16}{25}\cd\f1{r^2_-}, 
\end{align*}
valid for $\z \in D_{\f{|r_-|}4}(r_-)$.  Because the arc
$\bdy\Delta\cap  D_{|r_-|/4}(r_-)$ passes through the center point $r_-$
of that disk, its $|d\z|$-length is at least $\f{|r_-|}2$.  But thus
implies that it has $ds_\Delta$-length of at least
$\f{|r_-|}2\cd\f{4C}{5|r_-|}$. Since the whole $ds_\Delta$-length of
$\bdy\Delta$ is $2\pi$, we see that $C<\f52\cd2\pi=5\pi$.

We have established that $C$ is finite, and we now turn our attention to
estimating $\e_2$. Consider the arc on the positive real axis connecting 
$\e_2$ and $r_+$. The $ds_\Delta$-length of this arc must be $d-1$, the
distance from $e_2$ to $\SC$ in the $|d\z|$-metric in the slit model. 
Since this also must be the distance in the $ds_\Delta$-metric we may
compute that
\begin{align} \label{distanced-1}
d - 1  &= d(\bdy\Delta, \e_2)\\
&= \int^{r_+}_{\e_2}ds_\Delta \nonumber\\
&= \int^{r_+}_{\e_2}C\biggm| 1 -
\f{\e_2}\xi\biggm|^{k-1}\f{|d\xi|}{|\xi|^2} \nonumber\\
&= C\int^{r_+}_{\e_2}\(1 - \f{\e_2}\xi\)^{k-1}\f{d\xi}{\xi^2}\nonumber\\
&= \f C{\e_2k}(1 - \f{\e_2}{r_+})^k\nonumber.
\end{align}
Now, since $\bdy\Delta\subset\ov{B_1(0)}$, we know that $0<\e_2<r_+<1$.
Hence
$$
d - 1 < \f C{\e_2k}(1 - \e_2)^k < \f C{\e_2k},
$$
which implies that
$$
\e_2  < \f C{k(d-1)} < \f{C'}{kd},
$$
where $C'= \f{Cd_0}{d_0 -1}$ is a constant that depends on the lower
bound, $d_0$ of $d$, which we are assuming to be larger than $1$.
Hence, $\e_2\to0$ as  either $d\to\infty$ or $k\to\infty$.

With the boundedness of $C$ and the decay of $\e_2$ established, we
easily prove the finer estimates of the lemma. We have normalized
$\Delta$ so that $\bdy\Delta$ meets the unit circle (in the $|\d\xi|$
metric) at at least one point. Consider such a point, say $\xi_0$. Then
since $\e_2\to0$, and $C$ is uniformly bounded, the metrics $ds_\Delta$
in (\ref{dsdelta}) subconverge on  a sequence $(k,d)$ and in a fixed
neighborhood of $\xi_0$ to
$$
C_\infty|\f{d\xi}\xi|^2,
$$
where $C_\infty$ is the limit of the constants $C$ in the chosen
subsequence.

Now, we would like to claim that $C_\infty=1$, but whatever it is, the
metric $C_\infty\f{|d\xi|^2}{|\xi|^2}$ will accord constant geodesic
curvature to the arc $\bdy\Delta$ passing through $\xi_0$ only if
$\bdy\Delta$ is the unit circle. But if $\bdy\Delta$ is  the unit
circle, the geodesic curvature will be equal to
one and the length $2\pi$ only if
$C_\infty=1$. This gives the estimate
$$
C = 1 + o(1),
$$
which is the first estimate of the lemma. 
(Here by $o(1)$, we are indicating a quantity which
tends to zero as either $k \to \infty$ or $d \to \infty$.)
In addition this argument 
implies that 
\begin{equation}
r_+ =1 + o(1). \label{rplus}
\end{equation}

To get the finer estimates on $\e_2$, we return to
\eqref{distanced-1}. For $d$ or $k$ large we must have $C=1+o(1)$, and
the estimate (\ref{rplus}). Thus, we can rewrite \eqref{distanced-1} as 
\begin{align*}
d - 1  &= \f{1+o(1)}{\e_2k}\(1 - \f{\e_2k}{k(1+o(1))}\)^k\\
&= \f{1+o(1)}\mu\(1 - \f{\mu(1+o(1))}k\)^k, 
\end{align*}
where we have written $\mu=\e_2k$ in the last term in order to show the
factor $(1-\f\mu k)^k$ to be uniformly bounded for large $k$. Expanding
the right-hand side of the above yields
$$
 d-1= \f1\mu\(1+o(1) - \mu(1+o(1)) + O(\mu^2)\),
$$
from which we conclude that
$$
d = \f1{\e_2k} + \f{o(1)}{\e_2k} + o(1) + O(\e_2k),
$$
or that
$$
\e_2 = \f1{kd} + o\(\f1{d}\)o\(\f1{k}\)
$$
as desired. (The error terms indicates 
a quantity that tends to zero when multiplied by either 
$k$ or $d$ as either $k$ or $d$ tends to infinity.)
This concludes the proof of the lemma.
\end{proof}

\subsubsection{Continuity along the right-hand edge away from the
bottom vertex}
\label{righthandedge}

We now use the estimates of Lemma~\ref{candetwoestimates} in the
previous  subsection to prove continuity along the right-hand edge
when  $d>2$. In fact it will be evident that the proof can be modified
to hold for any  point on that edge of the form $(\infty, d)$, $d>0$. In
particular, it holds at the the top right-hand vertex of $\SP$, the
point $(\infty, \infty)$.

We know from Lemma~\ref{compactnesslemma} that the underlying compact
rhombic tori associated to the slit model determined by $(k,d)$ also 
converge to a non-trivial, non-degenerate rhombus as $kd\to\infty$.
(We will assume  throughout that $d>2$.)

We will  prove first that, as $k\to\infty$, the points $E_1$ and $E_2$ 
coalesce.  This is required because the $(k,d)$ structures we have
defined in Section~4 on the right-hand edge, where $k=\infty$, 
have this property.  Note
that we will allow $(k,d)\to\infty$ along any path where $k\to \infty$ 
and $d>2$. In particular, we allow $d$ to tend to $\infty$ with $k$.

We have constructed in section \ref{uniformestimates}
an isometric model $\Delta$ of a domain $D$ 
on the slit model
containing points $\e_1$ and $\e_2$ corresponding to the points $E_1$
and $E_2$ (respectively). It is  evident from
Lemma~\ref{candetwoestimates} that $\e_2\to\e_1=0$ in this  model, as 
$kd\to\infty$. However, we need to prove that within the corresponding 
rhombic model, the corresponding points $E_1$ and $E_2$  are coalescing
to a point as $k\to\infty$.  To do this, we need  to relate the 
disk $\Delta$ with metric $ds_\Delta$ (which we will refer to as the 
$\Delta$-model) to a domain inside the rhombic model.

We do this by constructing yet another model of the surface, a hybrid of
the slit model and the $\Delta$-model. More precisely, we make a
two-step construction of a metric torus that we will denote by $N$.
First, we apply a quasi-conformal map (with quasiconformal constant
close to one) of the topological disk
$\Delta$ so that it becomes a round disk of radius one. Then, recalling that
the $\Delta$-model is a model for the exterior of the disk of radius one
in the slit model, we sew the interior of the disk of radius one in the
slit model to the disk $\Delta$ (with metric $ds_\Delta)$ along the
common round-circle boundary: here we require the point $r_+$ in
$\Delta$ to glue to the imaginary axis of the slit model. As the circle
is round, this determines the gluing completely. It is clear that if we
were to equip $\Delta$ with the $|d\z|$ metric and perturb it only
slightly quasi-isometrically, then 
the diameter of the resulting torus $N$ is
bounded.

In this construction, we alter the metric and conformal structure of the
surface by the initial quasi-conformal map, and so it is crucial that
this deformation be quite small. However, we have already seen that the
normalizing constant $C$ is $1+o(1)$ as $k\to\infty$ and, by
(\ref{rplus}), $r_+=1+o(1)$.  Thus the boundary $\bdy\Delta$ lies at
radius $1+o(1)$, with geodesic curvature $\kappa=1+o(1)$. Now, as
$\e_2=\f1{kd}+o(\f1{kd})$, we see that if we define our
quasi-conformal map to be the identity on a disk of radius $\f3{kd}$ and
a radial stretch (the stretching dependent on the polar angle)  on the
exterior of that disk, then the map is ($1+o(\f1{\log kd})$)-
quasiconformal.

The metric torus $N$ allows us to estimate the distance separating the
points representing $E_1$ and $E_2$. More precisely, consider the
extremal length of the curve class $\G^*$ consisting of all curves
freely homotopic to the circle $\SC$ in the slit model, a curve that 
encircles $e_1$ and $e_2$. Now the extremal length $\ext_\G(\Om)$ of a
curve system $\G$ in a domain $\Om$ will be dilated or contracted by a
factor of $K$ under a $K$-quasi-conformal map $F:\Om\to\Om'$, i.e.
\begin{equation} \label{distortion}
\f1K\ext_\G(\Om')\le\ext_{\G'}(\Om')\le K\ext_\G(\Om), 
\end{equation}
where $\G'=F(\G)$. Thus the extremal length of the particular class
$\G^*$ will be within a factor of $1+o(\f1{\log kd})$ of the extremal
length $\ext_{\G^*}(N)$  on the hybrid model $N$. But the extremal
length $\ext_{\G^*}(N)$ is easy to compute. From
Lemma~\ref{candetwoestimates}, we have $\e_2=\f1{kd}+o(\f1{dk})$. Using
this  estimate together with the fact that $N$ has finite diameter in
the sense described above we find that
$$
\ext_{\G^*}(N) = \f{2\pi}{\log kd\cd c_N} + o\(\f1{dk}\)
$$
as $k\to\infty$, with $c_N$ a constant that depends on the diameter of
$N$,  which is bounded above and below. (See Ohtsuka \cite{Oht70} Theorems
2.55, 2.80.) Taking into account the distortion of extremal length
caused by the quasi-conformal map of the slit model (which we denote
here by $M_{|d\zeta|}$) to $N$, we  see from the above equation that
\begin{eqnarray}\label{slitmodelestimate}
\ext_{\G^*}(M_{|d\zeta|})  &= \[1+o\(\f1{\log kd}\)\]
\[\f{2\pi}{\log (kd\cd c_N)} + o\(\f1{\log kd}\)\]\nonumber\\
&= \f{2\pi}{\log (kd\cd c_N)} + o\(\f1{\log kd}\)
\end{eqnarray}
as $k\to\infty$. Of course, as the slit model $M_{|d\z|}$ is
conformal {to the  rhombic model, say $M_{|dz|}$, we see 
that
\begin{equation} \label{Esystem}
\ext_{\G^*}(M_{|dz|}) = \f{2\pi}{\log kd\cd c_N} +
o\(\f1{\log kd}\)
\end{equation}
as $k\to\infty$; here, $\G^*$ refers to the system of curves encircling
$E_1$ and $E_2$ in the rhombic model. In particular,  this extremal 
length goes to zero as $k\to\infty$. Thus, $E_1$ and $E_2$ coalesce in
the limit rhombic model.

We now consider the limiting  behavior of the points $V_1$ and $V_2$ in
the  rhombic model. To do so, we will establish a quantitative version
of the  coalescence of $E_1$ and $E_2$ in order to use the relation
(\ref{Esystem}) to estimate  the separation of $V_1$ and $V_2$. This is
now straightforward, for we know that (again see \cite{Oht70}, Theorems
2.55, 2.80) that from (\ref{Esystem}) and (\ref{slitmodelestimate}))
$$
\ext_{\G^*}(M_{|dz|}) = \f{2\pi}{\log\f{c_M}{|E_1-E_2|}} +
o\(\f1{\log|E_1-E_2|}\),
$$
as  $k\to\infty$. Thus, combining the last two estimates for
$\ext_{\G^*}(M_{|dz|})$, we find
\begin{eqnarray} \label{endestimate}
|E_1-E_2|  &= \f1{kd}\cd\f{c_M}{c_N}(1+o(1))\\
&= \f C{kd} + o\(\f1{d}\)o\(\f1{k}\)\nonumber
\end{eqnarray}
as $k\to\infty$.

Now in the notation of (\ref{ab}), $|E_1-E_2|=2(1-b)$ and
$|V_1-V_2|=2a$, and using (\ref{abk}) or Proposition~\ref{a+kb2} (namely
$a+kb=k$) and the above estimate for $|E_1-E_2|$ we find
\begin{align} \label{vestimate}
|V_1-V_2| &= 2a\\ \nonumber
&= 2k(1-b)\\ \nonumber
&= k|E_1-E_2|\\ \nonumber
&= k\cd\(\f C{kd}+o\(\f1{d}\)o\(\f1{k}\)\)\\ \nonumber
&= \f Cd + o\(\f1d\)o\(1\)
\end{align} 
as $k\to\infty$. (Here the term $o(1)$ indicates a 
quantity that tends to zero as $k \to \infty$.)

Thus, not only do $E_1$ and $E_2$ coalesce as $(k\to\infty)$, but so do
$V_1$ and $V_2$, provided $d\to\infty$. This corresponds to the
definition of the twice-punctured torus at $(k,d) =(\infty, \infty)$. If
$d$ is  finite, it is clear from \eqref{vestimate} that the $V_i$  {\em
do not  coalesce} as $(k\to\infty)$. This is consistent with
Proposition~\ref{anotzero} which proves that $V_1$ and $V_2$ must be
distinct on the right-hand edge  when $0<d<\infty$. That same
Proposition shows, assuming that $E_1=E_2$ and  the underlying
conformal stucture of the compact rhombic torus is 
determined---exactly our situation---that the positions of the $V_i$
are  determined. Therefore the limiting data is precisely the data
we  specified on the interior of the right-hand edge.

\subsubsection{Continuity on the bottom edge}
Finally, we concern ourselves with continuity of the map
$\SP\to\ov{T_{1,4}}$ on the set where $k\in(\f12,\infty]$ and $d\geq 0$
is bounded above. The architecture of the argument is identical to that
of Section~\ref{righthandedge}, with only the details (and level of
complication) changing. In particular, we will again find a
homotopically trivial curve $\SC$ in the slit model and describe an
isometric model for the simply connected component $\Delta$ of its
complement. Most of our work will involve a careful study of the
asymptotics of the representatives of $V_1$, $V_2$, $E_1$ and $E_2$ 
inside a suitably modified model $\Delta$. We end the argument by
showing that the asymptotic relationships we found in $\Delta$ have
corresponding statements within the rhombic model.

\FigWeenydisk 

As before, we begin by showing that the convergence is easily
guaranteed when we bound $k$ from above. In particular, consider a
sequence $\{(k_n,d_n)\}$ with $k_n\le k_0<\infty$ and $d_n\to0$; with
no loss in generality, as we can always pass to a subsequence, we
assume that $(k_n,d_n)\to(k_\infty,0)$. We then need to show that $V_2$
and $E_2$ coalesce, as it will then be clear by symmetry that $V_1$ and
$E_1$ coalesce. Thus, it is enough to show that the curve system
consisting of the curves surrounding the segment $i[0,d_n]$ has small
extremal length, or equivalently, that we can embed, into the slit
model, annuli of large modulus whose core curve surrounds $i[0,d_n]$.
Note that since $id_n=e_2$ and $0=v_2$, this will show that $e_2$ and
$v_2$ coalesce. 
This core curve is
slightly awkward to describe: it consists of circles centered
at $0$, $\f12$ and $-\f12$ of radius $2d_n$, with the standard
identifications.  See Figure~\ref{figure:weenydisk}.

With this core curve defined, it is easy to find fat annuli with that
core: consider the annuli of outer radius $\f15$ and inner radius
$2d_n$ about $0$, $\f12$ and $-\f12$. For $\{k_n\}$ bounded, the moduli
of these annuli go to infinity as $d_n\rightarrow 0$. Therefore $V_2$
and $E_2$ (hence also $v_1$ and $e_1$) coalesce. 

Our goal then is to
show the coalescence of $v_1$, $v_2$, $e_1$ and $e_2$ as
$(k_n,d_n)\to(\infty,0)$.

We begin by describing the curve $\SC$ in the slit model. It is crucial
that we can draw a single curve $\SC$ in  the $\z-$plane of the slit
model, namely $\bfC\setminus[-1/2,1/2]$, which surrounds $v_1$, $v_2$,
$e_1$ and $e_2$, under the assumption that 
$(k,d)\in[1,\infty]\x[0,d_0]$.
\FigTheCurveC
 Such a curve is drawn in Figure~\ref{figure:thecurveC}.
(Formally, it is the union of four curves in the plane: circles of
radius $1/10$ centered at $-1/2$ and $1/2$, a semi-circle of radius $1/10$
 in the lower half plane centered at $0$, and a curve that begins at $1/10$,
enters the upper half plane before entering the lower half plane at
$1$, then leaves the lower half plane at $-1$ and finally meets the
$[-1/2,1/2]$ slit again at $-1/10$.) The reader can easily check that
one component, say $D$, of the complement contains $v_1$, $v_2$, $e_1$
and $e_2$, while the other component, say $T$, is a punctured torus.

Let $\Delta$ be a model disk for $D$; we construct a model metric
$ds_\Delta$ on $D$ with singular points $\nu_2$, $\e_1$, and $\e_2$
corresponding to the points $v_2$, $e_1$ and $e_2$, respectively, on
$D$. This metric is isometric to the $|gdh|-$metric (the metric induced
by $|d\z|$ restricted to $D$). We  define $ds_\Delta$ by
\begin{equation}
ds_\Delta = \Bigm|\f{C(\xi-\nu_2)^2(\xi-\e_2)^{k-1}}{(\xi-\e_1)^{k+1}}
d\xi\Bigm|.\label{ds-delta}
\end{equation}
We are permitted to normalize the domain $\Delta$ and, concomitantly,
the metric $ds_\Delta$, by composing our developing map with a
Euclidean isometry; furthermore, as we have already included an unknown
scaling constant $C$ in the form of the metric $ds_\Delta$, we see that
we are also permitted a composition by a homothety. (See 
Lemma~\ref{metricform} in Appendix B for
details.)

With these allowances in mind, we see that we may assume that $\nu_2=0$,
that $\e_1$ lies on the positive real axis and that, after a homothety,
the disk $\Delta$ is contained in the unit disk $B_1(0)$ while
$\bdy\Delta$ meets the unit circle at at least one point.

There is one further normalization that follows not from general facts
about isometries, but from the specific form of the metric $|gdh|$ on
the slit model: the $|gdh|-$metric when considered on the slit model
(i.e the $|d\z|-$metric)  admits a reflection that fixes pointwise
the imaginary axis and $v_2$. (See Lemma~\ref{symmetriesgdh}~(i)and its 
proof.) Thus, we may assume that $(\Delta,ds_\Delta)$ admits an isometry
that fixes $\nu_2=0$, $\e_2$, $\e_1$,  the shortest geodesic from
$\nu_2$ to $\e_2$, and a geodesic from $\e_2$ to $\e_1$.

Now, since the metric (\ref{ds-delta}) is written in terms of
the $|d\xi|$-distances $|\xi-\nu_2|, |\xi-\e_1|$ and $|\xi-\e_2|$,
we see that this isometry must preserve not just
$ds_{\Delta}$ distances from $\nu_2$, $\e_2$ and $\e_1$, but
also $|d\xi|$-distances (in the space $(\Delta,|d\xi|)$) from
$\nu_2$, $\e_2$ and $\e_1$; this is only possible if $\e_2$ is
real. Next we claim that we may assume that $0<\e_2<\e_1$. To see this,
note first  we must have $\e_2>0$, or else there would not exist a path
$\G$ (like the one on the imaginary axis in the slit model) which is a
geodesic from $\e_2$ to $\e_1$ and whose distance from $\nu_2$ increases
along the path. Then next, we see that we must have $\e_1<\e_2$, so that
there exists a geodesic from $\nu_2$ to $\e_1$ whose distance from
$\nu_2$ is an increasing function {\it and} which passes through $\e_2$.

This effect of symmetry simplifies the derivations,
contained in the next  lemma, of the asymptotics
of $C$, $\e_1$ and $\e_2$.  
However, the proof of this lemma is rather long and a bit intricate.
(Happily, though, the conclusion of the proof of 
Proposition~\ref{kdcontinuity} is then nearly immediate.)

\begin{lemma} Consider the metric on the model disk $\Delta$ given by
\label{peskyestimates}
\begin{equation}
ds_\Delta = \Bigm|\f{C\xi^2(\xi-\e_2)^{k-1}}{(\xi-\e_1)^{k+1}}
d\xi\Bigm| \label{ds-deltanormalized},
\end{equation}
which is  normalized as above (in particular, with $\nu_2=0$). 
Then
\begin{enumerate} 
\item[{\rm(i)}] $C\asymp 1$
\item[{\rm(ii)}] as $d\to0$, we have that $\e_2\to0=\nu_2$,
in fact $\e_2=O\(\f1{\log\f1d}\)$.
\item[{\rm(iii)}] as $k\to\infty$, we have that
$|\e_1-\e_2|=O\(\f1k\)$.
\end{enumerate}
\end{lemma}

Here the expression $C\asymp 1$ means that there are constants
$A>0$ and $B< \infty$ so that $A<C<B$. 

\begin{proof} We begin by proving that $C$ is bounded away from zero.
Recall that the distance in the $ds_\Delta$-metric between $\nu_2$ and
$\bdy\Delta$ is equal to the distance in the $|gdh|$ metric between
$v_2$ and $\SC$ in the slit model, which is equal to $d+\f1{10}$ by
construction. Since we are assuming that $d<d_0$, it follows that this
distance between $\nu_2$ and $\bdy\Delta$ is bounded above and below.  
Consider the restriction of $ds_\Delta$ to the  negative real axis.
Along that axis, for any $(k,d)$, we always have $|\xi|<|\xi-\e_1|$, and
$|\xi-\e_2|<|\xi-\e_1|$. Thus, we always have that
$ds_\Delta\bigm|_{\{x<0\}}\le C|d\xi|$. If $C\to0$, then the distance in
the $ds_\Delta$-metric between $\nu_2$ and $\bdy\Delta$  would also go
to $0$. However, the distance in the $ds_\Delta$-metric between $\nu_2$
and $\bdy\Delta$ is equal to the distance in the $|gdh|$  metric
between $v_2$ and $\SC$ in the slit model, and this distance is 
bounded away from $0$. Hence $C$ must be bounded below.

We next claim that there is an upper bound on $C$. This argument is
actually a bit involved, and so we separate it off as a claim.

\paragraph{Claim: $C$ is bounded above.}

To see that there is an upper bound on $C$, we proceed in several steps,
always using that the curve $\SC$ is of fixed length in the slit model,
hence also in $\Delta$. Recall that we have normalized $\Delta$ so that
there is one point, say $\th_0\in\bdy\Delta$, which meets the unit
circle.

To begin, suppose that $\ree\th_0\ge\e_1$, and note that this implies
that $|\th_0|>|\e_2-\th_0|>|\e_1-\th_0|$. Then, at $\th_0$
$$
\biggm|\f{\xi^2(\xi-\e_2)^{k-1}}{(\xi-\e_1)^{k+1}}\biggm| >1.
$$
If the portion of $\SC$ with $\ree\xi\ge\e_1$ has $|d\xi|$-length
bounded from below, then we see that the normalizing constant $C$ is
bounded above, as
\begin{align*}
C_1 &> \int_{\SC\cap\{\ree\xi\ge\e_1\}}ds_\Delta  \\
&=
\int_{\SC\cap\{\ree\xi\ge\e_1\}}
\biggm|\f{\xi^2(\xi-\e_2)^{k-1}}{(\xi-\e_1)^{k+1}}\biggm| |d\xi|\\
&> C\int_{\SC\cap\{\ree\xi\ge\e_1\}}|d\xi|\\
&\ge C\cd C_2
\end{align*}
for $C_1$ a length that is bounded above, and $C_2$ a length that is
bounded below.

In general then, we may assume that $\ree\th_0-\e_1$ has no lower
positive bound, for if it did, we would find a portion of $\SC$
with $\ree\xi > \e_1$ and with length bounded below, and end
with the same contradiction as in the last paragraph.
To further refine the conditions that concern us,
suppose for the moment that there is a lower bound on $\e_1-\e_2$, say
$\e_1-\e_2\ge\e_0>0$. Then in that case the quantity
$|(\xi-\e_2)|/|\xi-\e_1|>1+\e_0'$ for all $\xi$ within a distance of
$(\e_1-\e_2)/4$ of $\th_0$. Thus, in this case
\begin{align*}
C_1 &> \int_{\SC\cap\{|\xi-\th_0|_{ds_\Delta}<\f14(\e_1-\e_2)\}} 
\\
&>
C\cd(2\cd\f14(\e_1-\e_2))\cd(1+\e_0')^{k-1}\cd\int
\biggm|\f{\xi}{\xi-\e_1}\biggm|^2d\xi\\
&\to\infty
\end{align*}
as $k\to\infty$. Thus, we see that it is only possible for $C\to\infty$
if the $\th_0-\e_1$ does not have a positive lower bound and
$\e_1-\e_2\to0$. 

We next dispose of a minor case by noting that it is not possible for
$\th_0\to1$ while $C\to\infty$. This is because these
hypotheses guarantee that there is an arc, say $\SC'$, of $|d\xi|$-length
of at least $(\th_0-\e_1)/4$ along which $|\xi-\e_1|<|\xi-\e_2|$. Thus,
since on this arc, the modulus $|\xi|$ is bounded below, by, say $C_3$,
we see that this arc $\SC'$ has $ds_\Delta$-length of at least
$$
C\cd C^2_3\cd1^{k-1}\cd\f{(\th_0-\e_1)}4\int_{C'}\f{|d\xi|}{|\xi-\e_1|^2}.
$$
Since $|\xi-\e_1|<5/4|\th_0-\e_1|$, the claim follows.

So we may now suppose that $C\to\infty$, $\e_1-\e_2\to0$,
$\ree\th_0-\e_1$ has no positive lower bound and $\th_0$ is bounded away
from 1 (so that $\th_0$ is bounded away from $\e_1$). 
In that case, the method of two paragraphs back provides a lower
bound on the decay of $\e_1-\e_2$: in particular, we claim that
$k(\e_1-\e_2) \to \infty$. To see this, consider a portion
$\SC'$ of $\SC$ within $|d\xi|$-distance of $1/2$ from $\th_0$, which
lies in $\{\ree\xi<\e_1\}$ and has $|d\xi|$-length of $C_3>0$. (It is
sufficient to assume that if $C\to\infty$, such an arc must exist by the
arguments in the previous paragraphs after minor modifications.) Then if
the total length of $\SC$ is $L$, we have
\begin{align*}
L &> \int_{\SC'}\biggm|\f{\xi^2(\xi-\e_2)^{k-1}}{(\xi-\e_1)^{k+1}}\biggm|
|d\xi|\\
&> C\cd
C_3\cd(1/2)^2\inf_{\SC'}\biggm|\f{\xi-\e_2}{\xi-\e_1}\biggm|^{k-1}\cd
\f1{2^2}
\end{align*}
here using the trivial bound $|\xi-\e_1|<2$ which comes from both $\xi$
and $\e_1$ lying in the unit $\xi$-disk. We might as well assume that
$\inf_{\SC'}|\f{\xi-\e_2}{\xi-\e_1}|<1$, or else the 
assumption that $C\to\infty$ immediately
provides a contradiction. Thus we may continue with
\begin{align*}
L &> \f{C_3}{16}C\inf_{\SC'}\biggm|1+\f{\e_1-\e_2}{\xi-\e_1}\biggm|^{k-1}\\
&=
\f{C_3}{16}C\inf_{\SC'}\biggm|1+\f{k(\e_1-\e_2)}{k(\xi-\e_1)}\biggm|^{k-1}
\end{align*}
Now, if $k(\e_1-\e_2)< K_0$ for some fixed $K_0$, then since on
$\SC'$, we have $\ree(\xi - \e_1)<0$ and $|\xi-\e_1|> C_4$ by
hypothesis, we find that 
$$
L>Ce^{-\f{K_0}{C_5}}
$$
where $C_5$ depends on $C_4$, $C_3$ and $|\th_0-\e_1|$: since $\th_0$ may be
assumed bounded away from 1, we see that $|\th_0-\e_1|$ is bounded away
from zero.

From this last expression and our assumption that $C\to\infty$, we see
that it is impossible that $k(\e_1-\e_2)< K_0$ for some fixed $K_0$,
and so $\limsup_kk(\e_1-\e_2)=\infty$.

Finally, consider the point $r_+$ where $\bdy\Delta$ meets the positive
real axis, and consider the curve $\SC^*=\SC\cap B(r_+,(r_+-\e_1)/4)$
which is the intersection of the curve $\SC$ with the ball around $r_+$
of radius $(r_+-\e_1)/4$. Certainly $\SC^*$ has $|d\xi|$-length of at
least $2(1/4(r_+-\e_1))$ and 
we have both that $|\xi-\e_1|<\f{5}{4}|r_+-\e_1|<\f{5}{4}$ and 
$|\xi|^2>(\e_1+ \f{3}{4}(r_+-\e_1))^2> \f{3}{4}(r_+-\e_1))^2$, and so
\begin{align*}
L &> \int_{\SC^*}C\biggm|\f{\xi^2(\xi-\e_2)^{k-1}}{(\xi-\e_1)^{k+1}}\biggm|
|d\xi|\\
&> C\cd1/2(r_+-\e_1)(4/5)^2(\e_1+3/4(r_+-\e_1))^2\inf_{\SC^*}
 \biggm|\f{\xi-\e_2}{\xi-\e_1}\biggm|^{k-1}\\
&> C\cd9/50(r_+-\e_1))^3\inf_{\SC^*}
 \biggm|\f{\xi-\e_2}{\xi-\e_1}\biggm|^{k-1}.
\end{align*}
Thus
\begin{align*}
L  &> \f{9}{50}C(r_+-\e_1)^3\inf_{\SC^*}
\biggm|1+\f{\e_1-\e_2}{\xi-\e_1}\biggm|^{k-1}\\
&> \f{9}{50}C(r_+-\e_1)^3
\biggm|1+\f{k(\e_1-\e_2)}{k(\xi-\e_1)}\biggm|^{k-1}
\end{align*}
Now, choose a subsequence of values of $k$ along our path so that
$k(\e_1-\e_2)$ is monotone and tends to infinity. 
Now, by construction, $\ree(\xi-\e_1) \in (\f{3}{4}(r_+ - \e_1),
(\f{5}{4}r_+ - \e_1))$ (so that in particular, $\ree(\xi-\e_1)>0$).
Since we can restrict to a portion of our path so that
$k(\e_1-\e_2)\ge1$, we find
$$
\biggm|1+\f{k(\e_1-\e_2)}{k(\xi-\e_1)}\biggm|^{k-1} \ge C_6e^{C_7
\ree\f{1}{\xi-\e_1}} \ge C_6e^{C_7 \f{1}{5/4(r_+-\e_1)}} 
$$
Thus 
$$
L>\f{9}{50}C C_6(r_+-\e_1)^3e^{\f{4C_7}{5(r_+-\e_1)}}
$$
Next, choosing a
further subsequence so that $r_+-\e_1$ converges, we see that --
independently of the (bounded) limit of $r_+-\e_1$ -- the right-hand
side of the inequality above must tend to infinity, a contradiction.
This proves the claim. \qed

We have established the statement (i) of the lemma. For statements (ii)
and (iii) to hold, it is necessary that along any path in $\SP$ for
which $(k,d)\to(\infty,0)$, we have that $\e_1$, $\e_2$ and $\nu_2$
all coalesce at $\xi=0$. 

We next assert that $\e_1$ is bounded away from $\bdy\Delta$.
This is because, once we know that $C$ is uniformly bounded on sequences
in our domain, we see that if $d_{|d\xi|}(\e_1,\bdy\Delta)\to0$, then
$\ell_{ds_\Delta}(\bdy\Delta)\to\infty$ as the line element on the limit
point of $\e_1$ would necessarily blow up uniformly: here note that
since $\e_2<\e_1$ and the pole of $|\xi-\e_1|^{(k+1)}$ has order two
greater than that of the zero of $|\xi-\e_2|^{(k-1)}$, this blow-up is
independent of the asymptotics of $\e_2$.

We prove statement (iii) of the lemma, before later
turning our attention to statement (ii). Consider the arc from $V_2$
to $\SC$ in the $|gdh|$-model $M_{|gdh|}$ represented by the segment
$i[-1/10,0]$ in the $|gdh|$ model.  This arc makes an angle of $3\pi$ at
$V_2$ with the arc $i[0,d]$ connecting $V_2$ and $E_2$ in that model
$M_{|gdh|}$.  As the arc $i[0,d]$ in $M_{|gdh|}$ is represented by
$[0,\e_2]$ in $\Delta$, we see that the arc $i[-1/10,0]$ in $M_{|gdh|}$
is represented by a path $\G$ on the negative real axis from $0$
to $\bdy\Delta$ in the model $\Delta$: this is true either by symmetry
or because we check that this path $\G$ is geodesic and makes an angle
of $3\pi$ with the positive real ray emanating from $0\in\Delta$.
Suppose that $\G$ meets $\SC$ at a point $-\d$. Then, computing in the
coordinates $\xi=x+iy$, because of our isometry between $ds_\Delta$ and
$|gdh|$, we know that
\begin{align*}
\f1{10} = d_{|gdh|}(V_2,\SC)  &= \int^0_{-\d}ds_\Delta\\
&= \int^0_{-\d}\biggm|\f\xi{\xi-\e_1}\biggm|^2
\biggm|\f{\xi-\e_2}{\xi-\e_1}\biggm|^{k-1} d\xi\\
&= \int^0_{-\d}\(\f x {x-\e_1}\)^2
\(\f{x-\e_2}{x-\e_1}\)^{k-1}dx\\
&< \d^2\int^0_{-\d}\f1{(x-\e_1)^2}
\(\f{x-\e_2}{x-\e_1}\)^{k-1}dx\\
&= \d^2\int^0_{-\d}\f1{(x-\e_1)^2}\(1 +
\f{\e_1-\e_2}{x-\e_1}\)^{k-1}dx,
\end{align*}
with the inequality coming from $x^2<\d^2$ on the domain of
integration. Now, this last integral we can integrate explicitly to get
\begin{align*}
\f1{10}  &< \d^2\(\f{-1}{\e_1-\e_2}\)\f1k\(1 +
\f{\e_1-\e_2}{x-\e_1}\)^k\bigg|^0_{-\d}\\
&= \f{\d^2}{(\e_1-\e_2)k}\bigg[\(1
-\f{\e_1-\e_2}{\e_1+\d}\)^k
- \(1 - \f{\e_1-\e_2}{\e_1}\)^k\bigg]\\
&< \f{\d^2}{(\e_1-\e_2)k}.
\end{align*}
The final inequality follows from the observation that 
both terms in the difference are positive (using $\e_1>\e_2>0$) and
the first is clearly less than the second. But recall that we
normalized $\Delta$ to be included in the unit disk, so $\d<1$. We
conclude that
$$
\f1{10} < \f1{(\e_1-\e_2)k}
$$
or that 
\begin{equation}
(\e_1-\e_2) < \f{10}k.\label{E1E2decay}
\end{equation}
This establishes statement(iii); note also that the
estimate \eqref{E1E2decay}
is independent of the behavior of $d$. We use this estimate
(\ref{E1E2decay}) on $\e_1-\e_2$ to establish  statement (ii). Here we
proceed similarly, but focus our attention on the path $\ov{V_2E_2}$ in
$M_{|gdh|}$, and correspondingly on the arc $[0,\e_2]$ in 
$\Delta$. From the isometry we constructed, since $\ov{V_2E_2}$ has length $d$
in $M_{|gdh|}$, we have
\begin{align*}
d = d_{|gdh|}(V_2,E_2)  &= \int^{\e_2}_0 ds_\Delta\\
&= \int^{\e_2}_0\(\f x{x-\e_1}\)^2
\(\f{x-\e_2}{x-\e_1}\)^{k-1}dx\\
&> \int^{\e_2}_{\e_2/2}\(\f x{x-\e_1}\)^2
\(\f{x-\e_2}{x-\e_1}\)^{k-1}dx,
\end{align*}
where here we restrict the integral of a positive quantity to
a subinterval $[\e_2/2,\e_2]$. Of course, on that subinterval, we have
that $x>\e_2/2$ and so our estimate becomes
$$
d > \(\f{\e_2}2\)^2\int^{\e_2}_{\e_2/2}\f1{(x-\e_1)^2}
\(1 + \f{\e_1-\e_2}{\xi-\e_1}\)^kdx.
$$
We integrate explicitly as before to find
\begin{align*}
d  &> \f{-\e^2_2}{4k(\e_1-\e_2)}
\(1 + \f{\e_1-\e_2}{x-\e_1}\)^k\bigg|^{\e_2}_{\e_2/2}\\
&= \f{\e^2_2}{4k(\e_1-\e_2)}
\(1 + \f{2(\e_1-\e_2)}{\e_2-2\e_1}\)^k\\
&= \f{\e^2_2}{4k(\e_1-\e_2)}
\(1 - \f{2(\e_1-\e_2)}{\e_1+(\e_1-\e_2)}\)^k.
\end{align*}
Now we invoke the previous estimate (\ref{E1E2decay}) that
$k(\e_1-\e_2)<10$ to obtain
\begin{align*}
d &> \f{\e^2_2}{40}\(1 - \f{2\cd10}{k(\e_1+(\e_1-\e_2))}\)^k\\
&> \f{\e^2_2}{40}\(1 - \f{20}{k\e_2}\)^k
\end{align*}
since the relevant denominator
$k(\e_1+(\e_1-\e_2))>k(\e_2+0)$ in our construction. Now, as
$k\to\infty$, the term $(1-20/k\e_2)^k\to e^{-20/\e_2}$ which
dominates $\e^2_2$. Thus, we conclude that for $k$ large we
have
$$
d > e^{-21/\e_2}
$$
and so
$$
\e_2 < \f{21}{\log 1/d}
$$
for $k$ large, proving statement (ii) of the lemma. This then
concludes the proof of the lemma.  
\end{proof}

We now proceed with the proof of continuity of the mapping 
$(k,d)\rightarrow$\tkd/ at $(\infty, 0)$. The argument is analogous to
the argument in Section~\ref {righthandedge}. It is evident that as
$(k_n,d_n)\to(\infty,0)$, the boundary $\bdy\Delta$ converges to a curve
in the $\xi$-disk which avoids a uniform neighborhood of zero. This is
because the metric $ds_\Delta$ subconverges to the metric
$|d\xi|/|\z|^2$, into which the curve $\SC$ must develop as a curve of
uniformly bounded length. Moreover, by  Lemma~\ref{peskyestimates} for
$(k,d)$ near $(\infty,0)$, we know that $\e_1$ is within a distance of
$O\(\f1{\log\f1d}\)+O\(\f1k\)$ of the origin.  We then create a new
surface $N$ by gluing the complement of $D$ to a
$1+o\(O\(\f1{\log\f1d}\)+O\(\f1k\)\)$ quasi-conformal image of
$\Delta$: here $\SC$ gets joined to $\bdy\Delta$ at corresponding
points. As in subsection~\ref{righthandedge},  we then conclude that
$|E_1-E_2|=O\(\f1k\)$ and $|E_1-V_2|=O\(\f1{\log\f1d}\)$ in the rhombic
model. This concludes the argument for continuity result along the edge
where $d=0$ and $k\in[\f12,\infty]$.

This concludes the proof of Proposition~\ref{kdcontinuity}.

\FigKdrectangletwo

\subsection{Estimates for the height function on the top and bottom of
the rectangle $\SP$ and the solution of the period problem for \Hk/.}
\label{solveperiodproblem}

We begin this section with estimates for the the height function $h$
near  the boundary strata $(1/2,\infty]\x\{0\}$ and
$(1/2,\infty]\x\{\infty\}$. In particular, we prove

\begin{proposition} \label{hnotzero} The function $h$ on the locus
$(\f12,\infty]\x\{0\}$ is positive, i.e. $h\bigm|_{(\f12,\infty]\x\{0\}}>0$.
The function $h$ on the locus $(\f12,\infty]\x\{\infty\}$ is negative,
i.e. $h\bigm|_{(\f12,\infty]\x\{\infty\}}<0$.
\end{proposition}

We postpone the proof of this Propositon briefly, preferring to present
the following Proposition, actually a corollary to
Propositon~\ref{hnotzero}, which will be important in the proof of the
main result (Theorem~\ref{MainThm}) of this paper.

\begin{proposition}\label{ivthypothesis} 
There exists a continuous path $\G\subset (\f12,\infty]\x (0,\infty)$
with  the following properties: for every $k\in (\f12,\infty]$ it
crosses the vertical  segment $\{k\}\x(0,\infty)$;  it passes through
the point $(1,d_1)$ that  corresponds to the Weierstrass data for
\Hone/, the singly periodic genus-one helicoid; every point $(k,d)\in
\G$ corresponds to Weierstrass data for \Hk/$/\sigma_k$ for which 
both the horizontal and vertical period problems are solved. 
\end{proposition}

\begin{proof}({\em of Proposition~\ref{ivthypothesis}}) 
It is immediate from Proposition~\ref{hnotzero} and 
Proposition~\ref{hcontinuity} that there is a neighborhood in $\SP$
of $(\f12,\infty]\x\{0\}$ on which the function $h$ is positive, and  a
neighborhood of $(\f12,\infty]\x\{\infty\}$ on which it is negative.
Therefore, there exists a boundary component of the region $\{h\le0\}$
which separates $(\f12,\infty]\x\{0\}$ from $(\f12,\infty]\x\{\infty\}$
in $(\f12,\infty]\x[0,\infty]$. This in turn implies that there is a
connected  curve, say $\G$, on which $h$ vanishes and with $\G$ having the
following two  properties: $\G$ cuts every vertical line segment
$(k,\infty)$, $k>\f12$ at least once, and $\G$ terminates at a point
on the right-hand edge of the form $(\infty, d_{\infty})$, with
$d_{\infty}\neq 0$. Since $h$ is analytic in $k$ and $d$, we may
assert that $\G$ is piecewise smooth. Each point $(k,d)\in\G$
corresponds to a solution of the vertical period  problem for the
Weierstrass data defined at $(k,d)$, data for which the horizontal 
period problem is solved by construction. Therefore, for every 
$k> \f12$, there exists a solution \Hk/ corresponding to a point on
$\G$, and we can assert the existence of a continuously varying 
(although not necessarily monotonic in $k$) family of such solutions.
From our discussion in Section~3 (in particular, Proposition~\ref{rmk2}) of
the uniqueness of the singly periodic  genus one helicoid, \Hone/, we
may assert  that there is precisely one value of $d$, say $d_1$,  for
which the vertical  period problem for the data at $(1,d)$ has a 
solution. Thus, whatever choice of path $\G$ we make, it passes through
this  distinguished point on the vertical segment where $k=1$.
\end{proof}

\begin{proof}({\em of Proposition \ref{hnotzero}}) Our argument
parallels that of Section~\ref{solvevpc}. We begin with a discussion of
the  ``top stratum" $(\f12,\infty]\x\{\infty\}$. By definition (see
Section~4.4)  the underlying compact torus on this stratum is \T/, and
the ends coincide  and the vertical points coincide: that is,
$E_1=E_2=E$ and  $V_1=V_2=V$. In the notation of (\eqref{ab}),
we have that $a=0$ and $b=1$. Note that this is  independent of choice
of $k$, and hence $h$ is constant on the top stratum. However, in the
proof of Proposition~8, in particular the proof of claim
(\ref{anegative}), we showed that $h$ is negative when $k=1$ and $a=0$.
The situation described there is precisely the one we have here;  the
underlying torus is \T/, the ends coincide and the vertical points 
coincide. This is equivalent to $k=1$ and $d=\infty$. 
This establishes the second statement of the Proposition.

We next consider the ``bottom'' stratum $(\f12,\infty)\x\{0\}$. Here, of
course we have $a=b=k/(k+1)$ and $E_i=V_i$. Notice that for the moment
we  have excluded the point $(\infty, 0)$.  Consider then a neighborhood
of $v_2=e_2$  in the slit model for values of $(k,d)$ converging to
$(k_0, 0)$. For each fixed $(k,d)$, the conformal map $w$ from the slit
model to itself that takes the $|gdh|$ structure to the $|\f1gdh|$
structure takes this neighborhood into a neighborhood containing 
$v_1=e_1$. From Figure~\ref{figure:divisors5} we see that, 
in that neighborhood, the total
curvature of
the cone metric $\{|gdh||w^*\f1gdh|\}^{1/2}$ vanishes.
Thus, using  Proposition~\ref{kdcontinuity}, we can conclude that the
limiting structure $|dh|$ is flat and non-singular, and
a (bounded, nondegenerate) homothety of the rhombic
model (because the $|gdh|$ and $|\f1gdh|$ structures converge as $d\to0$).
The argument then follows exactly as in the proof of
Proposition 8, especially the part that precedes the claim
(\ref{anegative}). In particular, the value of $h$ is given by half the 
length of the vertical axis in the rhombic model. Because the
collection  of underlying closed tori associated to $\SP$ form a
compact set by Lemma~\ref{compactnesslemma}, the
lengths of the vertical axes are uniformly bounded away from zero. Thus,
the values of $h$ along $[k_0,\infty)\x\{0\}$, for any $k_0>\f12$,
are also uniformly bounded away from zero. By continuity of $h$ 
(Proposition~\ref{hcontinuity}), they are positive and bounded away from 
zero on an open neighborhood of the closed  set  $[k_0,\infty]\x\{0\}$
in $\SP$.
\end{proof}

\section{The proof of Theorem~\ref{MainThm}} \label{mainproof} 

We will now prove the main theorem stated in Section~1, namely
\setcounter{theorem}{2}
\begin{theorem}
For every $k>1$, there exists a complete,
$\sigma_k$-invariant, properly embedded minimal surface,
\Hk/, whose quotient by $\sigma _k$ satisfies conditions (4).
As $k\rightarrow\infty$, a limit surface exists and is an {\em
embedded} \HEone/, i.e. a properly embedded minimal surface
satisfying conditions (1).
\end{theorem}

\begin{proof}
Proposition~\ref{ivthypothesis} at the end of
Section~\ref{solveperiodproblem} asserted the existence of a continuous
curve, $\G\in \SP$, of Weierstrass  data for which the vertical and
horizontal period conditions are satisfied. The curve $\G$ crosses all
verticals of the form $\{k\}\times (0, \infty)$, $\f12<k<\infty$.  For
any value of $k\in(\f12,\infty)$ the  Weierstrass data at $(k,d)\in\SP$
were  constructed to guarantee that the conditions (4) hold provided
the vertical period condition $h(k,d)=0$ is satisfied. Therefore,
Proposition~\ref{ivthypothesis} gives us a continuous family of
properly {\em immersed} minimal surfaces \Hk/  satisfying (4). To
complete the proof of the theorem we must show that each surface \Hk/,
$\f12<k<\infty$ in the family defined by $\G$ is embedded  and that the
limit surface---the surface corresponding to the Weierstrass data
$(\infty,d_{\infty})$, is an embedded \HEone/ satsifying the conditions
of (1). 

\begin{lemma}\label{twistedsurfacesembedded} Let $\G \subset\SP$
be the path described above and in Proposition~\ref{ivthypothesis}. Suppose 
$(k,d) \in \G$ with $k\neq\infty$. Then the surface \Hk/ corresponding
to $(k,d)$ is embedded.
\end{lemma}

Lemma~\ref{twistedsurfacesembedded} completes the proof of embededness
of  the continuous family of \Hk/. 

\begin{lemma}\label{limits} At a point $(\infty, d_*)$, $0<d_*<\infty$, 
where $h(\infty, d_*)=0$, the Weierstrass data defines an \HEone/, a
surface satisfying the conditions (1) of Section~1.2:

\begin{enumerate}
\item [(i)] \HEone/ is a properly immersed minimal surface;
\item [(ii)]\HEone/ has genus one and one end, 
that end being asymptotic to the helicoid;
\item [(iii)]\HEone/ contains a single vertical line ({\em the axis})
and a single horizontal line.
\end{enumerate}
\end{lemma}

We will prove both lemmas after completing the proof of the theorem.

Given Lemma~\ref{limits}, all we need to show at this point is that
the limit \HEone/ is embedded. To this end,
consider the limit point of $\G$, a point of the form $(\infty, 
d_{\infty})$, $d_{\infty}\neq 0,\infty$, satisfying the requirement of
Lemma~\ref{limits}, namely $h(\infty, d)=0$.  According to
Lemma~\ref{limits}, its Weierstrass data produce an \HEone/, and
according to  Lemma~\ref{twistedsurfacesembedded}, that \HEone/ is the
limit of embedded  minimal surfaces. By statement (ii) of
Lemma~\ref{limits}, \HEone/ is embedded  outside of a compact set.
Because we may approximate the part of this \HEone/ inside the compact
set by embedded minimal surfaces, the limit  surface has no transverse
intersections in this compact set. Invoking the maximum  principle for
minimal surfaces, we conclude that the only other possibility  is that
the Weierstrass immersion defining this \HEone/ is a multiple covering 
of  its image.  But as pointed out above, this \HEone/ is embedded
outside of  the compact set,  so it must be embedded. (See the proof of
Lemma~\ref{twistedsurfacesembedded} for a similar argument  in  
more detail.) This completes the proof of Theorem~3.
 \end{proof}

\begin{proof}[Proof of Lemma \ref{twistedsurfacesembedded}] From 
Proposition~\ref{ivthypothesis} we know that the path $\G$ must pass
through the point $(1,d_1)$ whose data defines the unique singly
periodic genus-one  helicoid \Hone/.  From  Theorem~\ref{PG1H} we know
that \Hone/ is embedded. Since  each \Hk/ is asymptotic to a helicoid
and therefore embedded outside of a suitably large cylinder about its
vertical axis---and the radius of that  cylinder can be chosen to be a
a continuous function on $\G$---it follows along the lines of arguments that 
are now becoming standard in the subject that all the \Hk/ defined by data 
on $\G$ must be embedded. (There does not seem to be a generally stated
argument in the  literature that applies directly to our case, so we
will give a proof here that all the \Hk/ are embedded.)

Let $\SE$ denote the points of $(k,d)\in\G^{\circ}$ that correspond to
an  embedded \Hk/. (The interior $\G^{\circ}\subset\G$ consists of 
$\G$ minus the right-hand  endpoint where $k=\infty$.) Because
$\G^{\circ}$ contains the point, $(1,d_1)$, representing \Hone/, and
\Hone/ is embedded, the set $\SE$ is nonempty. Since $\G^{\circ}$ 
is connected, showing that $\SE$ is both open and closed will prove the 
lemma.

Consider a point $(k,d)\in\G^{\circ}$. It corresponds to an \Hk/, a
minimal surface that is invariant under the action of $\sigma_{k}$  and 
whose quotient, modulo $\sigma_{k}$, 
has two ends asymptotic to  the ends of the helicoid. 
All such surfaces have a vertical axis and we have
normalized  the family so that the vertical axis of any \Hk/ is the
$x_3$-axis. In particular we may assert that: 
\vspace{.06in}

{\em
(i) Any \Hk/ is embedded outside of some cylinder of sufficiently large 
radius about $x_3$-axis. 
}
\vspace{.06in}

\noindent By Proposition~\ref{kdcontinuity} and the continuity of $\G$
we  may also assert that:
\vspace{.06in}

{\em  (ii) Given any connected subset ${\cal U}\subset\G^{\circ}$ on
which $k$ is bounded away from $\f12$ and $\infty$, there exists an
$R>0$ such that every  minimal surface \Hk/ corresponding to 
$(k,d)\in{\cal U}$ is embedded outside of the vertical cylinder of
radius $R$ around the $x_3$-axis. In fact, they are  uniformly (in $k$)
asymptotic to the ends of the same helicoid.  }

\vspace{.1in}

To see this note that such a subset $\SU$ is compact in $\G^{\circ}$,
and if we have a bound $R$ that holds at $p\in\G^{\circ}$, 
then the bound
$R/2$ holds in a small neighborhood of $\G^{\circ}$ near $p$.

{\em $\SE$ is open in $\G^{\circ}$.} Let $(k_0,d_0)\in \SE$. Choose a  
neighborhood ${\cal U}\subset\G^{\circ}$  of $(k_0,d_0)$ on which $k$ is 
bounded away from $\f12$ and $\infty$ and select $R>0$ according to
statement (ii) above.  The \Hk/ with $k=k_0$ corresponding to
$(k_0,d_0)$ is embedded, so there must be a (possibly smaller)
neighborhood ${\cal U'}$ of $(k_0,d_0)$ in ${\cal U}\subset\G^{\circ}$ 
for which every \Hk/ with $(k,d)\in {\cal U'}$ is  embedded inside the
vertical cylinder of radius $R$. By assertion (ii), they are all 
embedded outside the cylinder of radius $R$. Hence they are all
embedded.

{\em $\SE$ is closed in $\G^{\circ}$.} The maximum principle forbids
a sequence of embedded minimal surfaces from developing a
self-intersection in the limit unless the limit surface is a (branched)
cover of a \ms. Since the limit is an \Hk/, the latter possibility cannot 
occur in our case because of assertion (i) above.

\end{proof}

\begin{proof}[Proof of Lemma \ref{limits}]
We will not use the hypothesis $h(\infty,d)=0$ until later in the
argument.  That is, we begin by discussing the properties of the
interior of the right-hand edge of $\SP$.
At {\em any} point $(\infty,d)$, with $d_*\neq 0, \infty$,
the structure of ${\cal T}_{\infty}(d)$ is defined by the slit model 
for \T/  with a cone, $S_{\infty}$, of simple exponential type  sewn in
along the positive imaginary axis with  vertex at  $v_2=di$. (For the
definition and properties of $S_\infty$ see  Sections~2.3 and 
2.4.) By Proposition~\ref{kdcontinuity} and
Lemma~\ref{compactnesslemma}, 
we know that the underlying 
conformal structure is that of a nondegenerate rhombic torus. From
Section~2.4, we  know that the points  $e_2=di$ and $e_1=\infty$ in
the slit model actually correspond  to the same point on the torus, a
point that we will label $e$. From Proposition~Œ\ref{anotzero}, we know
that the vertical point $v_2$  corresponding to the origin in the slit
model does not coincide with $v_1$. This means that  on the rhombic
model of ${\cal T}_{\infty}(d)$, the point $E$ correponding to $e$ 
must be at the vertex of the rhombus and points $V_i$ corresponding to 
the $v_i$ are distinct and---as usual---symmetrically placed with
respect to the  center of the rhombus. In the notation of that
Proposition~\ref{anotzero} and of (\ref{ab}), we must have $0<a<b=1$.

The Weierstrass data at $(\infty,d)$, with $d\neq 0$, $\infty$, is
defined in order  to solve the horizontal period problem. (Certainly we
have made this clear in  Section~4  for points $(k,d)$ on the interior
of $\SP$. The same argument  works when $k=\infty$. Alternatively, one
can use the continuity of structure in Proposition~\ref{kdcontinuity} to
conclude that the horizontal period  problem is also solved when
$k=\infty$.) Therefore, the Weierstrass data at $(\infty,d)$ defines a
multivalued minimal immersion of a rhombic torus into Euclidean  space
with periods (if any) that are all vertical.

In Lemma~\ref{limits}, (iii) we proved that at a $(k,d)\in \SP^\circ$, 
any  branch of the minimal immersion defined by the Weierstrass data
there has  the property that it maps the horizontal diagonal through
\O/ in the rhombic  model to two horizontal lines in Euclidean space,  and
the vertical diagonal through \O/ in the rhombic model to a vertical
line in Euclidean space. By  Proposition~\ref{kdcontinuity}, the same
is true for the Weierstrass data at $(\infty,d)$, with $d\neq 0$,
$\infty$.

There are only two possible sources of multiple values for the minimal 
immersion defined by the Weierstrass data at $(\infty,d)$, with $d\neq
0$, $\infty$. One, of  course, is the possibility that the vertical
period problem is not solved,  i.e. that $h(\infty,d)\neq 0$. The
second is the existence of a vertical period  at the end. We will show
in this paragraph that the latter does not happen. In  the rhombic
model, the periods at the ends, $E_1, E_2$, for  a $(k,d)$ structure in
the interior of $\SP$ are  vertical and of length $+2\pi k$ at $E_1$ and
$-2\pi k$ at $E_2$.  The  vertical period around a small cycle that
surrounds  $E_1$ and $E_2$ is zero. It follows  from the 
continuity of $(k,d)$ structures given in Proposition~\ref{kdcontinuity}  
and the fact that the ends coalesce as $k\rightarrow\infty$, that
the vertical period of a small cycle around the end $E$ is zero. 

The previous four paragraphs describe the geometric properties of the 
Weierstrass immersion associated to {\em any} $(\infty,d)$, where $d\neq
0$, $\infty$. (See Remark~\ref{kinfinity}.)

 We now use for the first time the
hypothesis that  we are at a point where $h(\infty,d_*)=0$. This
implies immediately that the  Weierstrass data associated to
$(\infty,d_*)$ defines a single-valued,  proper immersion  minimal 
immersion of a torus with one end into Euclidean space. The image must
contain a vertical line and a  horizontal line through the image, say
$p_0$, of \O/ (the center of the  rhombus). In particular, the tangent
plane to the surface at $p_0$ is a vertical plane,  implying that it
can't contain any other horizontal line through $p_0$. If  the image
surface contained any other horizontal or vertical line that did not
pass through the image of \O/ then the surface would be singly 
periodic, an impossibility since the surface must have genus equal to
one. This completes the proof of statements (i), (iii) and the first
part of  statement (ii) of the lemma.

We now turn our attention to the second part of statement (ii). We
will  give three different proofs that the end of this surface is
asymptotic to a  helicoid, the last proof assuming that the surface is
a limit of \mss\ (which is the case of interest).

{\em Proof using analysis of the special end structure.} We consider
the  Weierstrass data at the end, $E_1$, of the minimal surface
associated with the point $(\infty, d_*)$. This data is the limit of 
Weierstrass data at $(k,d)\in\SP^\circ$ with $k\rightarrow\infty$ and
$d\rightarrow d_*$. (Note that we  are {\em not} assuming that
$h(k,d)=0$ for this data. Nor are we assuming  that the limit is
taken over specially chosen values of $(k,d)$.) From Figure~19 we can
read off that $dg/g$ has a simple pole at the ends $E_1$ and $E_2$
where its residues are $k$ and $-k$, respectively. Using the fact,
discussed in the first paragraph of the proof of the lemma,  that 
$E_1$ and $E_2$ coalesce as $k\rightarrow\infty$, as well as 
Proposition~\ref{kdcontinuity}, it follows that the (well-defined) one
form $dg/g$ on the torus defined at $(\infty, d_*)$ has a double pole at
the end $E$. Now, as $V_1\neq V_2$, we see that $dg/g$ has poles at
both $V_1$ and $V_2$ with opposite residues. As $g$ is regular elsewhere
on the surface, we find that $dg/g$ has no residue at the end $E$.

Turning our attention to the one-form $dh$ for values of 
$(k,d)\in\SP^\circ$, a similar argument shows that the simple poles of
$dh$  at $E_1$ and $E_2$ coalesce, as $ (k,d)\rightarrow(\infty, d_*)$,
to a double  pole at the end $E$. Moreover, since $dh$ can have no poles
on the surface,  the residue theorem guarantees that the double pole
at the unique end of the  surface has no residue. Since we have already
established (statement (iii)) that there is one  horizontal and one
vertical line diverging into the end at $E$, we may use  
Proposition~\ref{embeddedend} to conclude that the end at $E$ is
asymptotic to  the end of a helicoid. 

{\em Proof using the asymptotic arguments of Section~\ref{righthandedge} to 
express $dg/g$ and $dh$.} We assume in this proof that there is a
family $\G\subset\SP^\circ$ with end point at $(\infty,d_\infty)$. This
is slightly stronger than the hypotheses of the lemma, but sufficient
for the application. By this method, we avoid the use of
Proposition~\ref{embeddedend}.

Let $(\infty,d_\infty)$ be a point of $\G$ on the line
$\{k=\infty\}\subset\ov\SP$, and let $(k,d_k)$ converge to
$(\infty,d_\infty)$. Simplifying the notation somewhat, let $M_k$
denote the (embedded) \ms\ associated to the point $(k,d_k)$. Each of
the surfaces $M_k$ contains a unique vertical line, and a collection of
horizontal lines; after translating, we can assume that for each $M_k$,
the origin is located at the intersection of the vertical line with a
horizontal line.

As the $|gdh|$ and $|\f1gdh|$ flat structures converge, as
$k\to\infty$, to the (non-degenerate) $|gdh|$ and $|\f1gdh|$ flat
structures associated to the point $(\infty,d_\infty)\in\SP$, we see
that surfaces $M_k$ converge, uniformly on compacta, to a \ms, say
$M_\infty$, whose $|gdh|$ and $|\f1gdh|$ structures are associated to
the point $(\infty,d_\infty)$. In particular, $M_\infty$ is
topologically a torus: note here that it is crucial that the $|gdh|$
model $T_k(d)$ has a slit whose length is constant, hence bounded away
from zero and infinity.

It is elementary at this point to see that the horizontal and vertical
lines are unique; as the surface is invariant with respect to the
group of rotations about parallel lines and so could not be of finite
but non-trivial topology if it were to include a pair of parallel lines.

Finally, we need to show that $M_\infty$ has an end asymptotic to a
helicoid. Naturally, it is enough to check that the Weierstrass data of
the end $E$ (the limit point of the ends $E_1$ and $E_2$ as
$k\to\infty$) has leading terms which agree with those of a helicoid.
For this, it is enough to recall the estimates of 
Lemma~\ref{candetwoestimates}:
$$
|E_1-E_2| = \f C{kd} + o(\f1{kd})
$$
and
$$
C = 1 + o(1),
$$
where $E_i$ referred to the ends in the rhombic model and the estimate
was valid as $k\to\infty$ for points $(k,d)$ with $d$ bounded away from
zero. As in the proof of Proposition~\ref{kdcontinuity}, 
we take $\z$ to be the variable
on the rhombus and we normalize the $E_i$ to be real, then the form
$gdh$ may be represented on $M_{k,d}$ (near, say, $E_1$) 
in this notation as
\begin{align*}
gdh &= -iC\f{(\z-E_2)^{k-1}(\z-V_2)^2}{(\z-E_1)^{k+1}}d\z\\
&=-iC(1 + \f{(E_1-E_2)}{(\z-E_1)})^{k-1}\f{(\z-V_2)^2}{(\z-E_1)^2}d\z
\end{align*}
with the factor of $-i$ coming from the rotation of the $|gdh|$ model
to this model. Substituting in the estimates from  Lemma 5 to this
description and setting $E_1=0$, one finds
\begin{align*}
gdh &= -ic(1 + o(1))(1 + \f{\f{1+o(1)}d+o(1)/d}{k\z})^{k-1}\f{d\z}{\z^2}\\
&\rightarrow-ice^{\f{1}{d}\f{1}{\z}}\f{d\z}{\z^2}
\end{align*}
as $k\to\infty$ when $\z$ is centered at the limit point of the ends
$E_i$ (and $c$ is some constant).
 The same estimates applied to $\f{dg}g$ yields from 
Figure~\ref{figure:divisors5} that
\begin{align*}
\f{dg}g  &= \lb\{\f k{\z-E_1} - \f k{\z-E_2}\rb\}d\z\\
&= \f{k(E_1-E_2)d\z}{(\z-E_1)(\z-E_2)}\\
&\rightarrow\f{1}{d}\f{d\z}{\z^2}
\end{align*}
for $\z$ as above. It is easy to check that these limits provide
Weierstrass data for the helicoid as described in section 2.6.
\end{proof}

\appendix

\section{Appendix:  The Weierstrass data for \Hk/ in terms of
theta functions} Consider a rhombic torus $T_\tau=\bfC /\Lambda$, 
where $\La=\{1,\tau\}$,  on
which we desire to write down Weierstrass data  for the \Hk/ as presented in
Figure~\ref{figure:divisors3}. In the rhombic model described in 
Section~\ref{rhombmodel}, we specified in (\ref{ab}) the location of the 
geometrically 
important points, i.e. the ends and the vertical points.
The center of the torus is located at \O/$=\frac{1+\tau}2$. In
Section~\ref{placementofends}, we saw that the ends $E_i$ and the vertical 
points  $V_i$ could be placed at
\begin{eqnarray} \label{EVab}
 E_1&={\cal O}-b\frac{1-\tau}2, \mbox{\,\,\,} E_2&={\cal 
O}+b\frac{1-\tau}2\\   V_1&={\cal O}-a\frac{1-\tau}2, 
\mbox{\,\,\,}V_2&={\cal O}+ a\frac{1-\tau}2, \nonumber
\end{eqnarray} where $a<b<1$ by assumption.
In particular $\frac12<b<1$.

We will use the theta function
$$\th(z)=\th(z,\tau)=\sum^\infty_{n=-\infty}
e^{\pi(n+\frac12)^2\tau+2\pi i(n+\frac12)(z+\frac12)}$$ to
express $g$ and $dh$. (This theta function is $\th_{1,1}$ in
Mumford  \cite{Mumford83}, pages 17-19.
It has the following properties:
\begin{eqnarray} \label{A1}
\th(0)=&0\text{, a simple zero;}&\notag\\
\th(z + 1)=& \th(z);&\\
\th(z + \tau)=& e^{2\pi i(z+\frac{\tau+1}2)}\th(z);&\notag\\
\th(z)=& \ \ \text{has no poles and no other zeros in a}&
\notag\\
&\text{fundamental domain of }\ T_\tau.\notag
\end{eqnarray}
We may use $\theta(z)$ to write down (perhaps
multivalued) meromorphic functions on $T_\tau$.

\begin{lemma} Let $a_i$ $b_i\in\bfC$, and let $\alpha_i$,
$\beta_i\in\BR$, with $\sum\alpha_i=\sum\beta_i$. Then
$$
f(z) =
\prod^n_{i=1}\frac{\th(z-a_i)^{\alpha_i}}{\th(z-b_i)^{\beta_i}}
$$
has a zero of order $\alpha_i$ at $a_i$, a pole of order $\beta_i$ at
$b_i$,and, modulo $\La=\{1,\tau\}$, no other poles or zeros. Furthermore, 
$f$ and satisfies
\begin{eqnarray}\label{fplus}
&f(z + 1)  & = f(z);\notag\\
&f(z + \tau)   &= e^{2\pi i(\sum\alpha_ia_i-\beta_ib_i)}f(z).
\end{eqnarray}
\end{lemma}

The lemma follows directly from the properties of $\theta$ listed in 
(\ref{A1}).
Using  Figure~\ref{figure:divisors3}, we may express the data 
$g$ and $dh$ in terms of $\theta$ as follows.

\FigDivisorsThree

\begin{align*}
dh   &= e^{it}\frac{\th(z-V_1)}{\th(z-E_1)}\
\frac{\th(z-V_2)}{\th(z-E_2)}dz\\
g(z)   &= \rho\frac{\th(z-V_2)}{\th(z-V_1)}\
\frac{\th(z-(E_2+\tau))^k}{\th(z-E_1)^k}.
\end{align*}
The factor $e^{it}$ in $dh$ is determined not by the divisor, but
rather by 
the requirement
that $dh$ be real on the vertical diagonal. Similarly, the real factor $\rho$
in $g(z)$ is determined by the requirement that $g$ be
unitary on the vertical diagonal. 

The presence of the shift
by $\tau$ (i.e. $E_2+\tau$ instead of $E_2$) in the term in
the numerator of the expression for $g$ is determined by
our desire for $g$ to have the correct transformation behavior, namely:
$$g(z+\tau)=e^{-2\pi ik}g(z).$$
A straightforward computation using the definition of $g$ and (\ref{fplus}) 
gives
 $g(z+\tau)=e^{2\pi iW}g(z)$,
where $W  = V_1 + k(E_1) - V_2 - k(E_2 + \tau)$. Using (\ref{EVab}), we get
\begin{align*}
W &= -2(a + kb)  \frac{1-\tau}2 - k\tau\\
&= -(a + kb) + (a + kb - k)\tau.
\end{align*}

Hence $W=-k$, as required geometrically on a surface invariant under 
the screw motion $\sigma_k$, if and only if $a+kb=k$. 
So in order to obtain the correct transformation behavior,
we must assume that $a+kb=k$.  
This also gives another derivation both of the relationship
$a+kb=k$ proved in Proposition~\ref{a+kb2} and of the 
turning of $g$ along the vertical axis of \Hk/ discussed in 
Remark~\ref{dz}.

\section{Appendix: Existence and uniqueness of flat cone metrics}
We prove here Propositions 3 and 4 from Section~2.3, as well
as a local representation lemma used in the Section~5. Some of these
results extend foundational results of Troyanov \cite{Tr86} in
the case when the cone angles are finite and positive.

\setcounter{proposition}{2}

\begin{proposition} Let $M$ be a compact Riemann surface,
and $\{p_1\dots p_r\dots,p_{r+\ell}\} \subset M$ 
a collection of distinct
points, with $r>0$, $\ell\ge0$. Suppose $\{a_1\dots a_r\}$ is a
collection of real numbers satisfying (\ref{GBe})
\begin{equation}
\sum^r_{j=1}a_j = -(2 - 2genus(M)) + r + 2\ell.\label{res1}
\end{equation}
Then there exists a cone metric on $M$ with finite cone points $p_j$
with cone angles $a_j$ ($j=1,\dots, r$) and exponential cone points $p_k$,
($0\le k\le\ell$) of simple type.
\end{proposition}

\begin{remark} When $\ell=0$, the Proposition is a statement about cone
metrics all of whose cone points are finite and 
\eqref{res1}or (\eqref{GBe}) is the
Gauss-Bonnet condition (\ref{GB}). 
\end{remark}

\begin{proof} We recall a fundamental theorem on the existence of
holomorphic  one-forms with presribed singularities, an eloquent
statment of which can be found in Royden's article \cite{Roy67}:

{\bf Lemma}  {\em Let $M$  be a Riemann surface, $E$ be a closed set
in $M$, $\SO$ an open set containing $E$ and $G$ a bounded
open set with smooth
boundary $\G$ such that $E \subset G$ and $\overline{G}\subset\SO$.  Let
$\Omega$ be an  analytic differential in $\SO \sim E$.  Then there
is an analytic differential $\omega$ in $M\sim E$ such that $\Omega -
\omega$ has an analytic extension  to all of $\SO$ if and only if (the
flux condition) $\int_{\Gamma} \Omega =0$.}

To apply this lemma, we begin by choosing a holomorphic one-form
$\omega_0$ on $M$. (If $M$ has genus zero, a simple modification of the
following  construction will work, and we leave that to the reader.)
Let $\{q_1,\dots,q_s\}$ be the zeros of  $\omega_0$ with the order of
the zero at $q_i$ equal to $b_i$ . Let $E=\{q_1,\dots,q_s\} \cup 
\{p_1,\dots,p_{r+\ell}\}$. Let $\SO$ be the union of disjoint
coordinate  neighborhoods of the points of $E$, with each point
corresponding to $z=0$  in its respective neighborhood. In each
coordinate neighborhood we  specify the analytic differential $\Omega$
on $\SO -E$ as follows:
\begin{itemize}
\item In the punctured neighborhood of $q_i\,$ $1\leq i\leq s$,
$\Omega=-b_i\frac{dz}{z}$;
\item In the punctured neighborhood of $p_j\,$, $1\leq j\leq r$,
$\Omega=(a_j -1)\frac{dz}{z}$;
\item In the  punctured neighborhood of $p_k\,$, $r+1\leq k\leq r+\ell$,  
$\Omega = - (\frac{1}{z^2} +\frac{2}{z})dz$.
\end{itemize}

Letting $G$ be the open set consisting of the union of slightly smaller 
coordinate neighborhoods, we can use the Lemma to assert the existence
of   meromorphic differential $\omega$ on $M$ with principal parts
specified near  the singular set $E$ provided the flux condition
$$
 -\sum_{i=1}^s{b_i} +\sum_{j=1}^r({a_j-1}) -2\ell =0
$$
is satisfied. Since $\omega_0$ is holomorphic,
$-\sum_{i=1}^s{b_i}=2(1-genus(M))$, from which it follows that the flux
condition is precisely  our assumption  (\ref{res1}).  We may also
choose $\omega$ to have purely imaginary periods on $M$. This can be
done by adding to $\omega$ a holomorphic differential specified by
having all its periods on $M$ equal to the negative of the real part of
the periods of $\omega$. The resulting meromorphic one-form has the same
poles and principal parts as $\omega$ and all of its periods imaginary
on $M-E$.

Define $f=e^{\int \omega},$ a multivalued function on $M-E$. Because 
$\omega$ has imaginary periods, the function 
$|f|$ is well-defined on $M-E$, and of
course, is never equal to zero. Let
$$
\eta = f\omega_0.
$$
Note that $\eta$ is regular at the the zeros of $\omega_0$, that
$|\eta|$ is a well-defined metric away from 
$\{p_1\dots p_r\dots,p_{r+\ell}\}$ and that $|\eta|$ has cone points
with  cone angles $a_j$ at the points $p_j$, $1\leq j\leq r$ and
exponential cone points  of simple type at the points $p_k$, $r+1\leq
k\leq r+\ell$.
\end{proof}

\begin{remark}
If one or more of the zeros of the chosen holomorphic one-form
$\,\omega_0$  coincides with a desired cone point, it can be easily
verified that the  construction still produces a cone metric with the
prescribed cone points and cone angles. Simply add the prescribed
residues $b_i$ and $a_j -1$ at those points.
\end{remark}

\begin{proposition} A cone metric on a compact Riemann surface with cone
points with finite cone angles is determined up to scaling by the
location of these cone points and their cone angles. The same result is
true if one or more of the cone points is an exponential cone point of
simple type, provided that these cone points are asymptotically
isometric.
\end{proposition}

\begin{proof} Let $|\mu_1|$ and $|\mu_2|$ be cone
metrics on $M$ with the same cone points and cone angles, and let $M'$
be the Riemann surface $M$ with the  cone points removed. Because
$K=\frac{-\Delta\log|\mu_i|}{|\mu_i|^2} =0$ on $M'$, the function
$\log(|\frac{\mu_2}{\mu_1}|)$ is harmonic on $M'$. In a neighborhood  of
a cone point $p$ with finite cone angle $\alpha$, we may write
$|\mu_1|=|z^{\alpha-1} dz|$ and $|\mu_2|=|w^{\alpha-1} dw|$ with $z$
and $w$ local coordinates, the cone point  $p$ corresponding to $0$ in
both $z$ and $w$ coordinates. If $w=w(z)$ with $w'(0) =c \neq 0$, then
$$
lim_{q\rightarrow p} \frac{|\mu_2|}{|\mu_1|} =c^\alpha.
$$
In particular, $\log(|\frac{\mu_2}{\mu_1}|)$ is bounded in a
neighborhood of $p$. If all the cone points are finite, then
$\log(|\frac{\mu_2}{\mu_1}|)$ is a harmonic function on $M'$ that is
bounded in a neighborhood of each  puncture. Therefore
$\log(|\frac{\mu_2}{\mu_1}|)$ extends to a bounded  harmonic function on
$M$ and hence is the constant function: $\mu_2$=$c\mu_1$ for some
positive constant $c$.

We now do a similar analysis in the neighborhood of a exponential cone
point  of simple type. Suppose there are coordinates in a punctured
neighborhood of $z=0$ in which $|\mu_1|=  |e^{\f1z} \f{dz}{z^2}|$ and
coordinates in a punctured neighborhood of $w=0$ for which $|\mu_2|=
|e^{\f1w} \f{dw}{w^2}|$, the common cone point  corresponding to the
puncture in each disk. Let $w=w(z)$ be a conformal  change 
of coordinates with $w(0)=0$. We may write $|\mu_2|= |e^{\f1w} 
\f{w'dz}{w^2}|$.

Then, for $w=z+o(z^2)$, we have
\begin{align}
\bigm|\f{\mu_2}{\mu_1}\bigm| &= \bigm|e^{\f1w-\f1z}\f{z^2}{w^2}w'\bigm|\\
&= e^{o(1)}(1+o(z))(1+o(z)).
\end{align}
The remainder of the proof goes through as in the case of finite cone
points.

\end{proof}

\begin{remark}
If $w=cz+o(z^2)$ where $c\neq1$, then the ratio above acquires a term
of size $e^{\f1r}$ as $r\to0$, and the proof fails. A counterexample in
the case where the hypothesis doesn't hold is given in Section 2.4.
\end{remark}

We next prove a lemma that was useful in our estimates in Section 5.
While the statement extends to more general settings, we restrict
ourselves here to the situation our estimates require. This restriction
then allows for a simpler proof.

\begin{lemma}\label{metricform} Let $D$ be a topological disk 
equipped with a flat
singular metric $ds_D$ with up to two cone points (say, $p_1$ and
$p_2$) with positive cone angles 
(say, $2\pi\alpha_1$ and $2\pi\alpha_2$)  and up to a single cone point (say
$p_0$) with negative cone angle (say, $2\pi\alpha_0$). 
Next suppose that the cone angle at
$p_1$ exceeds the cone angle at $p_2$, i.e. $\alpha_1 > \alpha_2$). 
Suppose also that a geodesic
$\G$ between $p_2$ and $p_0$ passes through $p_1$, and that this
geodesic bisects the angle found by removing all possible extensions of
$\G$ to a geodesic from $p_2$ to $p_0$. Suppose finally that the sum of
the cone angles vanish. Then $D$ admits a conformal and isometric
development onto a topological disk $\Delta$ equipped with a metric
$ds_\Delta$ of the form
$$
ds_\Delta = C\prod^2_{i=0}|z-q_i|^{\a_i-1},
$$
where $q_i$
is the point in $\Delta$ corresponding to $p_i\in\Delta$. This
representation is unique up to a rigid motion of the plane containing
$\Delta$; a homethety of the plane containing $\Delta$ is an isometry
of $ds_\Delta$ after a corresponding change in the normalizing constant
$C$.
\end{lemma}

\begin{proof} We place $q_0$ at the origin, $q_1$ at the point
$1\in\BC$, and consider $q_2>q_1$ as a variable point on $\BR_+$. The
possible locations of $q_2$ are then a segment. We consider the family
of flat singular metrics on the family of domains $\Delta(q_2)$ with
$$
ds_{\Delta(q_2)} = \prod^2_{i=0}|z-q_i|^{\a_i-1}
$$
where the cone angle at $p_i$ is $2\pi\a_i$, 
as described in the statement of the lemma. The metric
is then determined by the location of $q_2$, as all choices of $q_1$
and $q_2$ on $\BR$ determine a metric whose geodesic connecting $q_2$
with $q_0$ passes through $q_1$; we then seek a metric where
$d_{\Delta(q_2)}(q_1,q_2)=\dist_\Delta(p_1,p_2)$. In fact, in the
present simplified situation of at most three singular points, this
is fairly straightforward, as we merely consider the function
$d(q_2)=d_{\Delta(q_2)}(q_1,q_2)$. Certainly $d(q_2)\to0$ as $q_2\to
q_1$ and $d(q_2)\to\infty$ as $q_2\to\infty$. So by the continuity of
$d$ on the variable $q_2$ and the intermediate value theorem, we
conclude that it is possible to find $q_2\in(1,\infty)$ with
$d(q_2)=\dist_\Delta(p_1,p_2)$ as required.

Note next that when we restrict $ds_{\Delta(q_2)}$ to the complement
$\Delta'(q_2)=\Delta(q_2)-\{q_0,q_1,q_2\}$, we obtain a flat metric. We
may then develop the flat metric $ds_\Delta$ on
$\Delta'=\Delta-\{p_0,p_1,p_2\}$ onto $\Delta'(q_2)$ (with $p_i$ being
sent to $q_i$). This developing map then develops $\bdy D$ onto a
Jordan curve in $\Delta'(q_2)$.

We next turn to uniqueness. Certainly after fixing the positions of
$q_0$ and $q_1$, the position of $q_2$ is determined by the singularity
at $q_2$ and the development of the geodesics from $q_1$ to $q_0$ and
$q_2$ respectively. So suppose we have two metrics, say $ds_1$ and
$ds_2$ on $D$ with the same singularities at $q_0$, $q_1$, and $q_2$.
Then the function $h=\log ds_1/ds_2$ is harmonic, as each $ds_i$ is
flat, and vanishes on $\bdy D$, as the development of $\Delta'$ onto
$\Delta'(q_2)$ was isometric. Thus $h$ vanishes identically, proving
the required uniqueness.
\end{proof}

\newpage

\bibliography{nsf}
\bibliographystyle{plain}

\end{document}